%% file: direction.tex
\documentclass[11pt]{article}
\usepackage{a4wide,amssymb,latexsym,epsfig,amsmath,amsfonts,mathtools,array,graphicx,verbatim}
\usepackage[english]{babel}
\newtheorem{thm}{Theorem}
\newtheorem{lem}[thm]{Lemma}
\newtheorem{prop}[thm]{Proposition}
\newtheorem{cor}[thm]{Corollary}

\newcommand{\sect}[1]{\section{#1}\setcounter{thm}{0}\setcounter{equation}{0}\setcounter{figure}{0}}
{\begin{trivlist}\item[]{\rm\sc{Proof.\ }}%
}{\Qed\end{trivlist}}
\renewcommand{\thesection}{\arabic{section}}

\renewcommand{\theequation}{\thesection.\arabic{equation}}
\newcommand{\be}{\begin{equation}}
\newcommand{\ee}{\end{equation}}
\newcommand{\ba}{\begin{array}}
\newcommand{\ea}{\end{array}}

\newcommand{\noi}{\noindent}
\newcommand{\pa}{\partial}
\newcommand{\na}{\nabla}
\newcommand{\De}{\Delta}
\newcommand{\ti}{\times}

\newcommand{\cpl}{{\mathbb C}}

\newcommand{\eee}{{\bf e}}

\newcommand{\ro}{\varrho}
\newcommand{\roh}{{\widehat{\varrho}}}
\newcommand{\rot}{{\widetilde{\varrho}}}
\newcommand{\vi}{\varphi}
\newcommand{\vii}{\psi}
\newcommand{\ppsi}{{\psi\kern-0.65em \psi}}
\newcommand{\ze}{\zeta}
\newcommand{\zet}{{\tilde \zeta}}
\newcommand{\ga}{\gamma}
\newcommand{\al}{{\alpha}}
\newcommand{\la}{\lambda}
\newcommand{\lat}{{\tilde \lambda}}
\newcommand{\lah}{{\hat \lambda}}
\newcommand{\lala}{{\lambda\kern-0.51em \lambda}}
\newcommand{\Om}{\Omega}

\newcommand{\FF}{\overline{ F'}}

\newcommand{\disp}{\displaystyle}
\newcommand{\ints}{{\mathbb Z}}

\newcommand{\nat}{{\mathbb N}}

\newcommand{\R}{{\rm R\kern-0.1em e\,}}
\newcommand{\RC}{{\rm R\kern-0.1em e_c}}
\newcommand{\RE}{{\rm R\kern-0.1em e_E}}
\newcommand{\Rx}{{\rm R\kern-0.1em e_E^x}}
\newcommand{\Ry}{{\rm R\kern-0.1em e_E^y}}
\newcommand{\qed}{\hfill\hspace*{2mm}\hfill $\Box$}
\newcommand{\Qed}{\nopagebreak\qed}
\newcommand{\real}{{\mathbb R}}
\newcommand{\ra}{\rightarrow}
\newcommand{\rd}{{\rm d}}

\newcommand{\du}{{\delta\! u}}
\newcommand{\dv}{{\delta\! v}}
\newcommand{\aaa}{{\bf a}}

\newcommand{\xx}{{\bf x}}

\newcommand{\nn}{{\bf n}}

\newcommand{\BC}{{\cal B}}
\newcommand{\BtC}{{\widetilde{{\cal B}}}}

\newcommand{\DC}{{\cal D}}
\newcommand{\DCt}{\widetilde{\cal D}}
\newcommand{\FC}{{\cal F}}
\newcommand{\FCt}{\widetilde{\cal F}}
\newcommand{\HC}{{\cal H}}
\newcommand{\LC}{{\cal L}}
\newcommand{\NC}{{\cal N}}
\newcommand{\NCt}{\widetilde{\cal N}}

\newcommand{\BB}{{\bf B}}

\renewcommand{\FF}{{\bf F}}
\newcommand{\Bt}{\widetilde{B}}
\newcommand{\BBt}{\widetilde{{\bf B}}}

\newcommand{\BBh}{\widehat{{\bf B}}}
\newcommand{\DD}{{\bf D}}
\newcommand{\DDh}{{\widehat{\bf D}}}

\newcommand{\Dh}{{\widehat D}}

\newcommand{\Ch}{{\widehat C}}
\newcommand{\phih}{{\widehat \phi}}

\newcommand{\Ct}{{\widetilde C}}
\newcommand{\e}{\epsilon}
\newcommand{\de}{\delta}
\newcommand{\dt}{{\tilde \delta}}
\newcommand{\gt}{{\widetilde \gamma}}
\newcommand{\bt}{{\tilde \beta}}
\newcommand{\rt}{{\tilde r}}
\newcommand{\bti}{{\tilde b}}
\newcommand{\ct}{{\tilde c}}
\newcommand{\ft}{{\tilde f}}

\newcommand{\qt}{{\tilde q}}
\newcommand{\ut}{{\tilde u}}
\newcommand{\vt}{{\tilde v}}
\newcommand{\xt}{{\tilde x}}
\newcommand{\yt}{{\tilde y}}
\newcommand{\zt}{{\tilde z}}
\newcommand{\ztb}{{\overline{\tilde z}}}
\newcommand{\wt}{{\widetilde w}}
\newcommand{\vit}{\widetilde{\psi}}
\newcommand{\chit}{\widetilde{\chi}}
\newcommand{\Rt}{{\widetilde R}}
\newcommand{\At}{{\widetilde A}}
\newcommand{\Qt}{{\widetilde Q}}
\newcommand{\Pht}{{\widetilde \Phi}}
\newcommand{\Pst}{{\widetilde \Psi}}

\newcommand{\zb}{{\overline{z}}}
\newcommand{\fb}{{\overline{f}}}
\newcommand{\hb}{{\overline{h}}}

\newcommand{\Phb}{{\overline{\Phi}}}

\begin{document}
\title{\bf Axisymmetric solutions in the geomagnetic direction problem} 
\author{Ralf Kaiser${}^+$ and Tobias Ramming${}^\#$\\[2ex]
         \small ${}^+$ Fakult\"at f\"ur Mathematik und Physik,
          Universit\"at Bayreuth, D-95440 Bayreuth, Germany\\[1ex]
     \small ${}^\#$ TNG Technology Consulting, D-85774 Unterf\"ohring, Germany\\[1ex]
     \small ralf.kaiser@uni-bayreuth.de, tobias.ramming@gmx.net} 
\date{July, 2019}
\maketitle
%
%
%
%
%
\begin{abstract}
The magnetic field outside the earth is in good approximation a harmonic vector field determined by its values at the earth's surface. The direction problem seeks to determine harmonic vector fields vanishing at infinity and with prescribed direction of the field vector at the surface. In general this type of data does neither guarantee existence nor uniqueness of solutions of the corresponding nonlinear boundary value problem. To determine conditions for existence, to specify the non-uniqueness, and to identify cases of uniqueness is of particular interest when modeling the earth's (or any other celestial body's) magnetic field from these data.

Here we consider the case of {\em axisymmetric} harmonic fields $\BB$ outside the sphere $S^2 \subset \real^3$. We introduce a rotation number $\ro \in \ints$ along a meridian of $S^2$ for any axisymmetric H\"older continuous direction field $\DD \neq 0$ on $S^2$ and, moreover, the (exact) decay order $3 \leq \de \in \ints$ of any axisymmetric harmonic field $\BB$ at infinity. Fixing a meridional plane and in this plane $\ro - \de +1 \geq 0$ points $z_n$ (symmetric with respect to the symmetry axis and with $|z_n| > 1$, $n = 1,\ldots,\rho-\de +1$), we prove the existence of an (up to a positive constant factor) unique harmonic field $\BB$ vanishing at $z_n$ and nowhere else, with decay order $\de$ at infinity, and with direction $\DD$ at $S^2$. The proof is based on the global solution of a nonlinear elliptic boundary value problem, which arises from a complex analytic ansatz for the axisymmetric harmonic field in the meridional plane. The coefficients of the elliptic equation are discontinuous and singular at the symmetry axis, which requires solution techniques that are adapted to this special situation.


\noindent MSC-Classification (2010): 35J65, 86A25.
\end{abstract}
%
%
%
%
%
%
%
%
%
%
%
%
%
\sect{Introduction}
The standard boundary value problems for harmonic vector fields in exterior domains prescribe besides asymptotic conditions at infinity either the normal components or the tangential components of the sought-after field at the boundary. Well-posedness of these problems, solution methods, and corresponding results on existence and uniqueness are well-known (see e.g.\ Martensen 1968). Concerning the {\em geomagnetic} field, however, these types of data are not always available or expensive to provide. In fact,
archaeomagnetic, palaeomagnetic, and even historical magnetic data sets up to the 19th century contain either exclusively information about the direction of the magnetic field vector or provide the directional information more reliably than information about the magnitude of the field vector (for more information about the significance of the direction problem for geomagnetism, we refer to 
Merrill \& McElhinny 1983, Proctor \& Gubbins 1990, Gubbins \& Herrero-Bervera 2007 and references therein). 
In view of the (meanwhile well-established) fact that the geomagnetic field differed in its history drastically from its present form (especially during ``pole reversals"), a general solution theory for large data of the direction problem would be of considerable interest.

When accepting some simplifications such as approximating the earth's surface by the unit sphere $S^2$, neglecting additional sources of the harmonic field in the exterior space $E \subset \real^3$ of the unit ball, and assuming discrete boundary data to be continuously interpolated all over $S^2$, the essence of the direction problem may be formalized as follows: Let $\DD : S^2 \ra \real^3$ be a
nonvanishing continuous vector field (the ``direction field'')
and $\de \in \nat\setminus\{1, 2\}$ (the ``decay order'' of the harmonic field at infinity). Given a direction field and a decay order, the direction problem $P_{\DD}$ asks for all vector fields $\BB \in C^1 (E) \cap C(\overline{E})$ for which a positive continuous function $a: S^2 \ra \real_+$ (the ``amplitude function'') exists such that the conditions
\begin{equation}\label{1.1}
\left.
\begin{array}{cc}
\nabla \times \BB = 0,\ \ \nabla \cdot \BB = 0 & \mbox{ in }E, \\[1ex]
|\BB (\xx)| = O(|\xx|^{-\de}) & \mbox{ for } |\xx|\ra \infty, \\[1ex]
\BB = a \, \DD & \mbox{ on } S^2 
\end{array} \right\}
\end{equation}
are satisfied. The decay order $\de$ is called ``exact" if $|\BB (\xx)| = O(|\xx|^{-\de})$ but not $|\BB (\xx)| = O(|\xx|^{-(\de + 1)})$ for $|\xx| \ra \infty$; it will play a crucial role in the classification of solutions.
Note, however, that usually the exact decay order is not part of the data in the direction problem and $P_\DD$ is formulated with $\de = 3$, which is the lowest possible decay order for magnetic fields vanishing at infinity. 

\begin{figure}
\begin{center}
\includegraphics[width=0.45\textwidth]{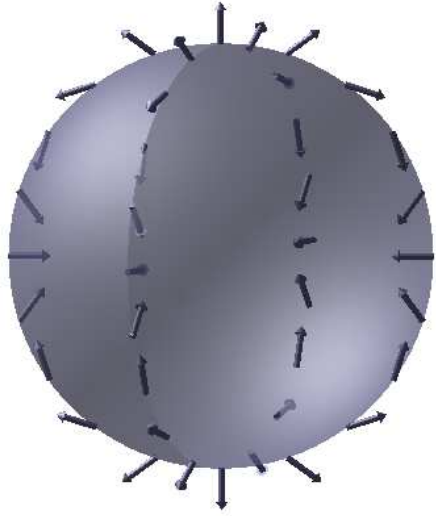}
\caption{Example of an axisymmetric direction field (shown along two meridians on a transparent sphere) with rotation number $\ro = 3$.}
\end{center}
\end{figure}

As is obvious from (\ref{1.1})$_3$ the direction problem is nonlinear in the sense that there is no linear relation between solution $\BB$ and boundary data $\DD$, which means that the usual solution techniques for the above-mentioned standard boundary value problems are not at our disposal and, moreover, that the solution set $S_{\DD}$ for a given direction field $\DD$ is not a linear space. However, the problem can be slightly relaxed so that the enlarged solution set becomes a linear space: dropping the positivity of the amplitude function defines the ``unsigned" direction problem $P^u_{\DD}$ with solution space $L_{\DD} \supset S_{\DD}$.\footnote{Sometimes, if the distinction between ``signed'' and ``unsigned'' direction problem is to be stressed, we use the notation $P_{\DD} =  P_{\DD}^s$ as opposed to $P_{\DD}^u$.} 

The boundary condition (\ref{1.1})$_3$ resembles a well-investigated boundary condition, viz., 
$$
\DD \cdot \BB =b \qquad \qquad \mbox{ on }\ S^2
$$
with given direction field $\DD$ and scalar field $b$ on $S^2$. The ``oblique'' case, i.e., $\DD$ is nowhere tangential to the boundary surface, is well understood and has much in common with the standard boundary value problems (see, e.g., Lieberman 2013). For the ``Poincar\'e problem'', where the obliqueness condition is violated in some part of the boundary, only partial results are so far available and the problem seems not yet to be well understood (see, e.g., Paneah 2000). Only in two dimensions, where the Poincar\'e problem for harmonic fields is known as ``Riemann-Hilbert problem'', there is a close relationship to the (unsigned) 2D-analogue of (\ref{1.1}) (see Kaiser 2010).

So far, for both, the signed and the unsigned, versions of (\ref{1.1}) there are only a few results concerning existence and uniqueness of solutions:
non-uniqueness is known by examples for the (signed and unsigned) direction problem in the axisymmetric case (Proctor \& Gubbins 1990) and in the non-axisymmetric case (Kaiser 2012). For the unsigned direction problem there is, furthermore, an upper bound on the dimension of the solution space $L_{ \DD}$ in terms of
the number $l_{\DD}$ of ``poles'' of the direction field $\DD$ (loci on $S^2$ with vanishing tangential components):
$$
\dim L_{\DD} \leq l_{\DD} -1
$$
(Hulot et al.\ 1997); in general, however, this bound is not sharp (Kaiser 2012). In the axisymmetric situation a better bound has been formulated in terms of rotation number and decay order (Kaiser 2010), and that this bound is sharp will be a corollary of the present work. Concerning the existence of solutions there is a small-data result in the axisymmetric case (Kaiser 2010, see below) 
and some results for special direction fields (Kaiser \& Neudert 2004). The approach in this latter reference is based on $L^2$-expansions in spherical harmonics, a method, which works well if the direction field is itself a single spherical harmonic.   

The present paper provides a complete solution of the {\em axisymmetric} direction problem. The method is inspired by the solution of the two-dimensional version of this problem (Proctor \& Gubbins 1990, Kaiser 2010), which used methods of complex analysis. Axisymmetry leads - in cylindrical coordinates $\rho$, $\ze$ in a meridional plane $\theta = const$ - likewise to a two-dimensional problem, which, although complicated by the coordinate singularity at $\rho = 0$, is amenable to a complex formulation and corresponding ans\"atze. With $(B_\ze, B_\rho)$ representing the nonvanishing components of the axisymmetric harmonic field we make the ansatz
\be \label{1.3}
B_\ze (\ze , \rho) + i\, B_\rho (\ze, \rho) = h(\ze + i \rho) \exp \Big(\frac{1}{2} \big(p(\ze, \rho) + i\, q(\ze, \rho)\big)\Big)
\ee
with the given holomorphic ($\triangleq$ complex analytic) function $h$ representing the zeroes and the asymptotic behaviour of the harmonic field, and the ``correction functions'' $p$ and $q$. Contrary to the two-dimensional case the axisymmetric problem then requires the solution of the following boundary value problem for the (angle-type) variable $q$ in $A_\infty := \{(\ze, \rho) \in \real^2 : \ze^2 + \rho^2 > 1 \} \subset \real^2$:
\begin{equation}\label{1.4}
\left.
\begin{array}{rll}
- \Delta\, q\!\!\!& = \disp - \pa_\ze \Big( \frac{1}{\rho} \cos (q - \Psi)\Big) + \pa_\rho \Big(\frac{1}{\rho} \sin (q - \Psi) \Big)  & \mbox{ in } A_\infty , \\[1ex]
q \!\!\!& = \, \phi & \mbox{ on } S^1 , 
\end{array} \right\}
\end{equation}
where the bounded but discontinuous (angle-type) function $\Psi$ is derived from $h$. $\Psi$, the boundary function $\phi$, and the (weak) solution $q$ are assumed to be antisymmetric with respect to the variable $\rho$. The key problem with (\ref{1.4}), which prevents the application of more or less standard solution methods, is clearly the singular coefficients on the right-hand side of (\ref{1.4})$_1$. In particular the second term on the right-hand side behaves near the symmetry axis $\{\rho =0\}$ like a second order derivative, which has to be controlled by the left-hand side. In (Kaiser 2010) this problem could be bypassed by a suitable embedding of the problem in $\real^5$ that eliminated the coordinate singularity; however, at the price of a then unbounded nonlinearity. Accordingly only a small data result could be achieved (by the Banach fixed point principle). Unfortunately, the smallness assumptions were depending on constants whose numerical values are unknown; so, the compatibiliy of these assumptions with the physical data of the problem (direction field and decay order) remained an open question. 

In this paper the two-dimensional, singular, but nonlinearly bounded problem (\ref{1.4}) is directly attacked via Schauder's fixed point principle. Without smallness assumptions the structure of the nonlinear terms must now more carefully be utilized. Key to success is a suitable choice of auxiliarily introduced parameters at the linear as well as at the nonlinear level of the solution procedure. At the linear level weighted Hardy-type inequalities of the form  
\be \label{1.5}
\int_{\real \ti \real_+} \Big|\frac{f}{\rho}\Big|^2\, \frac{\rd \ze \rd \rho}{\rho^\ga} \leq \Big( \frac{2}{1 + \ga}\Big)^2 \int_{\real \ti \real_+} |\pa_\rho f|^2\, \frac{\rd \ze \rd \rho}{\rho^\ga} \; , \quad f \in C_0^\infty (\real\ti \real_+)\, ,\; \ga \neq -1
\ee
allow the control of the right-hand side in (\ref{1.4})$_1$ by the left-hand side. The solution of a linearized version of (\ref{1.4}) then proceeds by a weighted Lax-Milgram-type solution criterion, whose applicability depends on two coercivity-type constants which in turn depend on the weight. The optimal constants are determined by min-max problems depending on $\ga$ and further parameters. Assisted by numerical computations we derive rigorous lower bounds on these constants with the result that the criterion works but only in a small ``window'' of $\ga$-values. At the nonlinear level, the crucial point is to devise a space large enough to comprise the solutions of the linearized problem but not too large in order not to loose control of the nonlinear terms. Here we make use of a weighted $L^p$-space and, again, only the subtle balance between $p$ and the weight makes Schauder's principle work.

In a first step, this program is carried out in the bounded regions $A_n := \{(\ze, \rho) \in \real^2 : 1 < \ze^2+ \rho^2 <n^2 \}$,  $n \in \nat\setminus \{1\}$ with artificial conditions at the exterior boundaries. The corresponding sequence of solutions turns out to be uniformly bounded; so, in a second step, by means of a ``diagonal argument'', one then obtains a solution of (\ref{1.4}) in the unbounded region $A_\infty$. Finally, the function $p$ is determined from $q$ up to a constant $p_0$, and by substitution into (\ref{1.3}) one obtains the complete set of solutions of the (signed as well as unsigned) direction problem (\ref{1.1}). 
\sect{Reformulation of the problem, results, and sketch of proof} 
This section provides the mathematical framework for our treatment of the direction problem, we review results from (Kaiser 2010) as far as they are relevant for the following, reduce problem (\ref{1.1}) to (\ref{1.4}), and present our results in the theorems 2.1 -- 2.5. Finally, we give an outline of the proofs.

Axisymmetric harmonic vector fields $\BB$, expressed in cylindrical coordinates $(\rho,\theta, \zeta)$, have just two nontrivial components $B_\rho$ and $B_\ze$, depending on $\rho$ and $\ze$, which satisfy the system  
\be \label{2.1}
\left. \ba{c} \pa_\ze B_\rho - \pa_\rho B_\ze = 0 , \\[1ex] 
\disp \pa_\ze B_\ze + \pa _\rho B_\rho + \frac{1}{\rho}\, B_\rho = 0 . 
\ea \right\}
\ee
So far, $\BB$ is defined on the half-plane $H = \{(\ze, \rho) \in \real^2: \rho >0\}$ bounded by the symmetry axis $\{ \rho = 0\}$. It is convenient, however, to extend the domain of definition to $\real^2$ by (anti)symmetric continuation: 
\be \label{2.2}
\left.
\ba{l}
B_\ze (\ze, -\rho) := B_\ze (\ze, \rho) , \\[1ex]
B_\rho (\ze, -\rho) := - B_\rho (\ze, \rho) 
\ea \right\}
\qquad (\ze, \rho) \in H .
\ee
Note that (\ref{2.1})$_2$ implies (if defined) $B_\rho (\ze, 0) = 0$. 
On $\real^2$ also polar coordinates $(r, \vi) \in (0, \infty) \ti (-\pi , \pi]$ with basis vectors $\eee_r$ and $\eee_\vi$ are useful.
They are related to $(\ze, \rho)$ by 
\be \label{2.3}
\ze = r \cos \vi\, , \quad \rho = r \sin \vi
\ee
(see Fig.\ 2.1). 
\begin{figure}
\begin{center}
\includegraphics[width=0.45\textwidth]{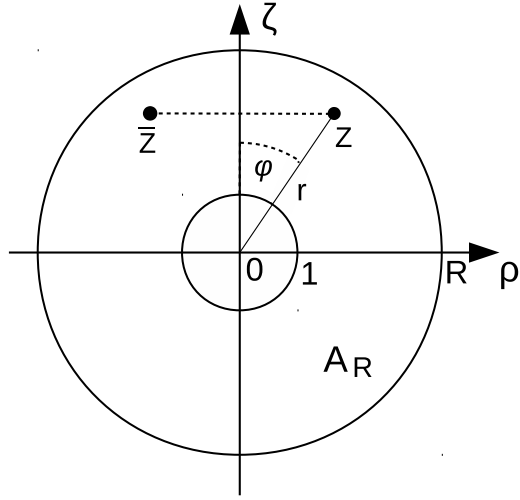}
\caption{Various coordinates in the meridional cross section $A_R$.}
\end{center}
\end{figure}
In these coordinates condition (\ref{2.2}) takes the form
\be \label{2.4}
\left.
\ba{l}
\Bt_\ze (r, -\vi) = \Bt_\ze (r, \vi) , \\[1ex]
\Bt_\rho (r, -\vi) = - \Bt_\rho (r, \vi), 
\ea \right\}
\ee
where $\BB(r \cos \vi, r \sin \vi) =: \BBt (r, \vi)$. Harmonic fields are associated with harmonic potentials $\Upsilon (r, \vi)$ by 
$$
\BBt = \na \Upsilon = \pa_r \Upsilon\, \eee_r + \frac{1}{r}\, \pa_\vi \Upsilon\, \eee_\vi ,
$$
and these potentials have well-known series representations in $A_\infty$, which is a cross-section of the exterior space $E$ through the symmetry axis:
$$
\Upsilon (r, \vi) = \sum_{n= \dt -2}^\infty \frac{c_n}{r^{n+1}}\, P_n (\cos \vi)\, , \quad c_{\dt -2} \neq 0\, , \quad \dt \in \nat\setminus \{1,2\}
$$
(see, e.g., Folland 1995, p.\ 144). Here $P_n$ is the Legendre polynomial of order $n$ and $\dt \in \nat \setminus \{1,2\}$ is the exact decay order of the associated magnetic field:\footnote{ {\em Exact} decay orders are often denoted by 
$\dt$ as opposed to $\de$ for decay orders in general.}
\be \label{2.5}
\ba{rl}
\BBt(r,\vi) = \BBt^{\dt}(r, \vi) &\!\!\! \disp = \na \Upsilon (r, \vi) = - \sum_{n= \dt -2}^\infty \frac{c_n}{r^{n+2}}\, \DD^{n} (\vi) \\[3ex]
& \!\!\!\disp = -\, \frac{c_{\dt-2}}{r^\dt}\, \DD^{\dt-2} (\vi) + O(r^{-(\dt+1)}) \qquad \mbox{for } r \ra \infty ,
\ea
\ee 
where
$$
\DD^{n} (\vi) := (n+1) P_n (\cos \vi)\, \eee_r + P_n'(\cos \vi) \sin \vi \, \eee_\vi
$$
is the exterior axisymmetric $2^n$-pole field restricted to the unit circle. These series are converging uniformly and absolutely for any $r>1$.

Axisymmetric harmonic fields have only a finite number of isolated zeroes with finite negative indices (``x-points'') in $A_\infty$. Let these zeroes be contained in the annulus $A_R := \{(\ze, \rho) \in \real: 1 < \ze^2 + \rho^2 < R^2 \} \subset A_\infty$, then its number $\nu (\BB, A_R)$ can be computed by 
\be \label{2.6}
\nu = \nu(\BB, A_R) = \ro - \roh\, ,
\ee
where $\ro$ and $\roh$ mean the rotation numbers of $\BB$ along the circles $S_1$ and $S_R$ of radii $1$ and $R$, respectively.\footnote{ We skip henceforth the upper index $1$ at $S$ indicating the dimension of the ``sphere''.} The rotation number $\ro \in \ints$ counts the number of turns, the field vector makes when circling once around $S_1$. In $\nu$ zeroes are counted as often as indicated by its index.  Equation (\ref{2.6}) is clearly an analogue of the argument principle in complex analysis and it has likewise some invariance properties with respect to continuous deformations (see section 9); in particular, $\roh$ is constant in the limit $R\ra \infty$, which yields for $A_\infty$ the relation  
\be \label{2.7}
\nu = \nu(\BB, A_\infty ) = \ro - \dt + 1\, .
\ee

As in the two-dimensional case a complex formulation of the (signed) axisymmetric direction problem is promising as it allows to view the direction of $\BB$ as the argument of a complex function $f$. With the identifications
\be \label{2.7a}
\Re z := \ze\, , \qquad \Im z := \rho\, , \qquad \Re f := B_\ze \, ,\qquad \Im f := - B_\rho\, ,
\ee
eqs.\ (\ref{2.1}) then take the form 
\be \label{2.8}
\pa_\zb f - \frac{1}{2}\, \frac{f - \fb}{z - \zb} = 0
\ee
and the symmetry condition (\ref{2.2}) amounts to 
\be \label{2.9}
f(z, \zb) = \fb (\zb , z)\, . 
\ee
We made there use of the Wirtinger derivatives $\pa_z := \frac{1}{2} (\pa_\ze - i \pa_\rho)$ and $\pa_\zb := \frac{1}{2} (\pa_\ze + i \pa_\rho)$
acting on functions $f: \cpl \ti \cpl \ra \cpl$, $(z, \zb) \mapsto f(z, \zb)$. Other than in two dimensions, where harmonic fields satisfy $\pa_\zb f = 0$ (i.e., $f$ is an analytic function), we are here left with the solution of the singular equation (\ref{2.8}).
 
Direction fields $\DD : S_1 \ra \real^2$, $\vi \mapsto \DD (\vi)$ are called symmetric, if
\be \label{2.9a}
D_\ze (\vi) = D_\ze(-\vi)\, ,\quad D_\rho (\vi) = - D_\rho(-\vi)\, ,
\ee
and fields $\BB$ and functions $f$ are called symmetric, if they satisfy (\ref{2.2}) and (\ref{2.9}), respectively.
Symmetry implies obviously $D_\rho (0) = D_\rho (\pi) = 0$ and hence, as $\DD \neq 0$, $D_\ze (0) > 0$ or $D_\ze (0) < 0$. 
The condition $D_\ze (0) > 0$ is no restriction in problem (\ref{1.1}) and is henceforth considered as implied by symmetry.
The {\em axisymmetric} direction problem $P_\DD = P_\DD (A_\infty)$ then reads:

\hspace{1ex}

\noi\textbf{Problem $P_{\DD}(A_\infty)$}: {\it Let $\DD \in C(S_1, \real^2)$ be a symmetric direction field and $\delta
\in \nat \setminus \{1, 2\}$. Determine all symmetric solutions $\BB \in C^1 (A_\infty) \cap C(\overline{A_\infty})$ of {\rm (\ref{2.1})}
with decay order $\de$ and boundary condition
\be \label{2.10}
\exists\, a \in C (S_1, \real_+) : \BB \big|_{S_1} = a\, \DD\, . 
\ee
}

\noi Equivalent is the following complex formulation:

\hspace{1ex}

\noi\textbf{Problem $P^c_{\DD}(A_\infty)$}: {\it Let $\DD \in C(S_1, \real^2)$ be a symmetric direction field and $\delta \in \nat \setminus \{1, 2\}$. Determine all symmetric solutions $f \in C^1 ( A_\infty) \cap C(\overline{A_\infty})$ of {\rm (\ref{2.8})}
with decay order $\de$ and boundary condition
\be \label{2.11}
\arg f \big|_{S_1} =  \arg D_c\, , 
\ee
where $D_c := D_\ze - i D_\rho$.\footnote{Note that the definitions of $D_c$ here and in ref.\ (Kaiser 2010) differ by complex conjugation.}
}

\hspace{1ex}

\noi The ambiguity in the arg-function is removed by the condition that $\vi \mapsto \arg D_c$ is continuous on $(-\pi, \pi)$. We then have $\arg D_c (0) = 0$ and $\arg D_c (\pm \pi) = \mp \ro\, \pi$.

A bounded version of the problem reads as follows:

\hspace{1ex}

\noi\textbf{Problem $P^c_{\DD, \DDh} (A_R)$}: {\it Let $\DD \in C(S_1, \real^2)$ and $\widehat{\DD} \in C(S_R, \real^2)$ be symmetric direction fields. Determine all symmetric solutions $f \in C^1 ( A_R) \cap C(\overline{A_R})$ of {\rm (\ref{2.8})} with boundary conditions
\be \label{2.12}
\arg f \big|_{S_1} =  \arg D_c\, , \qquad   \arg f \big|_{S_R} =  \arg \widehat{D}_c\, , 
\ee
where $D_c := D_\ze - i D_\rho$ and $\Dh_c := \Dh_\ze - i \Dh_\rho$ .
}

\hspace{1ex}

The basic idea to solve the direction problem is to extend the direction field from the boundary to the entire annulus and to replace the boundary value problem for the harmonic field by one for the direction field. The direction field, however, is in general multivalued and not well-defined at the zeroes of the harmonic field. This suggests for the harmonic field an ansatz with a given field describing the zeroes and by (\ref{2.6}) the rotation numbers at the boundaries, and a further zero-free field with well-defined directions on the entire annulus. For the bounded complex problem $P_{\DD, \DDh}^c (A_R)$ such an ansatz is  
\be \label{2.13}
f(z,\zb) = h(z)\, e^{g(z,\zb)} 
\ee
with the analytic function
\be \label{2.14}
h(z) := \prod_{n=1}^{\ro - \roh} (z - z_n)\, z^{-\ro}
\ee
and the exponential function $e^g$ with well-defined argument function $\Im g$ on $A_R$. Note that for given direction fields $\DD$ and $\DDh$ and hence rotation numbers $\ro$ and $\roh$, the number of zeroes in $A_R$ is fixed by (\ref{2.6}), so that the ansatz (\ref{2.13}) does not restrict the solution set of $P_{\DD, \DDh}^c (A_R)$. However, there is a (preliminary) restriction: in order that $h(z)$ is a symmetric function, the set $S$ of zeroes must be symmetric, i.e.\ $z \in S$ implies $\zb \in S$. If zeroes on the symmetry axis are not allowed (as it will be the case in the subsequent solution procedure), $\ro - \roh$ must be an even number (a restriction that is lifted in section 9). 

When inserting (\ref{2.13}) into (\ref{2.8}) one obtains
\be \label{2.15}
-\pa_\zb\, g = \frac{1}{2}\, \frac{1}{z-\zb} \Big( \frac{\hb}{h}\, e^{- 2i \Im g} - 1 \Big) ,
\ee
which by further differentiation can be reduced to a semilinear elliptic equation in the variable $2 \Im g =: q$  alone:
\be \label{2.16}
-\pa_z \pa_\zb\, q = \Im \Bigg\{\pa_z \Big[ \frac{1}{z-\zb} \Big( \frac{\hb}{h}\, e^{- i q} - 1 \Big)\Big]\Bigg\} .
\ee
As $|\hb/h| = 1$, an angle-type variable $\Psi$ may be introduced by $\hb/h =: e^{i \Psi}$, where $\Psi$ is a bounded but discontinuous function on $A_R$ (see appendix A). Using the real variables $(\ze, \rho)$ in $A_R$, eq.\ (\ref{2.16}) then takes the real form
\be \label{2.17}
\Delta q - \pa_\ze \Big( \frac{1}{\rho} \big(\cos (q -\Psi) -1\big) \Big) + \pa_\rho \Big(\frac{1}{\rho}  \sin (q - \Psi) \Big) = 0\, ,
\ee
where $\Delta := \pa_\ze^2 + \pa_\rho^2$. Boundary values for $q$ on $S_1$ arise from (\ref{2.12})--(\ref{2.14}):
$$
\arg D_c = \arg f|_{S_1} = \bigg(\sum_{n=1}^{\ro - \roh} \arg (z-z_n) - \ro \arg z + \Im g\bigg)\bigg|_{S_1} \, ,
$$
which give rise to the definition\footnote{The tilde denotes again dependence on polar coordinates: $\qt (r, \vi) := q(r \cos \vi, r \sin \vi)$.}
\be \label{2.18}
\phi (\vi) := \qt (1, \vi) = 2\bigg(\ro\, \vi - \sum_{n=1}^{\ro -\roh} \arg \big(e^{i \vi} - z_n\big) + \arg D_c (\vi) \bigg) ,
\ee
and analogously on $S_R$:
\be \label{2.19}
\phih (\vi) := \qt (R, \vi) = 2 \bigg(\ro\, \vi - \sum_{n=1}^{\ro -\roh} \arg \big(R\,e^{i \vi} - z_n\big) + \arg \Dh_c (\vi) \bigg) .
\ee
Note that by construction $\phi (\pi) = \phi (-\pi) = \phih (\pi) = \phih (-\pi) = 0$, i.e.\ $D_c \in C(S_1)$ implies $\phi \in C(S_1)$, and analogously for $\phih$. 

Once $q$ is determined, the real part $p:= 2 \Re g$ of $g$ is given by the other half of eq.\ (\ref{2.15}) up to a constant $p_0$:
\be \label{2.20}
-\pa_\zb\, p = i \pa_\zb\, q + \frac{1}{z-\zb} \big( e^{- i (q - \Psi)} - 1 \big) .
\ee

It is useful, in particular for a weak formulation of the problem, to transform to zero boundary conditions. To this end let $\Phi$ be a harmonic interpolation of the boundary functions, i.e.\ a solution of the (standard) boundary value problem 
\be \label{2.21}
\left.\ba{c}
\Delta \Phi = 0 \qquad \mbox{ in } A_R\, ,\\[1ex]
\Phi\big|_{S_1} = \phi\, ,\qquad \Phi\big|_{S_R} = \phih\, ,
\ea
\right\}
\ee
and define $u:= q - \Phi$, $\Om:= \Psi - \Phi$. In these variables eq.\ (\ref{2.17}) takes the form  
\be \label{2.22}
\Delta u + \pa_\ze \Big( \frac{1}{\rho} \big(1 - \cos (u - \Om) \big) \Big) + \pa_\rho \Big(\frac{1}{\rho}  \sin (u - \Om) \Big) = 0
\ee
or
\be \label{2.23}
\na \cdot \Big( \na u + \aaa[u - \Om]\, \frac{u}{\rho}\Big) = \na \cdot \Big( \aaa [u - \Om]\, \frac{\Om}{\rho} \Big)
\ee
with 
\be \label{2.24}
a_\ze [x] := \frac{1 - \cos x}{x}\, ,\qquad a_\rho [x] := \frac{\sin x}{x}
\ee
and 
\be \label{2.25}
u\big|_{\pa A_R} = 0\, .
\ee

Symmetry of $f$ and $h$ implies 
\be \label{2.26}
p (\ze , -\rho) = p(\ze, \rho) \, ,\qquad q(\ze, -\rho) = - q (\ze, \rho)\, ,
\ee
and by (\ref{2.18}), (\ref{2.19}), $\phi (-\vi) = - \phi (\vi)$, $\phih(-\vi) = - \phih(\vi)$, and hence (see appendix B), 
\be \label{2.27}
\Phi (\ze , -\rho) = -\Phi(\ze, \rho) \, ,\qquad \Psi(\ze, -\rho) = - \Psi (\ze, \rho)\, ,
\ee
which implies, finally,
\be \label{2.28}
u (\ze , -\rho) = -u(\ze, \rho) \, ,\qquad \Om(\ze, -\rho) = - \Om (\ze, \rho)\, .
\ee
The boundary value problem (\ref{2.23})--(\ref{2.25}), and (\ref{2.28}) for $u$ with given $\Om$ is henceforth called {\bf problem} $P_\Om (A_R)$.

The divergence structure of eq.\ (\ref{2.23}) and its non-smooth coefficients suggest a weak formulation of $P_\Om (A_R)$: a function $u \in H_{0,as}^1 (A_R)$ is called weak solution of $P_\Om (A_R)$, if $u$ satisfies
\be \label{2.29}
\int_{A_R} \na u \cdot \na \vii\, \rd \ze \rd \rho + \int_{A_R} \frac{u}{\rho}\, \aaa[u - \Om] \cdot \na \vii\, \rd \ze \rd \rho =
\int_{A_R} \frac{\Om}{\rho}\, \aaa[u - \Om] \cdot \na \vii\, \rd \ze \rd \rho
\ee
for all antisymmetric testfunctions $\vii \in C_{0, as}^\infty (A_R)$, where
$$
C_{0,as}^\infty (A_R) = \big\{\vii \in C_0^\infty (A_R) : \vii (\ze , -\rho) = -\vii(\ze, \rho)\big\}
$$
and 
$$
H_{0, as}^1 (A_R) = clos \big(C_{0, as}^\infty (A_R)\, , \, \|\na  \cdot\|_{L^2 (A_R)} \big) .
$$
Note that in the bounded domain $A_R$ by (\ref{1.5}) $\|\na\cdot\|_{L^2}$ is equivalent to the usual $H^1$-norm $\| \cdot\|^2_{H^1} = \| \cdot\|^2_{L^2} + \| \na \cdot \|^2_{L^2}$; $C_0^\infty (A_R)$ denotes as usual the space of infinitely differentiable functions compactly supported in $A_R$. In the unbounded case a function $u \in H_{loc , as}^1 (A_\infty)$ with $trace\, u|_{S_1} = 0$ is called a weak solution of {\bf problem} $P_\Om (A_\infty)$, if $u$ satisfies (\ref{2.29}) (with $A_R$ replaced by $A_\infty$) for every test function $\vii \in C_{0, as}^\infty (A_\infty)$. The following theorems assert the existence of unique weak solutions in bounded annuli $A_R$ and in the exterior plane $A_\infty$. 
\begin{thm}
 Let $R > 1$ and $\Om \in L^\infty (A_R)$ with bound
 \be \label{2.30}
 |{\widetilde \Om} (\cdot , \vi)| \leq K\, |\sin \vi| \qquad \mbox{in }\, A_R
 \ee
 for some constant $K > 0$, then problem $P_\Om (A_R)$ has a unique weak (in the sense of eq.\ (\ref{2.29})) solution $u \in H^1_{0,as} (A_R)$
 with bound
 \be \label{2.31}
 \int_{A_R} |\na u|^2\,  \rho^{-\beta }\, \rd \ze \rd \rho \leq C\, ,
 \ee
 where $\beta = 1/5$ and $C$ is some constant that depends on $K$, but does not depend on $R$. 
\end{thm}
\begin{thm}
 Let $\Om \in L^\infty (A_\infty)$ with bound
 \be \label{2.32}
 |{\widetilde \Om} (\cdot , \vi)| \leq K\, |\sin \vi| \qquad \mbox{in }\, A_\infty
 \ee
 for some constant $K > 0$, then problem $P_\Om (A_\infty)$ has a unique weak solution $u \in H^1_{loc,as} (A_\infty)$
 with $trace\, u\big|_{S_1} = 0$ and the bound
 \be \label{2.33}
 \int_{A_\infty} |\na u|^2\,  \rho^{-\beta }\, \rd \ze \rd \rho \leq C
 \ee
 with $\beta = 1/5$ and some constant $C > 0$.
\end{thm}
Based on these results the subsequent theorems (2.3)--(2.5) give answers to the direction problems $P^c_{\DD, \DDh} (A_R)$, $P_\DD (A_R)$, and $P^u_\DD (A_\infty)$, respectively. To this end recall that for any continuous, symmetric direction field the $\rho$-component vanishes by (\ref{2.9a})$_2$ at the symmetry axis. The more precise condition
\be \label{2.34}
D_\rho (\vi) = O (\vi) \quad \mbox{for }\, \vi \ra 0\, ,\qquad D_\rho (\vi) = O (\pi - \vi) \quad  \mbox{for }\, \vi \nearrow \pi
\ee
or, equivalently with $D_c = D_\ze - i D_\rho$,
\be \label{2.35}
\arg D_c (\vi) = O (\vi) \quad \mbox{for }\, \vi \ra 0\, , \qquad \arg D_c (\vi) + \ro \, \pi = O (\pi - \vi) \quad \mbox{for }\, \vi \nearrow \pi
\ee
turns out to be more appropriate in the following. Moreover, H\"older continuity of $\DD$ will help to establish continuity of the solution up to the boundary.
\begin{thm}
 Let $\DD$ and $\DDh$ be H\"older continuous, symmetric direction fields with rotation numbers $\ro$ and $\roh \in \nat \setminus \{1\}$, respectively, $\ro - \roh \geq 0$ and even, and satisfying condition (\ref{2.34}) at the symmetry axis. Let, furthermore, $\{z_1, \ldots,z_{\ro - \roh} \}$ be a symmetric set of points in $A_R$. Then, problem $P^c_{\DD, \DDh} (A_R)$ has a unique solution $f = B_\ze - i B_\rho \triangleq \BB$ vanishing at  $z_1, \ldots,z_{\ro - \roh}$ and nowhere else.
\end{thm}

\begin{thm}
 Let $\DD$ be a H\"older continuous, symmetric direction field with rotation number $\ro \in \nat \setminus \{1\}$ and satisfying condition (\ref{2.34}) at the symmetry axis. Let, furthermore, $\dt \in \nat \setminus \{1,2\}$, $\dt \leq \ro +1$, and $\{z_1, \ldots,z_{\ro - \dt +1} \}$ be a symmetric set of points in $A_\infty$. Then, problem $P^c_{\DD} (A_\infty)$ has a unique solution $f \triangleq \BB$ with exact decay order $\dt$ vanishing at  $z_1, \ldots,z_{\ro - \dt + 1}$ and nowhere else.
\end{thm}
When the zeroes of a solution do not matter, the total set of solutions (with unspecified zeroes) of the signed direction problem can best be \ ``counted'' by means of the unsigned problem. Let $\DD$ be a H\"older continuous, symmetric direction field with rotation number $\ro$ satisfying the axis-condition  (\ref{2.34}) and let $S_\DD^\dt$ and $C_\DD^\de$ be the solution sets of the signed direction problem $P_\DD^s (A_\infty)$ with exact decay order $\dt$
and (not necessarily exact) order $\de$, respectively. Let furthermore, $L_\DD^\de$ be the set of solutions of the unsigned problem $P_\DD^u (A_\infty)$, then we clearly have
$$
S_\DD^\dt \subset C_\DD^\de \subset L_\DD^\de \, ,\qquad C_\DD^\de = \bigcup_{\de \leq \dt \leq \ro +1} S_\DD^\dt\, ,
$$
where $3 \leq \de \leq \dt \leq \ro + 1$. Moreover, $L_\DD^\de$ is a linear space, whose dimension is determined by $\ro$ and $\de$. More precisesly the following theorem holds:
\begin{thm}
 Let $S_\DD^\de$, $C_\DD^\de$, and $L_\DD^\de$ be the solution sets of the signed direction problem with exact decay order $\de$, of the signed problem with decay order $\de$, and of the unsigned problem with decay order $\de$, respectively. If $3 \leq \de \leq \ro + 1$, we have 
 \be \label{2.36}
 L_\DD^\de = \langle C_\DD^\de \rangle = \langle S_\DD^\de \rangle \, ,
 \ee
 \be \label{2.37}
 \dim L_\DD^\de  = \ro - \de +2\, ,
 \ee
 where $\ro$ is the rotation number of $\DD$ and $\langle S\rangle$ denotes the real linear span of the set $S$. If $\de > \ro + 1$ we have $L_\DD^\de = \{0\}$.
\end{thm}
Some comments are in order:

1. Theorems 2.1 and 2.2 are robust in the sense that the theorems still hold for weights $\beta$ that vary in some interval around $1/5$. 

2. The restriction to an even number of zeroes in theorem 2.3 can supposedly be removed. The physically relevant case, however, is the unbounded one, where this restriction does not apply, and we saved us this effort.

3. In view of the regularity of the data $\Om$ (see appendix A), much more regularity than $u \in H^1_{loc}$ can not be expected.\footnote{ H\"older continuity of $u$, $q$, and $p$ is shown in section 9.}
However, when inserting $q = u + \Phi$ into (\ref{1.3}) the result will be more regular (by exponentiation  of $q$ and by multiplication by the zeroes of $h$). In fact, the harmonic vector field $\BB$ is known to be analytic in the exterior space $E$.

4. According to theorems 2.3 and 2.4, all the nonuniqueness of the signed direction problem is encoded in the arbitrary positions of the zeroes. Uniqueness (up to a positive constant factor) in $A_\infty$ is thus only guaranteed if $\de = \dt = \ro + 1$; if the direction field is the only data (and hence $\de = 3$) this requires $\ro =2$, which holds, e.g., for ``dipole-type'' direction fields.

5. The unsigned direction problem prescribes the direction of the field vector only up to a sign at $S_1$. Arbitrary linear combinations of different solutions for the same direction field $\DD$ are thus again solutions for $\DD$, forming the linear space $L_\DD^\de$. As to the signed problem only positive linear combinations are admissible, i.e. $C_\DD^\de$ has the representation  
\be \label{2.38}
C^\delta_\DD = \Big\{\sum^N_{n=1} \lambda_n\, \BB_n :\,\BB_n \in C^\delta_\DD\, , \; \lambda_n > 0\, , \; n =1, \ldots, N\, ,\; N \in
\nat\Big\} ,
\ee
which is a {\em cone} in $L_\DD^\de$.

As to the solution of problem $P_\Om (A_R)$ note that it does neither have a variational form nor does it have (obvious) monotonicity properties, but a favourable feature clearly is the boundedness of the nonlinear term. The method of choice to obtain global solutions of $P_\Om (A_R)$ is thus to solve a suitably linearized version and to define thereby a compact mapping to which Schauder's fixpoint principle applies. In the weak setting of eq.\ (\ref{2.29}), linearized by replacing $\aaa[u - \Om]$ by $\aaa[w - \Om]$ with given function $w$ , the Lax-Milgram criterion provides easily general solvability, if only (besides boundedness) some coercivity condition is satisfied. This latter condition amounts to
$$
\bigg| \int_{A_R} \frac{u}{\rho} \, \aaa \cdot \na u\, \rd \ze \rd \rho \bigg| \leq C\, \|\na  u\|^2_{L^2 (A_R)}
$$
for some $C < 1$. A superficial estimate of the left-hand side by (\ref{1.5}) with $\ga =0$ and $|a_\ze | \leq 0.73$, $|a_\rho | \leq 1$
(see Fig.\ 2.2), however, fails:
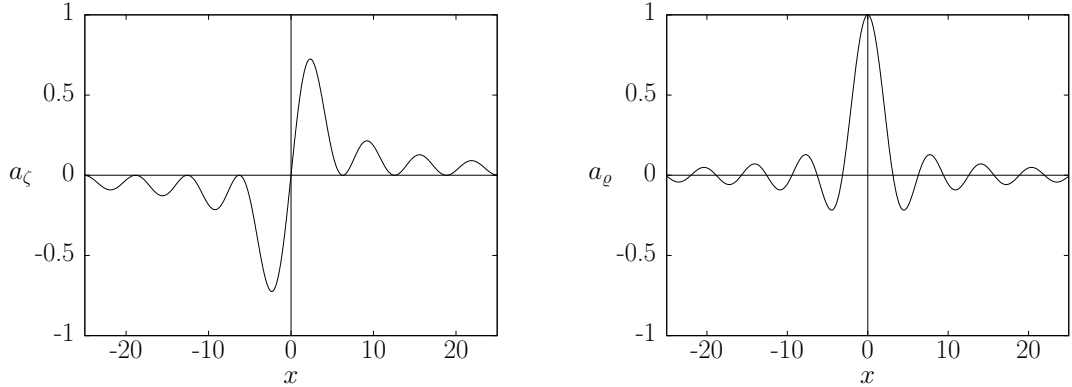
\begin{figure}
\begin{center}
\resizebox{0.48\textwidth}{!}{\input{Fig22cos}}
\resizebox{0.48\textwidth}{!}{\input{Fig22sin}}
\caption{Graphs of the functions $a_\ze : x\mapsto (1 - \cos x)/x$ and $a_\ro : x \mapsto \sin x /x$.}
\end{center}
\end{figure}
\be \label{2.39}
\bigg| \int_{A_R} \frac{u}{\rho} \, \aaa \cdot \na u\, \rd \ze \rd \rho \bigg| \leq 2\,\Big( 0.73\,  \|\pa_\ze u\|_{L^2 (A_R)}
\| \pa_\rho u\|_{L^2 (A_R)} + \| \pa_\rho u\|_{L^2 (A_R)}^2 \Big) .
\ee
A great deal of the present work is devoted to overcoming this problem by the combined effect of three measures: 

\noi (i) We introduce the variable $v:= \rho^\al u$ with suitable $\al > 0$, which has the effect that $a_\rho$ is shifted by $\al$ and allows a bound better than 1.

\noi (ii) A $\rho$-dependent weight can improve the optimal constant in the Hardy inequality (\ref{1.5}). But this requires a generalized ($\rho$-weighted) version of the Lax-Milgram criterion.

\noi (iii) $\pa_\ze u$ and $\pa_\rho u$ do not appear symmetrically in inequality (\ref{2.39}). A weighted gradient $\na_d := d\, \eee_\ze\, \pa_\ze + \eee_\rho\, \pa_\rho$, $d > 0$ can exploit this for further improvement.

Unfortunately, in the generalized Lax-Milgram criterion the coercivity condition is now governed by constants $C_c$ and $\Ct_c$, which are defined by min-max problems, viz.\ 
\be \label{2.40}
\disp C_c := \inf_v \, \sup_\vii\frac{\int_{A_R} \na v \cdot \na \vii\, \rho^{-\al}\, \rd \ze \rd \rho}{\Big( \int_{A_R} |\na_d\, v |^2\, \rho^{-\al - \ga}\, \rd \ze \rd \rho\Big)^{1/2} \Big( \int_{A_R} |\na_e\, \vii |^2\, \rho^{-\al + \ga}\, \rd \ze \rd \rho\Big)^{1/2}
}
\ee
and similarly for $\Ct_c$. $v$ and $\vii$ vary here in differently weighted versions of $H_0^1 (A_R)$. The case $\ga = 0$, $d = e =1$
corresponds to the ordinary  coercivity constant $C_c = 1$; increasing $\ga$ leads to decreasing $C_c$-values thus strengthening the coercivity condition. It needs a subtle balance of all involved parameters ($\al$, $\ga$, $d$, and $e$) to meet, finally, the coercivity condition of the generalized Lax-Milgram criterion.

To obtain sufficiently sharp lower bounds on $C_c$ we proceed as follows. Firstly the annular region $A_R$ is replaced by a rectangle $Q$, which allows us to split the two-dimensional problem in a sequence of one-dimensional problems in the variable $\rho$. The Euler-Lagrange equations for these problems amount to low-dimensional systems of ordinary linear differential equations, which constitute eigenvalue problems whose minimum eigenvalue bounds $C_c$. The numerical solution of these equations has heuristic value in that it allows us to identify appropriate values of the above parameters. A special case can be solved fully analytically and this solution hints to the kind of test function that, finally, yields rigorous lower bounds on $C_c$ sufficient for our needs. 

Once the solution of the linearized problem is established, a successful iteration depends crucially on the underlying space, which is here an $L^p$-space that is again suitably weighted by some power of $\rho$. $L^p$-generalizations of (\ref{1.5}) of the form
$$
\Big\|\frac{f}{\rho^\de} \Big\|^2_{L^p (A_R)} \leq C \int_{A_R} |\na_d\, f|^2\, \frac{\rd \ze \rd \rho}{\rho^\beta}
$$
with as large as possible values of $p$ and $\de$ for given $\beta$ play here a major role. Optimal estimates of this type in two dimensions for ``small'' and ``large'' values of $p$ are the key to prove continuity of the nonlinear iteration mapping. Compactness of the mapping then is a comparatively easy consequence of well-known embedding theorems. 

Schauder's fixed point principle does not provide uniqueness; so, this issue has to be faced separately. A heuristic consideration exploiting the divergence-character of eq.\ (\ref{2.23}) suggests the kind of test function that would exclude nontrivial solutions of the corresponding equation for the difference of two solutions. In fact, in order to prove rigorously uniqueness in $H_0^1 (A_R)$ it takes a whole sequence of test functions and special attention has to be given to the behaviour at the symmetry axis. In $A_\infty$ the weight $\ro^{-\beta}$ in (\ref{2.33}) cannot be neglegted. Starting point is now a formulation of eq.\ (\ref{2.23}) that incorporates the weight while keeping the divergence structure. The proof then proceeds quite analogously to the bounded case providing, finally, uniqueness 
for functions in $H_{loc,as}^1 (A_\infty)$ satisfying the bound (\ref{2.33}) for $0< \beta <1$.
Recall that all the nonuniqueness that is typical for the solutions $\BB$ of the direction problem is encoded in the zero-positions angle $\Psi$, which is part of the data in eq.\ (\ref{2.23}).

The transition from $A_R$ to the unbounded region $A_\infty$ proceeds by a suitable sequence of solutions $(u_n)$ defined on $A_n$ with boundary values on $S_n$ that are obtained from the exterior harmonic potential with boundary function $\phi$ on $S_1$. Restricting $u_n$ on $A_m$, $m \leq n$ yields actually a double sequence $(u_{m,n})$ with uniform bound (\ref{2.31}). A suitably defined diagonal sequence $(u^{(k)})$ then has the favourable property that $u^{(k+1)}$ extends $u^{(k)}$ defined on $A_k$ onto $A_{k+1}$. Thus $(u^{(k)})$ allows the definition of a function $u$ on $A_\infty$ that satisfies the bound (\ref{2.33}) and that is in fact a weak solution of $P_\Om (A_\infty)$. 
Once $u$ and hence $q$ are known, $p$ is determined by (\ref{2.20}) up to a constant $p_0$. By substitution into the ansatz (\ref{2.13}) one obtains, finally, $f$, which corresponds to a weakly harmonic field $\BB$. Higher interior regularity then follows from standard elliptic regularity theory, whereas continuity up to the boundary requires some more subtle arguments depending on the H\"older continuity of the boundary data.

So far zeroes on the symmetry axis are not allowed mainly for the technical reason not to loose favourable properties of the zero-positions angle $\Psi$ at the symmetry axis. So the number of zeroes in $A_R$ must be even, which means, e.g.\ for $\de =3$, a restriction of possible direction fields to those with even rotation numbers (see (\ref{2.7})). This limitation can be overcome by taking suitable linear combinations of ``even solutions''. The coefficients of a linear combination can be viewed as ``deformation parameters'' that govern the positions of the zeroes. Based on invariance properties of the degree of mapping with respect to continuous deformations (which are recalled in this context), single zeroes can be eliminated from $A_\infty$ by ``pushing'' them to infinity. A crucial point is here to keep control over the other zeroes, which will move but which had to avoid the boundary.

Finally, we characterize and, in particular, determine the dimension of the solution space of the {\em unsigned} direction problem. The presentation follows here largely the corresponding one in the two-dimensional case in (Kaiser 2010). The basic result is that $L_\DD^\de$ is generated by an arbitrary set of $\ro - \de + 2$ solutions of the signed problem with precisely $0, 1, \ldots , \ro - \de +1$ zeroes. In this sense no fundamentally new solutions appear in the unsigned problem  and $L_\DD^\de$ can be viewed as a convenient way to quantify the nonuniqueness of the direction problem.

The material just described is organized in the following sections: section 3 collects the various Hardy-type inequalities we make use of in the course of the proof. Uniqueness in the problem $P_\Om$, an issue that is independent of the rest of the paper, is proved in section 4. Sections 5 and 6 are devoted to the linearized problem and, as an essential part of it, to the min-max problem. Sections 7 and 8 present solutions of $P_\Om$ in annuli $A_R$ and in the exterior plane $A_\infty$, respectively. The direction problem itself, again in 
$A_R$ and in $A_\infty$, is solved in section 9 under the restriction that no zeroes are lying on the symmetry axis. 
Zeroes on the symmetry axis are discussed in section 10 and the unsigned problem in section 11. Some more technical estimates and some additional material are deferred to a number of appendices: appendices A and B contain estimates of the zero-positions angle $\Psi$ and of the boundary function $\Phi$, respectively. Appendix C contains a proof of the generalized Lax-Milgram criterion. Appendices D and E contain the analytic solution of a special one-dimensional min-max problem and numerical solutions of the general one-dimensional problem, respectively. Finally, appendix F contains an explicit exemplary solution of the 2D-direction problem that illustrates the migration of a zero.   
\sect{Hardy-type inequalities}
The results 3.1 -- 3.3 are formulated for (not necessarily bounded) domains $G$ in the $n+1$-dimensional half-space 
$$
H^{(n+1)} := \{(x_1, \ldots ,x_n, y) \in \real^{n+1} : y > 0\}\, ,\quad n \in \nat\, ,
$$ 
without causing additional effort. All other results hold for bounded domains $G \subset H$ with
$$
H:= H^{(2)} = \{(x, y) \in \real^{2} : y > 0\}\, .
$$ 
Note for subsequent applications the correspondence $(x,y) \triangleq (\ze, \rho)$. 

$y^\ga$-weighted $L^p$-norms, especially with $p = 2$, play in this paper a dominant role. The following notation is here useful:
\be \label{3.1}
\| \cdot \|_{p, \ga} := \bigg( \int_G | \cdot |^p\, \rd \mu_\ga \bigg)^{1/p} , \qquad \rd \mu_\ga := y^{-\ga} \rd x_1 \ldots \rd x_n \rd y
\ee
with given domain $G \subset H^{(n+1)}$, $p \geq 1$, and $\ga \in \real$. For $p= 2$ we use the simplified notation: $\|\cdot \|_{2 ,\ga} =: \| \cdot \|_{\ga}$, and for $p = \infty$: $\|\cdot\|_{\infty, 0} =: \|\cdot \|_\infty$. The following function spaces are associated to these norms:
\be \label{3.2}
\HC^{(c)}_{p, \ga} (G) := clos\big(C_0^\infty (G)\, ,\, \|\na_c\cdot \|_{p, \ga}\big)\, ,\qquad \HC^{(c)}_\ga (G) := \HC^{(c)}_{2, \ga} (G)
\ee
where
$$
\na_c := (c\, \na_x\, ,\, \pa_y )\, ,\quad c \geq 0\, .
$$
For {\em bounded} $G\subset H^{(n+1)}$ we have the well-known inclusion
$$
\HC^{(c)}_{p,\ga} (G) \subset \HC^{(c)}_{q,\ga} (G) \qquad \mbox{for } \; p \geq q\, ,
$$
as well as 
\be \label{3.3}
\HC^{(c)}_{p,\ga} (G) \subset \HC^{(c)}_{p,\de} (G) \qquad \mbox{for } \; \ga \geq \de\, ,
\ee
where the latter inclusion is in fact an equivalence in the case $G \Subset H$, i.e.\ that $G$ is compactly contained in $H$. With respect to
$c$ only the difference between $c=0$ and $c > 0$ matters:\footnote{To simplify the notation we sometimes omit the upper index $c$, which means that $c> 0$, or omit the indication of $G$ if the underlying domain is clear from the context.}
$$
\HC^{(c)}_{p,\ga} (G) = \HC^{(d)}_{p,\ga} (G)  \subset \HC^{(0)}_{p,\ga} (G)\, , \qquad c, d > 0\, .
$$
The following proposition is a weighted $L^p$-version of Hardy's inequality (see Hardy et al.\ 1952, p.\ 175) suitable to our needs.
\begin{prop}
 Let $G \subset H^{(n+1)}$ be a domain and $f \in \HC^{(c)}_{p,\ga} (G)$ with $\ga \neq 1-p$ and $p \geq 1$. Then, the following inequalities hold:
 \be \label{3.4}
 \Big\|\frac{f}{y}\Big\|_{p,\ga} \leq \frac{p}{|p + \ga -1|}\, \|\pa_y f\|_{p,\ga}
 \ee
 and, especially for $p=2$,
 \be \label{3.5}
 \Big\|\frac{f}{y}\Big\|_{\ga} \leq \frac{2}{|1 + \ga |}\, \|\pa_y f\|_{\ga}\, .
 \ee 
\end{prop}
\textsc{Proof:}
Let $\vii \in C_0^\infty (G) \subset C_0^\infty (H)$ and let us define $\vit \in C_0^\infty (\real_+)$ by fixing $ (x_1, \ldots, x_n) = \xx \in \real^n$ in $\vii$:
$$
\vit = \vit_\xx :\,  \real_+ \ra \real\, ,\quad y \mapsto \vii (\xx, y)\, .
$$
Integrating by parts then yields
$$
\ba{l}
\disp \int_{\real_+} \frac{|\vit|^p}{y^{p + \ga}}\, \rd y = \frac{-1}{p + \ga -1} \int_{\real_+} |\vit|^p \Big(\frac{1}{y^{p + \ga -1}}\Big)' \rd y \\[3ex]
\disp \quad = \frac{p}{p + \ga -1} \int_{\real_+} |\vit|^{p-2}\, \vit \,\vit'\, \frac{1}{y^{p + \ga -1}}\, \rd y 
\leq \frac{p}{|p + \ga -1|} \int_{\real_+} \frac{|\vit|^{p-1}}{y^{(p + \ga)(p-1)/p}}\, \frac{|\vit'|}{y^{\ga/p}}\, \rd y \\[3ex]
\disp \qquad \leq \frac{p}{|p + \ga -1|} \bigg(\int_{\real_+} \frac{|\vit|^{p}}{y^{p + \ga}}\, \rd y \bigg)^{\frac{p-1}{p}}
\bigg(\int_{\real_+} \frac{|\vit'|^p}{y^{\ga}}\, \rd y \bigg)^{\frac{1}{p}} ,
\ea
$$
where in the last line we applied H\"older's inequality. By cancellation one obtains
$$
\int_{\real_+} \Big|\frac{\vit}{y}\Big|^p y^{-\ga}\, \rd y \leq \bigg(\frac{p}{|p + \ga -1|}\bigg)^p \int_{\real_+} |\vit'|^p\, y^{-\ga}\, \rd y\, .
$$
Finally, applying Fubini's theorem in the form 
$$
\int_H \chi (\xx, y)\, \rd x_1 \ldots \rd x_n \rd y = \int_{\real^n} \bigg( \int_{\real_+} \widetilde{\chi}_\xx (y) \, \rd y \bigg) \rd x_1 \ldots \rd x_n
$$
yields the assertion for $\vii \in C_0^\infty (G)$ and by approximation for $f \in \HC^{(c)}_{p,\ga}$.
\qed\\[-2ex]

\noi {\em Remark}: In the case that $G$ is contained in the half-ball $B_R^+ \subset H^{(n+1)}$ of radius $R$ centered at the origin, inequality (\ref{3.4}) provides immediately a Poincar\' e-type inequality:
 \be \label{3.6}
 \|f\|_{p,\ga} \leq R\, \frac{p}{|p + \ga -1|}\, \|\pa_y f\|_{p,\ga} \leq R\, \frac{p}{|p + \ga -1|}\, \|\na_c\, f\|_{p,\ga}\, .
 \ee

The following lemma demonstrates that the constant in inequality (\ref{3.5}) cannot be improved in the case that $G$ contains a box that touches the symmetry axis.
\begin{lem}[box-criterion]
Let $G\subset H^{(n+1)}$ be a domain and $Q \subset G$ a box of the form $Q = Q^{(n)} \ti (a,b) \subset \real^n \ti \real_+$ with $0< a < b$. Then, any admissible constant $c$ in inequality (\ref{3.5}) satisfies 
\be \label{3.6a}
c \geq \bigg[\Big(\frac{1 + \ga}{2}\Big)^2 + \Big( \frac{\pi}{\ln (b/a)}\Big)^2 \bigg]^{-1/2} .
\ee
In particular, if $G$ contains a box of the form $Q^{(n)} \ti (0,b)$, the constant in (\ref{3.5}) is the smallest possible. 
\end{lem}
\textsc{Proof:} The smallest constant $c_{min}$ in (\ref{3.5}) is associated to a variational problem, viz.,  
\be \label{3.7}
\disp c_{min}^{-2} = \inf_{\vii \in C_0^\infty (G)}\, \frac{\int_G |\pa_y \vii |^2 \,\rd \mu_\ga}{\int_G |\vii /y |^2 \,\rd \mu_\ga}\, .
\ee
It suffices for our purposes to consider the following simpler, one-dimensional version:
\be \label{3.8}
\inf_{f \in H_0^1 ((a,b))}\, \frac{\int_a^b f'^2\, y^{-\ga}\, \rd y}{\int_a^b (f /y )^2\, y^{-\ga} \,\rd y} = :\la_{min}\, .
\ee
Problem (\ref{3.8}) is in fact a standard problem in the calculus of variation. The associated Euler-Lagrange equations,
\
$$
\ba{rl}
\disp f'' -\frac{\ga}{y}\, f' + \frac{\la}{y^2} \, f = 0 & \mbox{ in } \, (a,b)\, ,\\[2ex]
f = 0 & \mbox{ at }\, \{a, b\}\, ,
\ea
$$
constitute an eigenvalue problem in $\la$, which can explicitly be solved. The minimizer in (\ref{3.8}) is the eigenfunction with smallest eigenvalue:
$$
\ba{c}
\disp f_{min} (y) = y^{(1 + \ga)/2}\, \sin \Big( \pi\, \frac{\ln (y/a)}{\ln (b/a)} \Big) ,\\[3ex]
\disp \la_{min} = \Big(\frac{1 + \ga}{2}\Big)^2 + \Big( \frac{\pi}{\ln (b/a)}\Big)^2 .
\ea
$$
Let now $\Qt := \Qt^{(n)} \ti (a,b)$ be an extension of $Q := Q^{(n)} \ti (a,b)$ such that $Q^{(n)} \Subset \Qt^{(n)}$ and $|\Qt \setminus Q| < \e$. Let, furthermore, $f : Q \ra \real$ be given by $f(\cdot , y) = f_{min} (y)$ and $\ft$ be a $C^1$-extension onto $\Qt$ such that $\ft |_{\pa \Qt}  = 0$ and 
$$
\max_{\Qt} |\pa_y \ft | \leq \max_Q |\pa_y f| = \max_{[a,b]} |f'_{min}|\, .
$$
For $\ft$ we then have the estimate
\be \label{3.9}
\ba{l}
\disp \frac{\int_\Qt |\pa_y \ft |^2 \,\rd \mu_\ga}{\int_\Qt |\ft /y |^2 \,\rd \mu_\ga} \leq 
\frac{\int_Q |\pa_y f |^2 \,\rd \mu_\ga + \int_{\Qt \setminus Q} |\pa_y \ft|^2 \rd \mu_\ga}{\int_Q |f /y |^2 \,\rd \mu_\ga} \\[3ex]
\disp \quad \leq \frac{\int_a^b f'^{\, 2}_{min}\, y^{-\ga}\, \rd y}{\int_a^b (f_{min} /y )^2\, y^{-\ga} \,\rd y} 
+ \e \, \frac{\max_{[a,b]} \{f'^{\, 2}_{min}\, y^{- \ga}\}}{\int_Q |f/y»|^2\, \rd \mu_\ga} \leq \la_{min} + \e\, \frac{b^2 (b/a)^\ga}{|Q^{(n)}|}\, \frac{\max_{[a,b]} f'^{\, 2}_{min}}{\int_a^b f^{\, 2}_{min}\, \rd y}\, .
\ea
\ee
As $C_0^\infty (\Qt)$ is dense in $H_0^1 (\Qt)$ we can approximate $\ft$ on the left-hand side of (\ref{3.9}) by $\vii \in C_0^\infty  (\Qt) \subset C_0^\infty (G)$ and find the infimum in (\ref{3.7}) be bounded from above by the right-hand side in (\ref{3.9}). As $\e >0$ is arbitrary we thus obtain
$$
c_{min}^{-2} \leq \la_{min}\, ,
$$
which is (\ref{3.6a}).

If $G$ contains a box of the form $Q^{(n)} \ti (0,b)$, the first assertion holds for any $a \in (0,b)$, i.e.\ 
$$
\frac{2}{|1 + \ga|} \leq c_{min} \leq \frac{2}{|1+ \ga|}\, ,
$$
which is the second assertion of the box-criterion. 
\qed

Inequality (\ref{3.5}) allows some alternative characterizations of $\HC^{(c)}_\ga (G)$, which will be useful when dealing with the min-max problem. 
\begin{lem}
Let $G \subset H^{(n+1)}$ be a domain and $-1 \neq \ga \in \real$. On $\HC^{(c)}_\ga (G)$ we then have the equivalence of norms
\be \label{3.10}
\Ct_\ga\, \|\na_c \, \cdot \|_\ga  \leq \|\na_c \, (y^{-\ga/2}\, \cdot\, ) \|_0 \leq C_\ga\, \|\na_c\, \cdot \|_\ga
\ee
with constants
$$
C_\ga := 1 + \Big|\frac{\ga}{1 + \ga}\Big|\; ,\qquad \Ct_\ga := \min \Big\{ 1\, ,\, \frac{1}{|1 + \ga|} \Big\} .
$$
\end{lem}
\textsc{Proof:} By 
$$| \pa_y (y^{-\ga/2} \vii)|^2 = \Big[|\pa_y \vii|^2 - \ga \, \frac{\vii}{y}\, \pa_y \vii + \frac{\ga^2}{4} \Big(\frac{\vii}{y} \Big)^2\, \Big] y^{-\ga} ,
$$
and repeated use of (\ref{3.5}) one obtains for $\vii \in C_0^\infty (G)$:
$$
\ba{l}
\disp \| \na_c (y^{-\ga/2} \vii)\|^2_0 \leq \|\na_c \, \vii\|^2_\ga + |\ga|\, \Big\|\frac{\vii}{y}\Big\|_\ga  \| \pa_y \vii\|_\ga + \frac{\ga^2}{4}\, \Big\|\frac{\vii}{y}\Big\|^2_\ga \\[3ex]
\disp \qquad \leq \Big(1 + \frac{2\, |\ga|}{|1 + \ga|} + \frac{\ga^2}{(1+\ga)^2}\Big) \|\na_c \,\vii\|^2_\ga = \Big(1 + \frac{|\ga|}{|1 + \ga|}\Big)^2 \|\na_c \, \vii\|^2_\ga \, .
\ea
$$
Similarly, by
$$- \ga \int_G \vii\, \pa_y \vii\, y^{-(\ga + 1)} \rd y \rd x_1 \ldots \rd x_n = -\frac{1}{2}\, \ga (\ga +1) \int_G \vii^2 y^{-(\ga + 2)} \rd y \rd x_1 \ldots \rd x_n
$$
and again by (\ref{3.5}) one obtains 
$$
\ba{l}
\disp \| \na_c (y^{-\ga/2} \vii) \|^2_0 = \int_G \Big(c^2 |\na_x\, \vii|^2 + \frac{1}{(1+ \ga)^2} \, |\pa_y \vii|^2 \Big) y^{-\ga} \rd x_1\ldots  \rd x_n \rd y \\[2ex]
\disp \qquad \qquad \qquad \quad  + \Big( 1 - \frac{1}{(1+ \ga)^2}\Big) \|\pa_y \vii\|^2_\ga  - \frac{1}{2} \, \ga (\ga +1) \Big\| \frac{\vii}{y} \Big\|^2_\ga  + \frac{\ga^2}{4} \, \Big\| \frac{\vii}{y}\Big\|^2_\ga \\[3ex] 
\disp \qquad  \geq \min \big\{ 1\, , (1+ \ga)^{-2} \big\}\, \|\na_c \, \vii\|^2_\ga + \frac{1}{4}\, \big((1+ \ga)^2 - 1 - 2\, \ga(\ga +1) + \ga^2 \big) \Big\|\frac{\vii}{y}\Big\|^2_\ga  \\[3ex]
\disp \qquad \quad = \min \big\{ 1\, , (1+ \ga)^{-2} \big\}\, \|\na_c \, \vii\|^2_\ga\, .
\ea
$$
\qed

$H_0^1  (G)$ denotes the usual Sobolev space $clos \big(C_0^\infty (G) , \|\cdot \|_{H^1} \big)$ with 
$$
\|\vii \|^2_{H^1} = \int_G \big(|\vii|^2 + |\na_x\, \vii|^2 + | \pa_y \vii|^2 \big)\, \rd x_1 \ldots \rd x_n \rd  y\, .
$$
By (\ref{3.6}) we have thus the norm equivalence $\|\cdot \|_{H^1} \sim \|\na_c \cdot \|_0$ and hence the identity
$$
H^1_0 (G) = \HC _0^{(c)} (G) \, ,
$$
provided that $G$ is bounded and that $c>0$. More generally, for $\ga  \neq -1$, (\ref{3.6}) and (\ref{3.10}) imply under these conditions:
\be \label{3.11}
\ba{rl}
\HC^{(c)}_\ga  (G)\!\!\! &= clos \big( C_0^\infty (G)\, ,\, \|y^{-\ga/2} \cdot \|_{H^1} \big) = \big\{ \vii : y^{-\ga/2} \vii \in H_0^1 (G) \big\} \\[2ex]
& = \big\{ y^{\ga/2} \chi : \chi \in H^1_0 (G) \big\}\, ,
\ea
\ee
and (\ref{3.3}) implies 
\be \label{3.12}
\HC_\ga^{(c)}(G) \subset H_0^1 (G) \qquad \mbox{for }\; \ga \geq 0 \, . 
\ee
Note, finally, that $\HC_\ga^{(c)} (G)$, $c>0$ equipped with the scalar product
\be \label{3.12a}
(\vii , \chi) \mapsto \int_G \na_c \, \vii \cdot \na_c \, \chi \; y^{-\ga} \rd x_1 \ldots \rd x_n \rd y
\ee
is a Hilbert space, which implies in particular that $\HC_\ga^{(c)} (G)$ is a reflexive space.

The rest of this section is devoted to inequalities of type 
\be \label{3.13}
\Big\|\frac{f}{y^\de}\Big\|_{p,0} \leq C \, \|\na_c \, f \|_{q, \beta}\, ,
\ee
which will be necessary in section 7. $G$ is now a bounded domain contained in some half-ball $B_R^+ \subset H \subset \real^2$ and we suppose $c> 0$. The focus is now on ``optimal values'' of $p$ and $\de$ for given values of $q$, especially for $q =2$, and $\beta$. 

We start with some one-dimensional inequalities, which are comparatively easy to derive. Let $f \in H^1((0,R))$ with $f(0) = 0$. By the fundamental theorem of calculus and by H\"older's inequality one obtains:
$$
\ba{l}
\disp |f(y)| \leq \int_0^y |f'|\, \rd z \leq \bigg( \int_0^y z^{\beta/(q-1)}\, \rd z \bigg)^{1-\frac{1}{q}} \bigg( \int_0^y |f'|^q\, z^{-\beta} \rd z \bigg)^{\frac{1}{q}} \\[3ex]
\disp \qquad \leq \Big( 1+ \frac{\beta}{q -1} \Big)^{\frac{1}{q} -1}\, y^{1 + \frac{\beta -1}{q} } \bigg( \int_0^y |f'|^q\, z^{-\beta} \rd z \bigg)^{\frac{1}{q}} ,
\ea
$$
i.e.\
\be \label{3.14}
\Big\|\frac{f}{y^\de}\Big\|_{\infty,0} \leq \Big(1 + \frac{\beta}{q -1} \Big)^{-\big(1-\frac{1}{q}\big)} \|f'\|_{q, \beta}\, ,
\ee
where
$$
\de = 1 + \frac{\beta - 1}{q}\, ,\qquad 1 < q <\infty\, ,  \qquad  \beta > 1- q\, .
$$
Inequality (\ref{3.14}) contains the special case
\be \label{3.15}
\|f\, y^{-(1+ \beta)/2} \|_{\infty,0} \leq (1 + \beta)^{-1/2}\, \|f'\|_{\beta}\, , \qquad \beta > -1 \, ,
\ee
and (formally) the limit cases:
$$
\| f/y\|_{\infty,0} \leq \|f' \|_{\infty,0} \, ,
$$
\be \label{3.16}
\| f/y^\beta \|_{\infty,0} \leq \|f' \|_{1,\beta} \, .
\ee
Inequality (\ref{3.16}) is clear for $\beta = 0$. Otherwise we have 
$$
f(y)\, y^{-\beta} = \int_0^y f' z^{-\beta} \rd z - \beta \int_0^y f\, z^{-1-\beta} \rd z \, ,
$$
$$
f(y)\, y^{-\beta} = - \int_y^R f' z^{-\beta} \rd z + \beta \int_y^R f\, z^{-1-\beta} \rd z \, .
$$
Thus, by summation and using (\ref{3.4}) with $p= 1$ one obtains
$$
2\ |f(y)|\, y^{-\beta} \leq  \int_0^R |f'|\, z^{-\beta} \rd z + |\beta| \int_0^R (|f|/z)\, z^{-\beta} \rd z 
\leq 2 \int_0^R |f'|\, z^{-\beta} \rd z \, ,
$$
which is (\ref{3.16}). 
Finally, interpolation between (\ref{3.15}) and (\ref{3.5})\footnote{Note that the proof of proposition 3.1 works as well for functions $f: (0, R) \ra \real $ without zero-boundary-condition at $R$.}
yields a (one-dimensional) inequality of type (\ref{3.13}): 
$$
\ba{l}
\disp \int_0^R |f /y^\de |^p\, \rd y = \int_0^R \big|f\, y^{-(1 + \beta)/2} \big|^{p-2}\, |f/y|^2\, y^{-\beta} \rd y \leq \\[3ex]
\disp \qquad \leq \|f \,y^{-(1+ \beta)/2} \|_{\infty, 0}^{p-2}\ \|f/y\|^2_\beta \leq (1 + \beta)^{\frac{2-p}{2}}\, 2^2 (1 + \beta)^{-2} \, \|f'\|^p_\beta \, , \quad  \de := (1+\beta)/2 + 1/p \, ,
\ea
$$
i.e.\ 
\be \label{3.17}
\Big\|\frac{f}{y^\de}\Big\|_{p, 0} \leq 2^{\frac{2}{p}} \Big(\frac{1}{1 +\beta}\Big)^{\frac{p+2}{2p}} \|f'\|_\beta
\ee
with
$$
\de = \frac{1+\beta}{2} + \frac{1}{p}\, , \qquad 2 \leq p \leq \infty\, , \qquad \beta > -1\, .
$$
Inequalities (\ref{3.14}) and (\ref{3.17}) are optimal in the sense that $\de$ cannot be enlarged as can easily be seen by testing the inequalities with $f(y) = y^\al$. 

In 2 dimensions we must proceed differently since a result of type (\ref{3.14}) with $q =2$ cannot be achieved (not even for $\de = \beta =0$). We thus assume $p < \infty$ and distinguish, moreover, between ``small $p$'' and ``large $p$''. 
\begin{prop}[small-p-case]
Let $G$ be a domain contained in $B_R^+ \subset H \subset \real^2$ and $f \in \HC_\beta^{(c)} (G)$ with $\beta > -1$. Then, the following inequality holds
\be \label{3.18}
\Big\| \frac{f}{y^\de} \Big\|_{p, 0} \leq \Big(\frac{2}{1+ \beta}\Big)^{\frac{6-p}{2 p}} \Big(\frac{R}{c^2}\Big)^{\frac{p-2}{ 2p}}\, \| \na_c\, f\|_\beta
\ee
with 
$$
\de = \frac{\beta}{2} + \frac{6-p}{2 p} \, ,\qquad 2 \leq p \leq 4 \, ,\qquad \beta > -1 \, .
$$
\end{prop}
\textsc{Proof:} It is sufficient to prove (\ref{3.18})  for functions $\vii \in C_0^\infty (G)$. By (\ref{3.15}) we have 
$$
|\vii (x,y) |^2\, y^{-(1 + \beta)} \leq (1 + \beta)^{-1} \int_0^R |\pa_y \vii|^2\, y^{-\beta} \rd y
$$
and by the fundamental theorem
$$
|\vii (x,y) |^2 \leq \frac{2 R}{c^2}  \int_{-R}^R |c\, \pa_x \vii|^2\, \rd x \, .
$$
With these estimates one obtains
$$
\ba{l}
\disp \int_G |\vii|^4\, y^{-(1+ \beta)} \rd \mu_\beta = \int_0^R \int_{-R}^R |\vii|^4 \, y^{-(1 + 2\beta)} \rd x \rd y \\[3ex]
\disp \quad \leq \int_{-R}^R \sup_y \big\{|\vii (x,y)|^2\, y^{-(1+\beta)} \big\} \rd x \times \int_0^R \sup_x |\vii(x,y)|^2\, y^{-\beta} \rd y \\[3ex]
\disp \qquad \leq \frac{1}{1 + \beta} \int_{-R}^R \int_0^R |\pa_y \vii |^2\, y^{-\beta} \, \rd y \rd x \times \frac{2 R}{c^2} \int_0^R \int_{-R}^R |c \, \pa_x \vii|^2\, y^{-\beta} \rd x \rd y \\[3ex]
\disp \qquad \quad \leq \frac{2 R}{(1 + \beta) c^2} \bigg(\int_0^R \int_{-R}^R |\na_c\, \vii|^2\, y^{-\beta} \rd x \rd y \bigg)^2
= \frac{2 R}{(1 + \beta) c^2} \bigg(\int_G |\na_c\, \vii|^2\, \rd \mu_\beta \bigg)^2 .
\ea
$$
Interpolation with (\ref{3.5}) then yields
$$
\ba{l}
\disp \int_G |\vii/y^\de |^p\, \rd x \rd y = \int_G \big|\vii\, y^{-(1+ \beta)/4}\big|^{2(p-2)}\, |\vii/y|^{4-p} \, \rd \mu_\beta \\[3ex]
\disp \quad \leq \bigg(\int_G |\vii|^4\, y^{-(1+ \beta)} \rd \mu_\beta \bigg)^{\frac{p-2}{2}} \bigg( \int_G |\vii/y|^2\, \rd \mu_\beta\bigg)^{\frac{4-p}{2}} \\[3ex]
\disp \qquad  \leq \Big(\frac{2 R}{(1 + \beta) c^2} \Big)^{\frac{p-2}{2}} \Big(\frac{2}{1+ \beta}\Big)^{4 -p} \bigg(\int_G |\na_c \, \vii|^2 \,\rd \mu_\beta \bigg)^{\frac{p}{2}} ,
\ea
$$
where we set $\de := \beta/2 + (6-p)/2p$ and made use of H\"older's inequality in the second line. This is (\ref{3.18}). 
\qed

\begin{prop}[large-p-case]
Let $G$ be a domain contained in $B_R^+ \subset H \subset \real^2$ and $f \in \HC_{\frac{2 p}{2 + p} ,\bt}^{(c)} (G)$ with $p \geq 2$ and $\bt > 0$. Then, the following inequality holds
\be \label{3.19}
\Big\| \frac{f}{y^\dt} \Big\|_{p, 0} \leq \frac{p}{2 \sqrt{c}}\, \| \na_c\, f\|_{\frac{2 p}{ 2+p} , \bt}
\ee
with 
$$
\dt = \bt \,\frac{2 + p}{2 p} 
\, ,\qquad p \geq 2 \, ,\qquad \bt > 0 \, .
$$
Moreover, for $f \in \HC_\beta^{(c)} (G)$ with $\beta > -2 /p$ holds
\be \label{3.20}
\Big\| \frac{f}{y^\de} \Big\|_{p, 0} \leq \frac{p}{2 \sqrt{c}}\,\Big(\frac{2}{\e}\, R^{1+ \e} \Big)^{\frac{1}{ p}} \, \| \na_c\, f\|_\beta
\ee
with 
$$
\de = \frac{\beta}{2} + \frac{1-\e}{ p} \, ,\qquad p \geq 2 \, ,\qquad \beta > -\frac{2}{p}\, , \qquad \e >0 \, .
$$
\end{prop}
\textsc{Proof:} Let $\vii \in C_0^\infty (G)$. We proceed similarly as in the proof of proposition 3.4, starting this time, however, with (\ref{3.16}), i.e.\
$$
|\vii (x,y) |\, y^{-\beta} \leq \int_0^R |\pa_y \vii|\, y^{-\beta} \rd y
$$
and 
$$
|\vii (x,y) | \leq \int_{-R}^R |\pa_x \vii|\, \rd x \, .
$$
We then obtain
\be \label{3.21}
\ba{l}
\disp \int_G |\vii|^2\, y^{-\beta} \rd \mu_\beta = \int_0^R \int_{-R}^R |\vii|^2 \, y^{-2\beta} \rd x \rd y \\[3ex]
\disp \quad \leq \int_{-R}^R \sup_y \big\{|\vii (x,y)|\, y^{-\beta} \big\} \rd x \times \int_0^R \sup_x |\vii(x,y)|\, y^{-\beta} \rd y \\[3ex]
\disp \qquad \leq \int_{-R}^R \int_0^R |\pa_y \vii |\, y^{-\beta} \, \rd y \rd x \times \frac{1}{c} \int_0^R \int_{-R}^R |c \, \pa_x \vii|\, y^{-\beta} \rd x \rd y \\[3ex]
\disp \qquad \quad \leq \frac{1}{c} \bigg(\int_G |\na_c\, \vii|\, \rd \mu_\beta \bigg)^2 .
\ea
\ee

By approximation inequality (\ref{3.21}) holds for any $g\in \HC^{(c)}_{1, \beta} (G)$.
Inserting $g =: |f|^{p/2}$ and $2 \beta =: \bt (p/2 + 1)$, and using H\"older's inequality we obtain further
$$
\ba{l}
\disp \bigg( \int_G \big|f\, y^{-\bt/2}\big|^p \, \rd \mu_\bt \bigg)^{\frac{1}{2}} \leq \frac{1}{\sqrt{c}} \int_G \big| \na_c\, |f|^{p/2} \big|^2\, y^{-(\bt/2) (p/2 - 1)}\, \rd \mu_\bt \\[3ex]
\disp \quad \leq \frac{p}{2 \sqrt{c}} \int_G \big|f \, y^{-\bt/2} \big|^{p/2 -1} |\na_c\, f|\, \rd \mu_\bt \\[3ex]
\disp \qquad \leq \frac{p}{2 \sqrt{c}} \bigg(\int_G \big|f \, y^{-\bt/2} \big|^{p}\, \rd \mu_\bt \bigg)^{\frac{p -2}{2 p}} \bigg( \int_G |\na_c\, f|^{\frac{2 p}{2 + p}}\, \rd \mu_\bt\bigg)^{\frac{2 + p}{2 p}} .
\ea
$$
After cancellation we have 
$$
\bigg( \int_G \big|f\, y^{-\bt\frac{2 +p}{2 p}}\big|^p \, \rd x \rd y \bigg)^{\frac{1}{p}} \leq \frac{p}{2 \sqrt{c}} \bigg( \int_G |\na_c\, f|^{\frac{2 p}{2 + p}}\, \rd \mu_\bt\bigg)^{\frac{2 + p}{2 p}} ,
$$
which is (\ref{3.19}).

Using once more H\"older's inequality the right-hand side in (\ref{3.19}) can be linked up with the $L^2$-norm: 
\be \label{3.22}
\ba{l}
\disp \|\na_c \, f\|_{\frac{2p}{2 +p}, \bt} =  \bigg(\int_G |\na_c\, f|^{\frac{2p}{2 +p}} \, y^{-\bt} \, \rd x \rd y \bigg)^{\frac{2 + p}{2 p}} \\[3ex]
\disp \quad \leq \bigg( \int_G y^{-(1-\e)}\, \rd x \rd y\bigg)^{\frac{1}{p}} \bigg( \int_G |\na_c\, f|^2\, y^{-\beta}\, \rd x \rd y \bigg)^{\frac{1}{2}} \leq \Big(\frac{2 R\, R^\e}{\e}\Big)^{\frac{1}{p}}\, \| \na_c \, f\|_\beta
\ea
\ee
with $\e := 1 - \beta  + (2 + p)( \beta -\bt)/2$. The necessary condition $\e > 0$ implies $\bt < \beta + 2(1 -\beta)/(2 +p)$, which in turn by $\bt >0$ implies $\beta > -2 /p$. Combining (\ref{3.19}) with (\ref{3.22}) yields, finally, (\ref{3.20}) with $\de := \beta/2 + (1-\e)/p$. 
\qed
\sect{Uniqueness in the problems $P_\Om (A_R)$ and $P_\Om (A_\infty)$}
Uniqueness for fixed function $\Om$ means, in particular, uniqueness for fixed boundary values represented by the harmonic function $\Phi$ and fixed set of zero-positions represented by the angle $\Psi$. Together with the (up to a constant $p_0$) unique solution of (\ref{2.20}) for given $q$, this implies by (\ref{2.13}) an (up to a positive constant factor) unique solution $f$ of the direction problem.

Let $u_1$ and $u_2$ be weak solutions of $P_\Om (A_R)$, $R >1$. Then, $\du := u_1 - u_2$ satisfies weakly the following ``perturbation equation'':
\be \label{4.1}
\na \cdot \Big( \na\, \du + \aaa'\, \frac{\du}{\rho}\Big) = 0 \qquad \mbox{in }\, A_R\, , \qquad \quad \du\big|_{\pa A_R} = 0
\ee
with 
\be \label{4.2}
\left. \ba{l}
\disp a'_\ze := -\cos (u_1 - \Om ) + \cos (u_2 - \Om) = \sin \Big( \frac{1}{2}\, (u_1 + u_2 ) - \Om \Big)\, \frac{\sin \big( (u_1 - u_2)/2\big)}{(u_1 - u_2)/2} \, ,\\[3ex]
\disp a'_\rho := \sin (u_1 - \Om ) - \sin (u_2 - \Om) = \cos \Big( \frac{1}{2}\, (u_1 + u_2 ) - \Om \Big)\, \frac{\sin \big( (u_1 - u_2)/2\big)}{(u_1 - u_2)/2} \, .
\ea
\right\}
\ee
Note that given $u_1$ and $u_2$, $\aaa'$ is here considered as a given (bounded) vector field on $A_R$ and (\ref{4.1}) is thus a linear equation in $\du$.
\begin{prop}
 Let $u_1$, $u_2 \in H_{0, as}^1 (A_R)$ be weak solutions of problem $P_\Om (A_R)$ with $R>1$ and $\Om \in L^\infty (A_R)$ satisfying the bound (\ref{2.30}). Then
 $$
 u_1 = u_2 \qquad a.\, e.\ \mbox{in } \, A_R\, .
 $$
 \end{prop}
\textsc{Proof:} For functions $f \in H^1_{0, as} (A_R)$ we have by definition $trace\, f\big|_{\{\rho = 0\}} = 0$; it is thus sufficient to consider functions $f \in H_0^1 (A^+_R)$, where $A_R^+ := A_R \cap \{\rho >0\} \subset H$. 

Let us start with a heuristic consideration that motivates the choice of test functions we will make use of in the following. Let $\du^+ := \max \{\du, 0\}$ and integrate (\ref{4.1}) over $supp\,\du^+$:
\be \label{4.4}
\int_{supp\, \du^+} \na \cdot \Big(\na\, \du^+ + \aaa'\, \frac{\du^+}{\rho} \Big)\, \rd \ze \rd \rho = 0\, .
\ee
On the assumption that $\du^+$ and $supp\, \du^+$ are such that Gauss' theorem is applicable one obtains
\be \label{4.5}
\int_{\pa(supp\, \du^+)} \nn \cdot \na\, \du^+ \rd s + \int_{\pa(supp\,\du^+)} \nn \cdot \aaa'\, \frac{\du^+}{\rho} \, \rd s = 0\, ,
\ee
Where $\nn$ denotes the exterior normal at $\pa (supp\, \du^+)$. By definition we have $\du^+\big|_{\pa(supp\,\du^+)} = 0$ and $\nn \cdot \na\, \du^+\big|_{\pa(supp\,\du^+)} \leq 0$ and thus from (\ref{4.5}) we conclude:
\be \label{4.6}
\du^+ = \na\, \du^+ = 0 \qquad \mbox{on } \ \pa(supp\,\du^+)\, .
\ee
This conclusion continues even in the case that $\pa(supp\, \du^+) \cap \{\rho =0\} \neq \emptyset$. On that portion of $\{\rho = 0\}$ the second term in (\ref{4.5}) need not vanish, but it exhibits a sign that fits to that of the first term: 
$$
\lim_{\rho \ra 0} \nn \cdot \aaa' \, \frac{\du^+}{\rho}  = - a_\rho'\big|_{\rho = 0}\, \pa_\rho \du^+\big|_{\rho = 0} = - \pa _\rho \du^+\big|_{\rho = 0} \leq 0\, .
$$
At least in the real-analytic framework, (\ref{4.6}) then implies by Cauchy-Kovalevskaya's theorem $\du^+ \equiv 0$ (end of the heuristics).

Equation (\ref{4.4}) suggests a test function that is constant on $supp\, \du^+$. The sequence
$$
(\psi_n) := \Big(\frac{ n \, \du^+}{n\, \du^+\! + 1}\Big) = \Big(\frac{\du^+}{\du^+\! + 1 /n}\Big)
$$
approximates for $n \ra \infty$ such a test function. Note that with $\du \in H^1_0 (A_R^+)$ we also have $\du^+$ and $\psi_n \in H^1_0 (A_R^+)$ (see, e.g., Gilbarg and Trudinger 1998, p.\ 152f). Testing of (\ref{4.1}) by $\psi_n$ yields
\be \label{4.7}
\ba{l}
\disp 0 = \int_{A_R^+} \Big(\na \du + a_\ze'\, \frac{\du}{\rho}\, \eee_\ze + a_\rho'\, \frac{\du}{\rho}\, \eee_\rho \Big) \cdot \na \psi_n\, \rd \ze \rd \rho \\[3ex]
\disp \quad = \frac{1}{n} \int_{A_R^+} \bigg(
\frac{|\na \du^+|^2}{(\du^+\! + 1/n)^2} + \frac{a_\ze'}{\rho} \, \frac{\du^+\, \pa_\ze \du^+}{(\du^+\! + 1/n)^2} + \Big(\frac{a_\rho' - 1}{\rho} + \frac{1}{\rho}\Big) \frac{\du^+ \,\pa_\rho \du^+}{(\du^+\! + 1/n)^2} \bigg) \rd \ze \rd \rho\, ,
\ea
\ee
where we made use of 
$$
\na \psi_n = \left\{ 
\ba{cc}
\disp \frac{1}{n}\, \frac{\na \du^+}{(\du^+\! + 1/n)^2} & \quad \mbox{on }\, supp\,\du^+\, , \\[3ex]
\disp 0 & \quad \mbox{ else } .
\ea\right.
$$
By rearrangement, use of Cauchy-Schwarz's inequality and some more estimates, (\ref{4.7}) takes the form 
\be \label{4.8}
\ba{l}
\disp \int_{A_R^+} \frac{|\na \du^+|^2}{(\du^+\! + 1/n)^2} \, \rd \ze \rd \rho + \int_{A^+_R} \frac{1}{\rho} \, \frac{\du^+\, \pa_\rho \du^+}{(\du^+\! + 1/n)^2}\, \rd \ze \rd \rho \\[3ex]
\disp \leq \bigg[\bigg(\int_{A^+_R} \Big(\frac{a_\ze'}{\rho}\Big)^2 \rd \ze \rd \rho\bigg)^{\frac{1}{2}} + \bigg( \int_{A^+_R} \Big(\frac{a_\rho' - 1}{\rho} \Big)^2 \rd \ze \rd \rho \bigg)^{\frac{1}{2}} \bigg] \bigg( \int_{A^+_R} \frac{|\na \du^+|^2}{(\du^+\! + 1/n)^2}\,  \rd \ze \rd \rho \bigg)^{\frac{1}{2}} .
\ea
\ee
By (\ref{4.2}), (\ref{3.5}), (\ref{A.7}), and (\ref{B.4}) the $\aaa'$-related terms have the $n$-independent bound: 
$$
\ba{l}
\disp \int_{A^+_R} \Big(\frac{a_\ze'}{\rho}\Big)^2 \rd \ze \rd \rho \leq \int_{A^+_R} \Big[\frac{1}{\rho}\, \sin \Big(\frac{1}{2}\, (u_1 + u_2) - \Om \Big) \Big]^2 \rd \ze \rd \rho \\[3ex]
\disp \quad \leq \int_{A^+_R} \frac{1}{\rho^2}\, \Big(\frac{1}{2}\, (u_1 + u_2) - \Om \Big)^2 \rd \ze \rd \rho
\leq \int_{A^+_R} \Big[\Big(\frac{u_1}{\rho}\Big)^2 + \Big(\frac{u_2}{\rho}\Big)^2 \Big] \rd \ze \rd \rho + 2 \int_{A^+_R} \Big(\frac{\Om}{\rho}\Big)^2 \rd \ze \rd \rho \\[3ex]
\disp \qquad \leq 4 \int_{A^+_R} \big(|\na u_1|^2 + |\na u_2|^2 \big) \rd \ze \rd \rho + 4\, \pi\, K^2 R^2 =: C_{\aaa'}
\ea
$$
and similarly
$$
\ba{l}
\disp \int_{A^+_R} \Big(\frac{a_\rho' - 1}{\rho}\Big)^2 \rd \ze \rd \rho \leq \int_{A^+_R} \frac{1}{\rho^2}\, \bigg[1 - \cos \Big(\frac{1}{2}\, (u_1 + u_2) - \Om \Big) \bigg( \frac{\sin\big((u_1 - u_2)/2\big)}{(u_1 - u_2)/2}\bigg)\bigg]^2 \rd \ze \rd \rho \\[3ex]
\disp \qquad \leq \int_{A^+_R} \frac{1}{\rho^2}\, \bigg[1 - \cos \Big(\frac{1}{2}\, (u_1 + u_2) - \Om \Big) \\[1ex]
\disp \qquad \qquad \qquad \quad + \cos \Big(\frac{1}{2}\,(u_1 + u_2) - \Om\Big)\bigg(1 - \frac{\sin\big((u_1 - u_2)/2\big)}{(u_1 - u_2)/2}\bigg)\bigg]^2 \rd \ze \rd \rho \\[3ex]
\disp \qquad \qquad \leq \int_{A^+_R} \frac{2}{\rho^2}\, \Big[ \Big(\frac{1}{2}\, (u_1 + u_2) - \Om \Big)^2 + \frac{1}{4} (u_1 - u_2)^2 \Big]\rd \ze \rd \rho \, \leq \, 3\, C_{\aaa'} \, .
\ea
$$
The second term on the left-hand side of (\ref{4.8}) is finite by (\ref{4.2}) for any $n \in \nat$,
$$
 \int_{A^+_R} \frac{1}{\rho} \, \frac{\du^+\, \pa_\rho \du^+}{(\du^+\! + 1/n)^2}\, \rd \ze \rd \rho \leq 2 n^2 \int_{A^+_R} | \na \du^+|^2 \rd \ze \rd \rho < \infty\, ,
 $$
 and, moreover, by integrating by parts, turns out to be nonnegative (see Fig.\ 4.1):
\begin{figure}
\begin{center}
\resizebox{0.7\textwidth}{!}{\input{Fig41}}
\caption{Graph of the function $f : x\mapsto \ln (1 + x) - x/(x+1)$, $\, x \geq 0$.}
\end{center}
\end{figure}
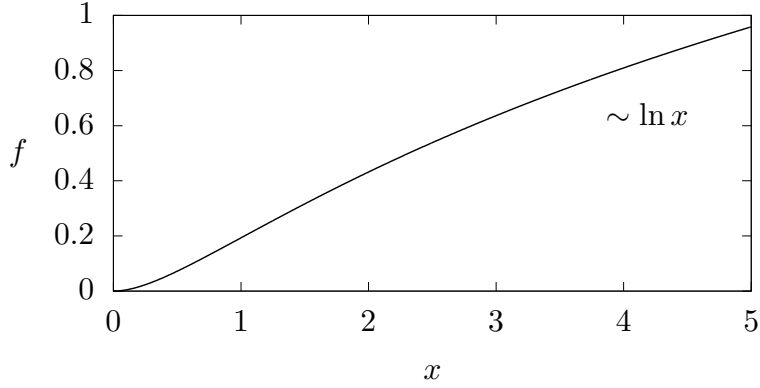
\be \label{4.9}
\ba{l}
\disp \int_{A^+_R} \frac{1}{\rho} \, \frac{\du^+\, \pa_\rho \du^+}{(\du^+\! + 1/n)^2}\, \rd \ze \rd \rho = \int_{A^+_R} \frac{1}{\rho} \,\pa_\rho \Big[ \frac{1/n}{\du^+\! + 1/n} - 1 + \ln \Big(\du^+ \! + \frac{1}{n}\Big) - \ln \frac{1}{n} \Big] \rd \ze \rd \rho \\[3ex]
\disp \qquad = \int_{A^+_R} \frac{1}{\rho^2} \Big[ \ln \big(n\, \du^+\! + 1\big) - \frac{n\, \du^+}{n\, \du^+\! + 1} \Big] \rd \ze \rd \rho \geq 0\, .
\ea
\ee
 We conclude therefore from (\ref{4.8})
$$
\int_{A_R^+} \frac{|\na \du^+|^2}{(\du^+\! + 1/n)^2} \, \rd \ze \rd \rho\leq 3 \sqrt{C_{\aaa'}} \bigg(\int_{A_R^+} \frac{|\na \du^+|^2}{(\du^+\! + 1/n)^2} \, \rd \ze \rd \rho \bigg)^{\frac{1}{2}} ,
$$ 
and, furthermore by (\ref{3.6}):
\be \label{4.10}
\ba{l}
\disp 9\, C_{\aaa'} \geq \int_{A_R^+} \frac{|\na \du^+|^2}{(\du^+\! + 1/n)^2} \, \rd \ze \rd \rho \geq \int_{A^+_R} \big| \na \ln ( 1 + n\, \du^+)\big|^2 \rd \ze \rd \rho \\[3ex]
\disp \qquad\quad \geq \frac{1}{4 R^2} \int_{A^+_R} \big( \ln (1 + n\, \du^+)\big)^2 \rd\ze \rd \rho\, .
\ea
\ee
The monotonically increasing sequence $\big( \ln (1 + n\, \du^+)\big)_{n \in \nat}$ thus has a bounded sequence of integrals, and Levi's theorem implies convergence a.e., which in turn means $\du^+ = 0$ a.e.
This argument works for $-\du^- := - \min \{\du, 0 \}$ as well, which completes the proof.
\qed

In the unbounded case we deal with solutions $u \in H^1_{loc, as} (A_\infty)$ of the problem $P_\Om (A_\infty)$ that satisfy a bound of the form 
\be \label{4.11}
\|\na u \|^2_\beta = \int_{A_\infty} |\na u|^2 \rho^{-\beta}\, \rd\ze \rd \rho \leq C\, , \qquad 0< \beta < 1\, 
\ee
with weight $\rho^{-\beta}$ that cannot be neglected. To incorporate the weight while keeping the divergence-character of the governing equation, we rewrite eq.\ (\ref{2.22}) in the variable $v:= \rho^\beta u$:
\be \label{4.12}
\na \cdot \big(\rho^{-\beta}\, \na v\big) + \pa_\ze \Big[ \frac{1}{\rho} \big(1 - \cos (\rho^{-\beta} v - \Om) \big) \Big] + \pa_\rho \Big[\frac{1}{\rho} \Big( \sin (\rho^{-\beta} v - \Om) -\beta\, \rho^{-\beta} v\Big) \Big] = 0\, .
\ee
From (\ref{4.12}) follows a perturbation equation for $\dv := v_1 - v_2$ analogous to (\ref{4.1}): 
\be \label{4.13}
\na \cdot \Big[ \rho^{-\beta} \Big( \na \dv + \big( a'_\ze\, \eee_\ze + (a_\rho' - \beta)\, \eee_\rho \big) \frac{\dv}{\rho}\Big) \Big] = 0 \quad \mbox{in }\, A_\infty \, , \qquad  \dv\big|_{\pa A_\infty} = 0
\ee
with $\aaa'$ given by (\ref{4.2}) and again considered as a given vector field on $A_\infty$. 
\begin{prop}
 Let $u_1$, $u_2 \in H_{loc, as}^1 (A_\infty)$ be weak solutions of problem $P_\Om (A_\infty)$ with vanishing trace on $S_1$ and bound (\ref{4.11}); let $\Om \in L^\infty (A_\infty)$ satisfy the bound (\ref{2.32}). Then
 $$
 u_1 = u_2 \qquad a.\, e.\ \mbox{in } \, A_\infty\, .
 $$
 \end{prop}
\textsc{Proof:} Based on eq.\ (\ref{4.13}) the proof proceeds quite analogously to that of proposition 4.1; we skip thus the details. With the sequence of test functions 
$$
(\psi_n) := \Big(\frac{\dv^+}{\dv^+\! + 1 /n}\Big)
$$
one obtains instead of (\ref{4.7})
\be \label{4.14}
0 = \int_{A_R^+} \bigg(
\frac{|\na \dv^+|^2}{(\dv^+\! + 1/n)^2} + \frac{a_\ze'}{\rho} \, \frac{\dv^+\, \pa_\ze \dv^+}{(\dv^+\! + 1/n)^2} + \Big(\frac{a_\rho' - 1}{\rho} + \frac{1 -\beta}{\rho}\Big) \frac{\dv^+ \,\pa_\rho \dv^+}{(\dv^+\! + 1/n)^2} \bigg) \rd \mu_\beta\, .
\ee
The $\aaa'$-related estimates are now done by (\ref{A.8}) and (\ref{B.5}):
$$
\ba{l}
\disp \int_{A^+_\infty} \Big(\frac{a_\ze'}{\rho}\Big)^2 \rd \mu_\beta 
\leq \int_{A^+_\infty} \Big[\Big(\frac{u_1}{\rho}\Big)^2 + \Big(\frac{u_2}{\rho}\Big)^2 \Big] \rd \mu_\beta + 2 \int_{A^+_R} \Big(\frac{\Om}{\rho}\Big)^2 \rd \mu_\beta \\[3ex]
\disp \qquad \leq \Big(\frac{2}{1 + \beta}\Big)^2 \big(\|\na u_1\|^2_\beta + \|\na u_2\|^2_\beta \big) + \frac{4\, \pi\, K^2}{\beta (1- \beta)} =: \Ct_{\widetilde{\aaa}'}\, ,
\ea
$$
and analogously for $a_\rho'$. The fourth term in (\ref{4.14}) is still nonnegative,
$$
(1 - \beta) \int_{A^+_\infty } \frac{1}{\rho} \, \frac{\dv^+\, \pa_\rho \dv^+}{(\dv^+\! + 1/n)^2}\, \rd \mu_\beta = (1 - \beta)(1 + \beta) 
\int_{A^+_\infty} \frac{1}{\rho^2}\, f(n\, \dv^+)\, \rd \mu_\beta \geq 0\, \, ,
$$
and (\ref{4.10}) takes the form
$$
\ba{l}
\disp 9\, \Ct_{\widetilde{\aaa}'} \geq \int_{A_\infty^+} \frac{|\na \dv^+|^2}{(\dv^+\! + 1/n)^2} \, \rd \mu_\beta \geq \int_{A^+_R} \big| \na \ln ( 1 + n\, \dv^+)\big|^2 \rd \mu_\beta \\[3ex]
\disp \qquad\quad \geq \frac{(1+ \beta)^2}{4 R^2} \int_{A^+_R} \big( \ln (1 + n\, \dv^+)\big)^2 \rd\mu_\beta
\ea
$$
with arbitrary $R > 1$. By $n \ra \infty$ we can thus again conclude that $\dv^+ = 0$ and hence $\du^+ = 0$ a.e.\ in $A_R^+$ for any $R > 1$. 
\qed
\sect{The linearized problem}
This section is devoted to the solution of eq.\ (\ref{2.29}) for given right-hand side and given coefficients $a_\ze[\cdot ]$ and $a_\rho[\cdot ]$ of type (\ref{2.24}) with measurable but otherwise not further specified argument. Let us start with the observation that by antisymmetry it is sufficient to consider (\ref{2.29}) in $A_R^+ = A_R \cap \{\rho >0\}$. Given $\Om$, $u \in H^1_0 (A_R^+)$ is called a weak solution of the problem $P_\Om  (A_R^+)$ if $u$ satisfies (\ref{2.29}), with $A_R$ replaced by $A_R^+$, for any test function $\vii \in C_0^\infty (A_R^+)$.  The equivalence with $P_\Om (A_R)$ is clear when observing the one-to-one correspondence between elements of $H_0^1 (A_R^+)$ and $H_{0, as}^1 (A_R)$ (by antisymmetric continuation and restriction, respectively).

The general solution strategy for equations of type
\be \label{5.1}
\int_{A_R^+} \na u \cdot \na \vii\, \rd \ze \rd \rho + \int_{A_R^+} \frac{u}{\rho}\, \aaa \cdot \na \vii\, \rd \ze \rd \rho =
\int_{A_R^+} \frac{\Om}{\rho}\, \aaa \cdot \na \vii\, \rd \ze \rd \rho \, ,
\ee
where $\aaa$ abbreviates $\aaa[f]$ with some measurable function $f : A_R^+ \ra \real$, is to apply a Lax-Milgram type criterion, which requires, in particular, coerciveness of the bilinear form on the left-hand side of (\ref{5.1}). In this situation as sharp as possible bounds on $\aaa$ play a crucial role. A first step in this direction is to introduce the variable $v := \rho^\al u$ that allows us to shift $a_\rho$ by $\al$ and thus to take advantage of the asymmetry between upper and lower bounds on $a_\rho$ (see Fig.\ 2.2). In fact, expressing (\ref{5.1}) by $v$ one obtains
\be \label{5.2}
\int_{A_R^+} \na v \cdot \na \vii\, \rd \mu_\al  + \int_{A_R^+} \frac{v}{\rho}\, \aaa^\al\! \cdot \na \vii\, \rd \mu_\al =
\int_{A_R^+} \frac{\Om}{\rho}\, \aaa \cdot \na \vii\, \rd \ze \rd \rho
\ee
with $\rd \mu_{\al} = \rho^{-\al} \rd \ze \rd \rho$ and $\aaa^\al := (a_\ze, a_\rho - \al)$. Proper choice of $\al$ ``centers'' $a_\rho$ with the effect of an improved bound: 
\be \label{5.3}
\|a_\rho^\al\|_\infty = \| a_\rho - \al\|_\infty \leq 0.61 \qquad \mbox{ for } \, \al = 0.39\, .
\ee
Obviously the bound 
\be \label{5.4}
\|a_\ze\|_\infty \leq 0.73
\ee
on the antisymmetric function $a_\ze$ cannot be improved this way.

In a second step we want to take advantage of the improved constant in the Hardy inequality (\ref{3.5}) when weighted by a factor $\rho^{-\ga}$, $\ga >0$. Note that by the box criterion 3.2 this constant is optimal in $A_R^+$. The Lax-Milgram criterion must now be adapted to this weighted situation. A suitable version is given in the following proposition whose proof is deferred to appendix C.
\begin{prop}[generalized Lax-Milgram criterion]
Let $\BC$ and $\BtC$ be reflexive Banach spaces and $B : \BC \ti \BtC \ra \real$ be a continuous bilinear form, i.e.\ for some $K> 0$ holds
\be \label{5.5}
B(u, \ut) \leq K\, \|u\|_\BC \|\ut\|_{\BtC} \qquad \mbox{for all }\, u \in \BC\, ,\; \ut \in \BtC\, .
\ee
Let, furthermore, $\BtC'$ be the dual space of $\BtC$ with norm
\be \label{5.6}
\|w \|_{\BtC'} := \sup_{0\neq \ut \in \BtC} \frac{\langle w, \ut\rangle}{\|\ut\|_\BtC}\, ,
\ee
where $\langle\cdot , \cdot \rangle$ denotes the dual pairing between $\BtC$ and $\BtC'$. Let, finally, $c> 0$ and $\ct >0$ be constants such that 
 \be \label{5.7}
 \sup_{0\neq \ut \in \BtC} \frac{B(u , \ut)}{\|\ut\|_\BtC} \geq c \, \|u\|_\BC \qquad \mbox{for all }\, u \in \BC\, ,
 \ee
 \be \label{5.8}
 \sup_{0\neq u \in \BC} \frac{B(u , \ut)}{\|u\|_\BC} \geq \ct \, \|\ut\|_\BtC \qquad \mbox{for all }\, \ut \in \BtC\, .
 \ee
 Then, the equation 
 \be \label{5.9}
 B(u, \ut) = \langle w, \ut \rangle \qquad \mbox{for all }\, \ut \in \BtC
 \ee
 with $w \in \BtC'$ has a unique solution $u \in \BC$ with bound 
 \be \label{5.10}
 \|u \|_\BC \leq \frac{1}{c}\, \|w \|_{\BtC'}\, .
 \ee
\end{prop}
Obviously, conditions (\ref{5.7}) and (\ref{5.8}) replace the coercivity condition $B(u,u) \geq c\, \|u\|^2$ of the ordinary Lax-Milgram criterion for a single (Hilbert) space.

To apply the criterion we set
$$
\BC := \HC^{(d)}_{\al + \ga} (A_R^+ )\, ,\qquad \BtC := \HC^{(e)}_{\al - \ga} (A_R^+ )\, ,
$$
with $\HC^{(c)}_\beta (G)$ defined in (\ref{3.2}) and with parameters 
\be \label{5.11}
\al \geq 0\, , \quad \ga \geq \al\, ,\quad 0 < d\leq 1\, ,\quad e\geq 1
\ee
yet to be fixed. The bilinear form $B$ is given by 
\be \label{5.12}
B[v, \psi] := B_0 [v, \psi] + B_1 [v, \psi] := \int_{A_R^+} \na v \cdot \na \psi \, \rd \mu_\al + \int_{A_R^+} \frac{v}{\rho}\, \aaa^\al \cdot \na \psi \, \rd \mu_\al
\ee
with 
\be \label{5.12a}
v \in \HC^{(d)}_{\al + \ga} (A_R^+ ) \subset \HC_{2 \al} (A_R^+) \subset H_0^1 (A_R^+) \, ,\qquad \psi \in \HC^{(e)}_{\al - \ga} (A_R^+ )
\supset H_0^1 (A_R^+)\, .
\ee
Condition (\ref{5.5}) may easily be checked by Cauchy-Schwarz's inequality and (\ref{3.5}):
\be \label{5.13}
\ba{l}
\disp |B[v, \psi]| \leq \int_{A_R^+} |\na v|\, \rho^{-(\al + \ga)/2}\, |\na \psi|\, \rho^{-(\al - \ga)/2}\, \rd \ze \rd \rho \\[2ex]
\disp \qquad \qquad \qquad + \int_{A_R^+} |v/ \rho|\, \rho^{-(\al + \ga)/2}\,\sqrt{2}\, |\na \psi|\, \rho^{-(\al - \ga)/2}\, \rd \ze \rd \rho \\[3ex]
\disp \qquad \qquad \leq \| \na v\|_{\al + \ga}\, \|\na \psi\|_{\al - \ga} + \sqrt{2}\, \|v/\rho\|_{\al + \ga }\, \| \na \psi \|_{\al - \ga } \\[2ex]
\disp \qquad \qquad \qquad \leq \Big(\frac{1}{d} + \frac{2 \sqrt{2}}{1 + \al + \ga} \Big) \|\na_{ d}\, v\|_{\al + \ga}\, \|\na_e \psi\|_{\al - \ga}\, .
\ea
\ee
As to conditions (\ref{5.7}), (\ref{5.8}) note that there is no useful information about $\aaa^\al$ other than the bounds (\ref{5.3}), (\ref{5.4}). In (\ref{5.12}), $B_1$ is thus considered as a perturbation of $B_0$. The optimal ($\triangleq$ largest possible) constants in (\ref{5.7}), (\ref{5.8}) with respect to $B_0$ then are determined by the following min-max problems:
\be \label{5.14}
\inf_{0 \neq v \in \HC_{\al + \ga}^{(d)} (A_R^+)}\, \sup_{0 \neq \psi \in \HC_{\al - \ga}^{(e)} (A_R^+)} \frac{ \int_{A_R^+} \na  v \cdot \na \psi \, \rd \mu_\al}{\| \na_d \, v\|_{\al + \ga}\, \|\na_e \, \psi \|_{\al - \ga}} = : C_c\, , 
\ee
\be \label{5.15}
\inf_{0 \neq \psi \in \HC_{\al - \ga}^{(e)} (A_R^+)}\, \sup_{0 \neq v \in \HC_{\al + \ga}^{(d)} (A_R^+)} \frac{ \int_{A_R^+} \na  v \cdot \na \psi \, \rd \mu_\al}{\| \na_d \, v\|_{\al + \ga}\, \|\na_e \, \psi \|_{\al - \ga}} = : \Ct_c\, . 
\ee
On the other side, $B_1$ may be estimated similarly as in (\ref{5.13}): 
\be \label{5.16}
\ba{l}
\disp \bigg|\int_{A_R^+} \frac{v}{\rho}\, \aaa^\al\! \cdot \na \psi \, \rd \mu_\al \bigg| \leq \|a^\al_\rho\|_\infty  \int_{A_R^+} \Big|\frac{v}{\rho} \Big|\, \rho^{-(\al+ \ga)/2}\, \big( |e\, \pa_\ze \psi| + |\pa_\rho \psi| \big) \rho^{-(\al -\ga)/2}\, \rd \ze \rd \rho \\[3ex]
\disp \qquad \leq \|a^\al_\rho\|_\infty\, \frac{2}{1 + \al + \ga}\, \|\pa_\rho v\|_{\al + \ga}\, \sqrt{2}\, \|\na_e \psi \|_{\al - \ga} \\[3ex]
\disp \qquad \qquad \leq \frac{2 \sqrt{2}\, \| a_\rho^\al\|_\infty }{1 + \al + \ga} \, \|\na_d \, v\|_{\al + \ga}\, \|\na_e \psi\|_{\al - \ga }\, ,
\ea
\ee
where we have set 
\be \label{5.17}
e := \frac{\|a_\ze\|_\infty}{\|a_\rho^\al\|_\infty}\, .
\ee
It is this estimate that profits by the introduction of the parameters $d$ and $e$; note that any $d> 0$ is allowed in (\ref{5.16}). 
Combining (\ref{5.12}), (\ref{5.14}), and (\ref{5.16}) one obtains for any $v \in \HC^{(d)}_{\al + \ga} (A_R^+)$:
\be \label{5.17a}
\ba{l}
\disp \sup_{0 \neq \psi \in \HC_{\al - \ga}^{(e)} (A_R^+)} \frac{B[v, \psi]}{\|\na_e \psi\|_{\al - \ga}}\ \geq \sup_{0 \neq \psi \in \HC_{\al - \ga}^{(e)} (A_R^+)} \bigg\{ \frac{B_0[v, \psi]}{\|\na_e \psi\|_{\al - \ga}} - \frac{2 \sqrt{2}\, \| a_\rho^\al\|_\infty }{1 + \al + \ga} \, \|\na_d \, v\|_{\al + \ga} \bigg\} \\[4ex]
\disp \qquad \qquad \qquad \geq \bigg( C_c - \frac{2 \sqrt{2}\, \| a_\rho^\al\|_\infty }{1 + \al + \ga} \bigg) \|\na_d \, v\|_{\al + \ga} \, ,
\ea
\ee
which is condition (\ref{5.7}), provided that 
\be \label{5.18}
C_c >  \frac{2 \sqrt{2}\, \| a_\rho^\al\|_\infty }{1 + \al + \ga} \, .
\ee
An analogous estimate yields
$$
\sup_{0 \neq v \in \HC_{\al + \ga}^{(d)} (A_R^+)} \frac{B[v, \psi]}{\|\na_d\, v\|_{\al + \ga}}\, \geq \bigg( \Ct_c - \frac{2 \sqrt{2}\, \| a_\rho^\al\|_\infty }{1 + \al + \ga} \bigg) \|\na_e \psi\|_{\al - \ga} \, ,
$$
which is (\ref{5.8}), provided that 
\be \label{5.19}
\Ct_c >  \frac{2 \sqrt{2}\, \| a_\rho^\al\|_\infty }{1 + \al + \ga} \, .
\ee
By (\ref{5.18}) and (\ref{5.19}) the generalized Lax-Milgram criterion 5.1 provides a solution $v \in \HC^{(d)}_{\al + \ga} (A_R^+)$ of eq.\ (\ref{5.2}), which by (\ref{5.6}) and (\ref{5.10}) satisfies the bound
$$
\ba{l}
\disp \|\na v\|_{\al + \ga} \leq \frac{1}{d}\, \|\na_d\, v\|_{\al + \ga} \leq \frac{1}{d}\, \frac{1}{\De} \sup_{0 \neq \psi \in \HC_{\al - \ga}^{(e)} (A_R^+)} \frac{\int_{A_R^+} \frac{\Om}{\rho}\, \aaa^\al\! \cdot \na \psi\, \rd \ze \rd \rho}{\|\na_e \psi\|_{\al - \ga}}\\[4ex]
\disp\qquad \qquad\ \leq \frac{1}{d}\,\frac{1}{\De}\, \Big\|\frac{\Om}{\rho}\,\aaa^\al \Big\|_{\ga - \al} \leq \frac{\sqrt{2}}{\De\, d} \, \|\Om\|_{2 + \ga - \al}
\ea
$$
with 
\be \label{5.20}
\De := C_c -  \frac{2 \sqrt{2}\, \| a_\rho^\al\|_\infty }{1 + \al + \ga} \, .
\ee
Having in mind problem $P_\Om (A_R^+)$ we call a function $v \in \HC_{2 \al} (A_R^+)$ solution of eq.\ (\ref{5.2}) with given measurable functions $a_\ze$, $a_\rho^\al$, $\Om : A_R^+ \ra \real$, if $v$ satisfies (\ref{5.2}) for any $\psi \in C_0^\infty (A_R^+)$. By (\ref{5.12a}), proposition 5.1 obviously provides a solution of this type. We summarize these results in the following proposition.
\begin{prop}
 Let $a_\ze$, $a_\rho^\al$ be bounded, measurable functions on $A_R^+$, satisfying (\ref{5.3}) and (\ref{5.4}),  
 and let $\al$, $\ga$, $d$, $e$, and $R$ be such that conditions (\ref{5.11}), (\ref{5.18}), and (\ref{5.19}) are satisfied, where $C_c$ and $\Ct_c$ are given by (\ref{5.14}) and (\ref{5.15}), respectively.
 Let, furthermore, $\Om : A_R^+ \ra \real$ be such that $\|\Om\|_{2 + \ga - \al}
  < \infty$. Then, eq.\ (\ref{5.2}) has a unique solution $v \in \HC_{\al + \ga}^{(d)} (A_R^+) \subset \HC_{2 \al} (A_R^+)$ satisfying the bound 
  \be \label{5.21}
  \|\na_d\, v\|_{\al + \ga } \leq \frac{\sqrt{2}}{\De\, d} \, \|\Om \|_{2 + \ga - \al}
  \ee
  with $\De$ given by (\ref{5.20}).
 \end{prop}
By lemma 3.2 the grad-norms of $u$ and $v$ are related by 
\be \label{5.21a}
\frac{\Ct_{\ga + \al}}{C_{\ga - \al}}\, \|\na v\|_{\al + \ga} \leq \|\na u\|_{\ga - \al} \leq \frac{C_{\ga + \al}}{\Ct_{\ga - \al}}\, \|\na v\|_{\al + \ga}\, .
\ee
The following corollary is thus an immediate consequence of the preceeding proposition.
\begin{cor}
 Under the conditions of proposition 5.2, eq.\ (\ref{5.1}) has a unique solution $u \in H_0^1 (A_R^+)$ with bound 
  \be \label{5.22}
  \|\na u\|_{\ga - \al } \leq \Ch \, \|\Om \|_{2 + \ga - \al}\, ,
  \ee
  where the constant $\Ch$ depends on $\al$, $\ga$, $d$, $e$, and $R$.
\end{cor}
\noi {\em Remark}: Inequality (\ref{3.4}) offers yet another possibility to improve the Hardy constant: higher $L^p$-spaces. The generalized Lax-Milgram criterion then is applied to conjugate $L^p$-spaces (see Kaiser 2010). However, the determination of the constants corresponding to (\ref{5.14}) and (\ref{5.15}) requires now the solution of systems of {\em nonlinear} partial differential equations as opposed to the {\em linear} systems associated to (\ref{5.14}) and (\ref{5.15}), which are considered and solved in section 6.  
\sect{A minimum-maximum problem}
In this section the min-max problems (\ref{5.14}) and (\ref{5.15}) are studied in some detail with the goal to establish parameter sets $\{\al, \ga, d, e\}$ that satisfy the sufficient conditions (\ref{5.18}) and (\ref{5.19}) for solvability of the linearized equation (\ref{5.2}). 
Proposition 6.1 allows us to relate the problem posed in $A_\infty^+$ to one posed in some finite rectangle $Q \subset H$. The rectangular geometry then allows the reduction of the two-dimensional problem to a sequence of one-dimensional ones (proposition 6.3).
We determine the Euler-Lagrange equations for the two-dimensional as well as for the one-dimensional problems and find a common lower bound on $C_c$ and $\Ct_c$ in terms of the minimum eigenvalue associated to these equations (proposition 6.2). 
The analytic solution of even the one-dimensional problems succeeds only in a special (the $
x$-independent) case either directly (see appendix D) or via Euler-Lagrange equations (see subsection 6.3). Finally, based on the analytic solution, we derive explicit rigorous lower bounds on $C_c$ and $\Ct_c$ (proposition 6.4), which are corroborated by the numerical solution of the general one-dimensional equations (see Appendix E).

Let us start with some simple observations. Using the notation 
\be \label{6.1}
\ba{l}
\disp \FC [v, \psi] := \FC [v, \psi ; G] := \frac{\int_G \na v \cdot \na \psi\, \rd \mu_\al}{\|\na_d\, v\|_{\al + \ga } \|\na_e \psi\|_{\al - \ga}} \\[3ex]
\disp  = \frac{\disp \int_{G} \big(\pa_x v\, \pa_x \vii + \pa_y v\, \pa_y \vii\big)\, y^{-\al}\, \rd x \rd y}{\disp \bigg( \int_{G} \big( (d\, \pa_x  v)^2
+ (\pa_y v)^2\big) y^{-(\al + \ga)} \rd x \rd y \bigg)^{1/2} \bigg( \int_{G} \big((e\, \pa_x \vii)^2 + (\pa_y \vii)^2\big) y^{-(\al - \ga)} \rd x \rd y \bigg)^{1/2}}
\ea
\ee
and 
\be \label{6.2}
C_c := C_c (G) := C_c (G; \al, \ga, d, e) := \inf_{0 \neq v \in \HC_{\al+ \ga}^{(d)} (G)}\ \sup_{0 \neq \psi \in \HC_{\al - \ga}^{(e)} (G)}\ \FC [v, \psi; G]\, ,
\ee
\be \label{6.2a}
\Ct_c := \Ct_c (G) := \Ct_c (G; \al, \ga, d, e) := \inf_{0 \neq \psi \in \HC_{\al - \ga}^{(e)} (G)}\ \sup_{0 \neq v \in \HC_{\al+ \ga}^{(d)} (G)}\ \FC [v, \psi; G]\, ,
\ee
where $G \subset H$ and $\al$, $\ga$, $d$, and $e$ satisfy (\ref{5.11}), one obtains by Cauchy-Schwarz's inequality:
\be \label{6.3}
\FC [v, \psi] \leq \frac{\|\na_{1/e}\, v \|_{\al + \ga} \, \|\na_e \psi\|_{\al - \ga}}{\| \na_d\, v \|_{\al + \ga}\, \| \na_e \psi\|_{\al - \ga}} = \frac{\|\na_{1/e}\, v\|_{\al + \ga}}{\|\na_d\, v\|_{\al + \ga}} \, .
\ee
Let $I_x \ti I_y \subset G$ some rectangle, $\beta \in \real$, $0 \not\equiv \chi_x \in C_0^\infty (I_x)$, and $(\chi_n) \subset C_0^\infty (I_y)$ a sequence of test functions  with
$$
\int_{I_y} \chi_n^2\, y^{- \beta} \rd y = 1\ , \quad \lim_{n \ra \infty} \int_{I_y} {\chi'_n}^2\, y^{- \beta} \rd y = \infty\, .
$$
Setting $v_n (x,y) := \chi_x (x)\, \chi_n (y) \in C_0^\infty (G)$ and $\beta:= \al + \ga$ we find 
\be \label{6.4}
\lim_{n \ra \infty} \frac{\|\na_{1 /e} \, v_n \|_{\al + \ga}}{\|\na_d\, v_n\|_{\al + \ga}} = 1\, ,
\ee
which implies $C_c \leq 1$. The lower bound $C_c \geq 0$ follows for given $v$ and $\psi$ by proper choice of the relative sign. 

In the case $\ga = 0$ by the choice $\psi := v$ we obtain the lower bound
$$
\ba{c}
\disp C^2_c \geq \inf_{0 \neq v \in \HC_\al (G)} \frac{\Big(\int_G \big((\pa_x v)^2 + (\pa_y v)^2\big) \rd \mu_\al \Big)^2}{
\int_G \big((d\, \pa_x v)^2 + (\pa_y v)^2\big) \rd \mu_\al\, \int_G \big((e\, \pa_x v)^2 + (\pa_y v)^2\big) \rd \mu_\al} \\[4ex]
\disp \geq \, \inf_{0\leq s < \infty}\, \frac{(1 + s)^2}{(d^2 + s)(e^2 + s)}\, ,
\ea
$$
where $s$ denotes the ratio $\int_G (\pa_y v)^2 \rd \mu_\al /\int_G (\pa_x v)^2 \rd \mu_\al$. By the condition 
\be \label{6.4a}
e^2 d^2 < \frac{d^2 + e^2}{2} < 1
\ee
it is easily checked that $C_c \geq 1$ and hence $C_c(G; \al, 0, d,e) = 1$.

Concerning $G$ the following scaling property of $\FC$ plays an important role. Let $G_\la := \{ (\la x, \la y) \in H : (x,y) \in G\}$ and 
\be \label{6.5}
S_\la : \HC_\beta (G) \ra \HC_\beta (G_\la ) \, , \quad f \mapsto f_\la := f(\la^{-1}\, \cdot\, )\, , \qquad \la >0\, .
\ee
Obviously we then have 
\be \label{6.6}
\FC [v, \psi; G] = \FC [v_\la , \psi_\la; G_\la]
\ee
and hence $C_c (G) = C_c(G_\la)$. In particular, applying $S_\la$ on $B_R^+$ and $A^+_{r_0, \infty} := \{(x,y) \in H : x^2+ y^2 > r_0^2 \}$ we find in the limits $\la \ra \infty$ and $\la \ra 0$, respectively,
\be \label{6.7}
C_c (B_R^+) = C_c (H) = C_c (A_{r_0, \infty}^+)
\ee
for any $R > 0$, $r_0 \geq 0$. Similary, for any rectangle $Q:= (-a, a) \ti (0, b) \subset H$ holds
\be \label{6.8}
C_c (Q) = C_c (H)\, .
\ee
The same or similar arguments apply to $\Ct_c$ with identical results.
We summarize these results in the following proposition.
\begin{prop}
 Let $G \subset H$ and let $\al$, $\ga$, $d$, and $e$ satisfy the conditions (\ref{5.11}) and (\ref{6.4a}). Then, the variational constants $C_c$ and $\Ct_c$ given by (\ref{6.1})--(\ref{6.2a}) satisfy the bounds
 \be \label{6.9}
 0 \leq C_c \leq 1\, , \qquad 0 \leq \Ct_c \leq 1\, .
 \ee
 For $\ga = 0$ holds 
 \be \label{6.10}
 C_c (G; \al, 0, d,e) = 1 = \Ct_c (G; \al, 0, d,e) \, ,
 \ee
 and $C_c (H)$ and $\Ct_c (H)$ coincide with $C_c$ and $\Ct_c$, respectively, for some (bounded and unbounded) standard domains such as $Q$, $B_R^+$, or $A_{r_0, \infty}^+$.
\end{prop}
\subsection{Euler-Lagrange equations}
Critical points of the variational expression (\ref{6.1}) are (weak) solutions of the associated Euler-Lagrange equations. We derive
these equations in two steps. First we fix $v \in \HC_{\al + \ga} (G)$ and vary the functional 
\be \label{6.1.1}
\int_G \na v \cdot \na \psi \ y^{-\al}\, \rd x \rd y
\ee
with respect to $\psi$ under the constraint $\|\na_e \psi \|_{\al - \ga} = \|\na_d \, v\|_{\al +\ga}$. Introducing the Lagrange parameter $\la \in \real$, this is equivalent to the variation of the extended functional 
$$
\int_G \na v \cdot \na \psi \ y^{-\al}\, \rd x \rd y -\la\, \big(\|\na_e \psi \|^2_{\al - \ga} - \|\na_d \, v\|^2_{\al +\ga}\big)
$$
with respect to $\psi$ and $\la$ (see, e.g., Courant and Hilbert 1953, vol.\ I, p.\ 216 ff). 
Setting this variation to zero one obtains 
\be \label{6.1.2}
\left.
\begin{array}{c}
\na \cdot (y^{-\al} \na v) - 2\,\la \na_e \cdot (y^{-(\al - \ga )} \na_e \psi ) =0\, , \\[2ex]
\|\na_e \psi \|_{\al - \ga} = \|\na_d \, v\|_{\al +\ga}\, . 
\end{array} \right\}
\ee
In the second step we vary (\ref{6.1.1}) under the differential constraint (\ref{6.1.2})$_1$ and $\|\na_d\, v \|_{\al + \ga} = 1$ with respect to $v$ and $\psi$. Introducing the Lagrange multipliers $p$, where $y^{-(\al -\ga)/2}\, p \in H^1(G)$, and $\mu \in \real$, this is equivalent to the variation of the extended functional
\be \label{6.1.2a}
\ba{c}
\disp \int_G \na v \cdot \na \psi \ y^{-\al}\, \rd x \rd y
+ \int_G \na p \cdot \na v\ y^{-\al} - 2\, \la \na_e\, p \cdot \na_e \psi\ y^{-(\al - \ga)}\, \rd x \rd y \\[2ex]
\disp -\mu\, \big(\|\na_d \, v\|^2_{\al +\ga} - 1\big)
\ea
\ee
with respect to $v$, $\psi$, $p$, and $\mu$. Setting again the variation to zero one obtains together with (\ref{6.1.2})$_2$:
\be \label{6.1.3}
\left.
\begin{array}{c}
\na \cdot (y^{-\al} \na \psi) - 2\,\mu \na_d \cdot (y^{-(\al + \ga )} \na_d \, v ) + \na \cdot (y^{-\al}\, \na p) =0\, , \\[2ex]
\na \cdot (y^{-\al} \na v) - 2\,\la \na_e \cdot (y^{-(\al - \ga )} \na_e\, p ) =0\, , \\[2ex]
\na \cdot (y^{-\al} \na v) - 2\,\la \na_e \cdot (y^{-(\al - \ga )} \na_e \psi ) =0\, , \\[2ex]
y^{-(\al -\ga)/2}\, p\big|_{\pa G} = 0\, , \\[2ex]
\|\na_e \psi \|_{\al - \ga} = \|\na_d \, v\|_{\al +\ga} = 1\, . 
\end{array} \right\}
\ee
This system is yet somewhat redundant. By (\ref{6.1.3})$_4$ and (\ref{3.11}) one finds $p\in \HC_{\al -\ga} (G)$; setting thus $p:= \psi$ makes eqs.\ (\ref{6.1.3})$_{2,3}$ equivalent. Furthermore, multiplying (\ref{6.1.3})$_1$ by $v$ and (\ref{6.1.3})$_3$ by $\psi$ and integrating over $G$ one finds by (\ref{6.1.3})$_5$ that $\mu = 2\, \la$. Finally, when discarding conditions (\ref{6.1.3})$_5$, we end up with the eigenvalue problem 
\be \label{6.1.5}
\left.
\begin{array}{c}
\na \cdot (y^{-\al} \na \psi) - \mu \na_d \cdot (y^{-(\al + \ga )} \na_d\, v ) =0\, , \\[2ex]
\na \cdot (y^{-\al} \na v) - \mu \na_e \cdot (y^{-(\al - \ga )} \na_e \psi ) =0\, ,
\end{array} \right\}
\end{equation}
whose eigenvalue $\mu$ is related to the critical point $(v, \psi)$ by 
\be \label{6.1.4}
\mu = \frac{\int_G \na v \cdot \na \psi\ y^{-\al}\, \rd x \rd y}{\|\na_d\, v\|_{\al + \ga } \|\na_e \psi\|_{\al - \ga}}\ . 
\ee
Note that stationarity of the functional (\ref{6.1}) is only a necessary condition for min-max-points. Comparing (\ref{6.1.4}) with (\ref{5.14}) we can thus only conclude that  
$$
C_c \geq \mu_{min}
$$
where
\be \label{6.1.5a}
\mu_{min} := \inf \, \{\, \mu : \mu \mbox{ eigenvalue of (\ref{6.1.5})}\, \}\, .
\ee

Repeating this procedure with $v$ and $\psi$ interchanged yields instead of (\ref{6.1.2})
\begin{equation}\label{6.1.6}
\left.
\begin{array}{c}
\na \cdot (y^{-\al} \na \psi) - 2\,\la \na_d \cdot (y^{-(\al + \ga )} \na_d\, v ) =0\, , \\[2ex]
\|\na_d\, v \|_{\al + \ga} = \|\na_e \psi\|_{\al -\ga}\, . 
\end{array} \right\}
\end{equation}
Therefore, in the second step we vary instead of (\ref{6.1.2a}) the functional
\be \label{6.1.7}
\ba{c}
\disp \int_G \na v \cdot \na \psi \ y^{-\al}\, \rd x \rd y
+ \int_G \na p \cdot \na \psi\ y^{-\al} - 2\, \la \na_d\, p \cdot \na_d \, v\ y^{-(\al + \ga)}\, \rd x \rd y \\[2ex]
\disp -\mu\, \big(\|\na_e \psi\|^2_{\al - \ga} - 1\big)
\ea
\ee
with respect to $v$, $\psi$, $p$, and $\mu$. Instead of (\ref{6.1.3}) one obtains  
\be \label{6.1.8}
\left.
\begin{array}{c}
\na \cdot (y^{-\al} \na \psi) - 2\,\la \na_d \cdot (y^{-(\al + \ga )} \na_d \, p ) =0\, , \\[2ex]
\na \cdot (y^{-\al} \na v) - 2\,\mu \na_e \cdot (y^{-(\al - \ga )} \na_e \psi ) - \na \cdot (y^{-\al} \na p) =0\, , \\[2ex]
\na \cdot (y^{-\al} \na \psi) - 2\,\la \na_d \cdot (y^{-(\al + \ga )} \na_d\, v ) =0\, , \\[2ex]
y^{-(\al +\ga)/2}\, p\big|_{\pa G} = 0\, , \\[2ex]
\|\na_d\, v \|_{\al + \ga} = \|\na_e\psi\|_{\al - \ga} = 1\, . 
\end{array} \right\}
\ee
Here, $p \in \HC_{\al + \ga} (G)$ and setting $p := v$ makes (\ref{6.1.8})$_{1,3}$ equivalent. As before we find $\mu = 2\,\la$ with the result that $v$, $\psi$, and $\mu$ satisfy precisely the eqs.\ (\ref{6.1.5}). We summarize these findings in the following proposition.
\begin{prop}
 Let (\ref{6.1.2a}) and (\ref{6.1.7}) be the extended functionals associated to the variational problems (\ref{5.14}) and (\ref{5.15}), respectively. Then both functionals share the same set of critical points
 and their min-max-related critical values $C_c$ and $\Ct_c$ satisfy the common lower bound 
 \be \label{6.1.9}
 C_c \geq \mu_{min}\, ,\qquad \Ct_c \geq \mu_{min}
 \ee
 with $\mu_{min}$ given by (\ref{6.1.5a}). 
 In the case of a unique critical value $\mu$ we have 
 \be \label{6.1.10}
 C_c = \mu = \Ct_c\, .
 \ee
\end{prop}
\subsection{One-dimensional min-max problems}
Let us now specialize to the rectangular domain $Q = (-a, a)\ti (0, b)$ (see Fig.\ 6.1), which suggests the following product ansatz for the variational functions $v$ and $\psi$:
\be \label{6.1.11}
v(x,y) = s_n (x)\, v_n (y)\, \qquad \psi(x,y) = s_n (x)\, \psi_n (y)
\ee
with
\be \label{6.1.12}
s_n (x) := \frac{1}{\sqrt{a}} \sin k_n (x+ a)\, ,\quad k_n := n\, \frac{\pi}{2 a}\, , \quad n \in \nat \, .
\ee
\begin{figure}
\begin{center}
\includegraphics[width=0.4\textwidth]{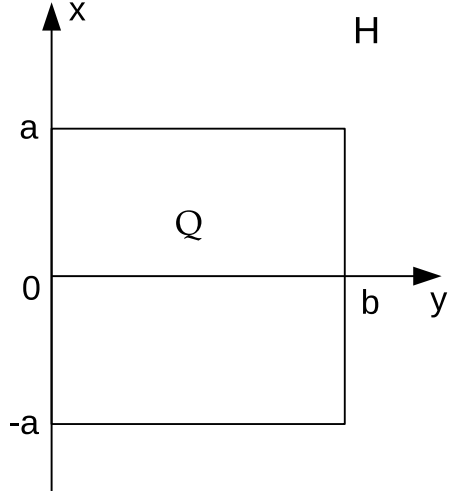}
\caption{The rectangle $Q$ in the half space $H$.}
\end{center}
\end{figure}
Note that $s_n$ obeys the boundary conditions $s_n (\pm a) =0$. Inserting (\ref{6.1.11}) into the two-dimensional variational expression (\ref{6.1}) yields a sequence of one-dimensional expressions. 
The following proposition relates the associated min-max problems with each other.
\begin{prop}
 Let $s_n$ as in (\ref{6.1.12}) and let $v_n$, $\psi_n : (0,b) \ra \real$ such that $s_n v_n \in \HC_{\al + \ga} (Q)$, $s_n \psi_n \in \HC_{\al - \ga} (Q)$ for any $n \in \nat$. Let, furthermore, $\FC [v, \psi; G]$, $C_c$, and $\Ct_c$ as given in (\ref{6.1})--(\ref{6.2a}). Then 
 \be \label{6.1.15}
 C_c (Q) = \inf_{n \in \nat}\ \inf_{s_n v_n \in \HC_{\al+ \ga} (Q)}\ \sup_{s_n \psi_n \in \HC_{\al - \ga} (Q)}\ \FC [s_n v_n, s_n \psi_n ; Q]\, ,
\ee
\be \label{6.1.15a}
 \Ct_c (Q) = \inf_{n \in \nat}\ \inf_{s_n \psi_n \in \HC_{\al - \ga} (Q)} \ \sup_{s_n v_n \in \HC_{\al+ \ga} (Q)}\ \FC [s_n v_n, s_n \psi_n ; Q]\, .
\ee
\end{prop}
\textsc{Proof:} Let us start with the observation that $\{s_n : n \in \nat\}$ is a complete orthonormal set in $L^2 ((-a, a))$. We may thus assume for $v$ and $\psi$ $L^2$--expansions in the variable $x$ with coefficients $v_n (y)$ and $\psi_n (y)$, respectively:
\be \label{6.1.16}
v(x,y) = \sum_{n = 1}^\infty  s_n (x) \, v_n (y)\  ,\qquad \psi (x,y) = \sum_{n =1}^\infty  s_n (x) \, \psi_n (y)\, .
\ee
Inserting (\ref{6.1.16}) into $\FC [v, \psi]$ and performing the $x$--integration one obtains
\be \label{6.1.17}
\FC [v,\psi] = \frac{\disp \sum_{n=1}^\infty \NC_n [v_n, \psi_n]}{\disp \Big(\sum_{n =1}^\infty \big(\DC_n^+ [v_n]\big)^2 \Big)^{1/2} \Big(\sum_{n =1}^\infty \big(\DC_n^- [\psi_n]\big)^2 \Big)^{1/2}}
\ee
with
$$
\ba{c}
\disp \NC_n [v_n, \psi_n] := \int_0^b \big(v'_n \psi'_n + k_n^2\, v_n \psi_n \big) y^{- \al} \rd y\, , \\[2ex]
\disp \DC_n^+ [v_n] := \bigg(\int_0^b \big({v'_n}\!^2 + d^2 k_n^2\, v_n^2 \big) y^{- (\al + \ga)} \rd y\bigg)^{\frac{1}{2}} , \\[2ex]
\disp \DC_n^- [\psi_n] := \bigg(\int_0^b \big({\psi'_n}\!^2 + e^2 k_n^2\, \psi_n^2 \big) y^{- (\al - \ga)} \rd y\bigg)^{\frac{1}{2}} .
\ea
$$
Fixing $v(x,y)$, i.e.\ now fixing $\{v_n (y) : n \in \nat \}$, and using Cauchy-Schwarz's inequality in the form 
$$
\frac{\disp \sum_{n = 1}^\infty a_n }{\disp \Big( \sum_{n=1}^\infty b_n^2 \Big)^{1/2} } \leq \bigg( \sum_{n= 1}^\infty \Big(\frac{a_n}{b_n} \Big)^2 \bigg)^{1/2}\, ,\qquad b_n > 0\, ,\quad n \in \nat
$$
we may estimate: 
\be \label{6.1.18}
\sup_{\{\psi_n :\, n \in \nat \}}\ \frac{\disp \sum_{n=1}^\infty \NC_n [v_n, \psi_n]}{\disp \Big(\sum_{n =1}^\infty \big(\DC_n^- [\psi_n]\big)^2 \Big)^{1/2}} \leq \Bigg( \sum_{n = 1}^\infty \bigg(
\sup_{\psi_n }\, \frac{\disp \NC_n [v_n, \psi_n]}{\disp \DC_n^- [\psi_n]}\bigg)^2 \Bigg)^{1/2}\, .
\ee
Inequality (\ref{6.1.18}) is in fact an equality as can be seen by a more careful consideration that makes use of an optimal choice of the relative amplitudes of the $\psi_n$. Let us demonstrate this by just two terms. Replacing $\psi_2$ by $\la\, \psi_2$ one obtains
$$
\ba{c}
\disp \sup_{\{\psi_1 ,  \psi_2\}}\, \sup_{\la \in\real}\ \frac{\disp \NC_1 [v_1 , \psi_1] + \NC_2 [v_2, \la\, \psi_2]}{\disp \Big( \big(\DC_1^- [\psi_1]\big)^2 + \big( \DC_2^- [\la\, \psi_2] \big)^2 \Big)^{1/2}} = \sup_{\{\psi_1 , \psi_2\}}\, \sup_{\la \in\real}\ \frac{\disp \NC_1 [v_1 , \psi_1] + \la\, \NC_2 [v_2, \psi_2]}{\disp \Big( \big(\DC_1^- [\psi_1]\big)^2 + \la^2 \big( \DC_2^- [\psi_2] \big)^2 \Big)^{1/2}} \\[5ex]
\disp = \Bigg( \sup_{\psi_1} \bigg(\frac{\disp \NC_1 [v_1, \psi_1]}{\disp \DC_1^- [\psi_1]}\bigg)^2 + \sup_{\psi_2} \bigg(\frac{\disp \NC_2 [v_2, \psi_2]}{\disp \DC_2^- [\psi_2]}\bigg)^2 \Bigg)^{1/2} ,
\ea
$$
as 
$$
\sup_{\la \in \real}\, \frac{a_1 + \la \, a_2}{\big(b_1^2 + \la\, b_2^2\big)^{1/2}} = \bigg( \Big(\frac{a_1}{b_1} \Big)^2 + \Big( \frac{a_2}{b_2}\Big)^2\bigg)^{1/2} ,
$$
where w.l.o.g.\ we may assume $a_{1} > 0$, $a_{2} > 0$, $b_{1} > 0$, $b_{2} > 0$. Iterating this argument yields (\ref{6.1.18}) as equality. 

Applying (\ref{6.1.18}) as equality to (\ref{6.1.17}) yields 
\be \label{6.1.19}
\ba{l}
\disp \inf_{v \in \HC_{\al+ \ga} (Q)}\ \sup_{\psi \in \HC_{\al - \ga} (Q)}\ \FC [v, \psi] \\[-2ex]
\disp \qquad \qquad = \inf_{\{v_n :\, n \in \nat \}}\, \frac{1}{\disp \Big( \sum_{n = 1}^\infty \big( \DC_n^+ [v_n]\big)^2 \Big)^{1/2}}\,
\sup_{\{\psi_n :\, n \in \nat \}}\ \frac{\disp \sum_{n=1}^\infty \NC_n [v_n, \psi_n]}{\disp \Big(\sum_{n =1}^\infty \big(\DC_n^- [\psi_n]\big)^2 \Big)^{1/2}} \\[4ex]
\disp \qquad \qquad \qquad = \inf_{\{v_n :\, n \in \nat \}} \, \Bigg( \sum_{n = 1}^\infty \frac{\big(\DC_n^+ [v_n] \big)^2}{\disp \sum_{m = 1}^\infty \big(\DC_m^+ [v_m]\big)^2}\ \bigg( \sup_{\psi_n }\, \frac{\disp \NC_n [v_n, \psi_n]}{\disp \DC_n^+ [v_n]\, \DC_n^- [\psi_n]}\bigg)^2 \Bigg)^{1/2} \\[6ex]
\disp \qquad \qquad \qquad \qquad \geq \inf_{n \in \nat}\ \inf_{v_n }\ \sup_{\psi_n }\, \frac{\disp \NC_n [v_n, \psi_n]}{\disp \DC_n^+ [v_n]\, \DC_n^- [\psi_n]}\, .
\ea
\ee
On the other hand, let $\big((v_n , \psi_n)\big)_{n \in \nat}$ be an optimizing sequence for the right-hand side of (\ref{6.1.19}). The sequence $\big((s_n v_n , s_n \psi_n)\big)_{n \in \nat}$, when inserted into $\FC [\cdot\, , \cdot ]$, clearly satisfies
$$
\FC[s_n v_n, s_n \psi_n] = \sup_\psi \FC[s_n v_n, \psi]\, ,
$$
and hence demonstrates that 
$$
\inf_v\, \sup_\psi \, \FC [v, \psi] \leq  \inf_{n \in \nat}\ \inf_{v_n }\ \sup_{\psi_n }\, \frac{\disp \NC_n [v_n, \psi_n]}{\disp \DC_n^+ [v_n]\, \DC_n^- [\psi_n]}\, .
$$
Inequality (\ref{6.1.19}) is thus an equality, too. By (\ref{6.2}) this is the assertion (\ref{6.1.15}).
With obvious modifications these arguments apply to (\ref{6.1.15a}) as well.
\qed
%
%
%
%
\subsection{A one-dimensional analytical test case}
An especially simple one-dimensional problem is the $x$-independent case. It allows a direct analytic solution of the min-max problem (see appendix D). Here we consider the corresponding Euler-Lagrange equations, which allow likewise an analytic solution. 
This solution is useful to estimate the min-max value and its dependence on the parameters also in the general case and, even more important, it gives the decisive clue how to obtain the sufficient lower bounds in subsection 6.4.

Let us introduce new variables, viz.,
\be \label{6.2.1}
w:= y^{- \beta_+} v := y^{-(\al + \ga)/2}\, v \, ,\qquad \chi := y^{- \beta_-} \psi := y^{-(\al -\ga)/2}\, \psi\, ,
\ee
which by (\ref{3.11}) have the advantage to be both in $H_0^1$. In the $x$-independent case system (\ref{6.1.5}) then takes the form
\be \label{6.2.2}
\left.\ba{c}
\disp \chi_0'' - \frac{\beta_+ - \beta_-}{y}\, \chi_0' - \frac{(\beta_+ + 1)\beta_-}{y^2} \,\chi_0 = \mu_0 \Big[w_0'' - \frac{(\beta_+ + 1)\beta_+}{y^2} \, w_0 \Big] \, ,\\[3ex]
\disp w_0'' + \frac{\beta_+ - \beta_-}{y}\, w_0' - \frac{(\beta_- + 1)\beta_+}{y^2} \, w_0 = \mu_0 \Big[\chi_0'' - \frac{(\beta_- + 1)\beta_-}{y^2} \, \chi_0 \Big] \, .
\ea \right\}
\ee
A power ansatz of type
$$
\left(\! \ba{c} w_0 \\ \chi_0 \ea \!\right) = \left(\! \ba{c} a \\ b \ea \!\right) y^\de ,\qquad a, b, \de \in \real
$$
is here successful and yields without difficulty 4 linear independent solutions:
$$
 \left(\! \ba{c} 1 \\ 0 \ea \!\right) y^{-\beta_+}\, , \quad  \left(\! \ba{c} 0 \\ 1 \ea \!\right) y^{-\beta_-}\, , \quad
 \left(\! \ba{c} \mu_0 \, \beta_0 \\ 2\, \beta_+ + 1 \ea \!\right) y^{\beta_+ + 1}\, , \quad
 \left(\! \ba{c} 2\, \beta_- + 1 \\ \mu_0 \,\beta_0 \ea \!\right) y^{\beta_- + 1}\, ,
$$
where $\beta_0 := \beta_+ + \beta_- + 1$. Implementing the boundary conditions 
$$
w_0(b_0) = \chi_0(b_0) = 0\, ,\qquad w_0(b_1) = \chi_0(b_1) = 0\, ,\qquad 0 \leq b_0 < b_1
$$
yields a homogeneous $4\ti 4$ system in the amplitudes, whose determinant must vanish:
\be \label{6.2.3}
\left| 
\ba{cccc}
b_0^{-\beta_+} & 0 & \mu_0 \, \beta_0 \, b_0^{\beta_+ + 1} & (2\, \beta_- + 1)\, b_0^{\beta_- + 1} \\[2ex]
0 & b_0^{-\beta_-} & (2 \, \beta_+ + 1)\, b_0^{\beta_+ +1} & \mu_0 \, \beta_0\, b_0^{\beta_- + 1} \\[2ex]
b_1^{-\beta_+} & 0 & \mu_0 \, \beta_0 \, b_1^{\beta_+ + 1} & (2\, \beta_- + 1)\, b_1^{\beta_- + 1} \\[2ex]
0 & b_1^{-\beta_-} & (2 \, \beta_+ + 1)\, b_1^{\beta_+ +1} & \mu_0 \, \beta_0\, b_1^{\beta_- + 1} 
\ea \right| = 0\, .
\ee
Equation (\ref{6.2.3}) determines a unique nonnegative eigenvalue $\mu_0$, viz.\ 
\be \label{6.2.4}
\ba{c}
\disp \mu_0^2 = \bigg(1 - \Big(\frac{\beta_+ - \beta_-}{\beta_+ + \beta_- + 1}\Big)^2 \bigg) \frac{\big(1 - (b_0/b_1)^{\beta_+ + \beta_- + 1} \big)^2}{\big(1 - (b_0/b_1)^{2\, \beta_+ + 1} \big) \big(1 - (b_0/b_1)^{2 \, \beta_- + 1} \big) } \\[4ex]
\disp \quad  = \bigg(1 - \Big(\frac{\ga}{1 + \al}\Big)^2 \bigg) \frac{\big(1 - (b_0/b_1)^{1 + \al } \big)}{\big(1 - (b_0/b_1)^{1 + \al + \ga} \big) }\, \frac{\big(1 - (b_0/b_1)^{1 + \al } \big)}{\big(1 - (b_0/b_1)^{1 + \al - \ga} \big) } \\[4ex]
\disp \quad = \frac{(1 + \gt)(1 - B)}{1 - B^{1 + \gt}}\, \frac{(1 - \gt)(1 - B)}{1 - B^{1 - \gt}}\, ,
\ea
\ee
where in the last line we have set 
$$
\gt := \frac{\ga}{1 + \al}\, ,\qquad B:= \Big(\frac{b_0}{b_1}\Big)^{1+ \al} .
$$
As can be read off (\ref{6.2.4}), $\mu_0$ depends only on the ratio $b_0/b_1$, which is in accordance with the scaling property (\ref{6.6}); $B\mapsto \mu_0$ is a monotonically increasing function on $[0,1]$ taking its minimum value $(1 - \gt^2)^{1/2}$ at $B= 0$, which corresponds to an interval that touches the origin or stretches up to infinity. $\gt \mapsto \mu_0$ is a monotonically decreasing function on $[0, \infty)$ with $\mu_0 = 1$ at $\gt =0$ in accordance with (\ref{6.10}). Figure 6.2 shows curves $\mu_0$ versus $\ga$ for some ratios $b_0/b_1$; note that $\ga > 1 + \al$ makes sense in the case that $b_0/b_1 > 0$.
\begin{figure}
\begin{center}
\input{Fig62}
\caption{$\mu_0$ versus $\ga$ for $\al = 0.39$ and $b_0/b_1 = 0$, $10^{-6}$, $10^{-4}$, $10^{-2}$ (from bottom to top).}
\end{center}
\end{figure}
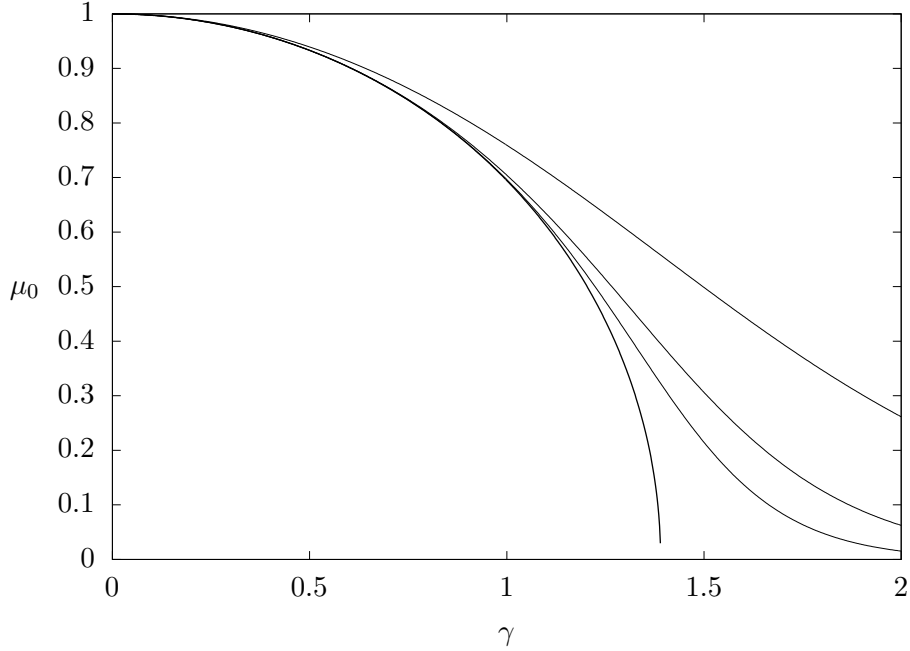

Computing $w_0$ and $\chi_0$ explicitly we find up to positive factors and a global sign the unique optimal pair
$$
\ba{c}
\disp w_0 = (1-B) \Big(\frac{y}{b_1}\Big)^{1 + (\al + \ga)/2} - \big(1 - B^{1 + \gt}\big) \Big(\frac{y}{b_1}\Big)^{1 + (\al - \ga)/2} + B \big(1 - B^\gt\big) \Big(\frac{y}{b_1}\Big)^{-(\al + \ga)/2}\  , \\[3ex]
\disp \chi_0 = \big(1-B^{1- \gt}\big) \Big(\frac{y}{b_1}\Big)^{1 + (\al + \ga)/2} - (1 - B) \Big(\frac{y}{b_1}\Big)^{1 + (\al - \ga)/2} - B \big(1 - B^{-\gt}\big) \Big(\frac{y}{b_1}\Big)^{-(\al - \ga)/2}
\ea
$$
with the remarkable property that $w_0$ and $\chi_0$ conicide in the limit $B = 0$:
\be \label{6.2.5}
w_0 = \Big(\frac{y}{b_1}\Big)^{1 + (\al + \ga)/2} - \Big(\frac{y}{b_1}\Big)^{1 + (\al - \ga)/2} = \chi_0 \quad \mbox{ on } \ (0, b_1)\, .
\ee
Note, by the way, that the uniqueness of the solution according to (\ref{6.1.10}) implies $C_c = \mu_0 = \Ct_c$ for the $x$-independent versions of problems (\ref{6.2}) and (\ref{6.2a}).
\subsection{Lower bounds on $C_c$ and $\Ct_c$}
Lower bounds on $C_c$ and $\Ct_c$ are obtained by replacing the maximizing function by a test function leaving us with a pure minimization problem. When using the ``symmetrical'' variables $w$ and $\chi$, the test case of the preceeding subsection suggests the kind of test function that yields sharp bounds -- at least for large domains.
Note, moreover, that differently to min-max values a pure minimum has the advantage to depend monotonically on the domain.

Let us start with rewriting the functional (\ref{6.1}) in the variables (\ref{6.2.1}). The $v$-related denominator then takes the form
$$
\ba{l}
\disp \DC^+ [v] := \bigg(\int_G \big( (d\, \pa_x v)^2 + (\pa_y v)^2\big) y^{-(\al + \ga)} \rd x \rd y \bigg)^{1/2} \\[3ex]
\disp \qquad \qquad = \bigg(\int_G \Big( (d\, \pa_x w)^2 + \Big(\pa_y w + \beta_+\, \frac{w}{y}\Big)^2\,\Big) \rd x \rd y \bigg)^{1/2} \\[3ex]
\disp \qquad \qquad \qquad = \bigg(\int_G \Big( |\na_d\, w|^2 + \de_+\Big(\frac{w}{y}\Big)^2\, \Big) \rd x \rd y \bigg)^{1/2} =: \DCt^+ [w]\, ,
\ea
$$
where we made use of the notation $\beta_\pm := (\al \pm \ga)/2$ and $\de_\pm := \beta_\pm (\beta_\pm + 1)$. In the last line we applied integration by parts. Similarly one obtains
$$
\ba{l}
\disp \DC^- [\psi] := \bigg(\int_G \big( (e\, \pa_x \psi)^2 + (\pa_y \psi)^2\big) y^{-(\al - \ga)} \rd x \rd y \bigg)^{1/2} \\[3ex]
\disp \qquad \qquad \ = \bigg(\int_G \Big( |\na_e\, \chi|^2 + \de_-\Big(\frac{\chi}{y}\Big)^2\, \Big) \rd x \rd y \bigg)^{1/2} =: \DCt^- [\chi]\, ,
\ea
$$
$$
\ba{l}
\disp \NC [v, \psi] := \int_G \big(\pa_x v\, \pa_x \psi + \pa_y v\, \pa_y \psi \big) \rd x \rd y \\[3ex]
\disp \qquad \qquad \ \ = \int_G \Big( \pa_x w\, \pa_x \chi + \Big(\pa_y w + \beta_+\, \frac{w}{y} \Big) \Big( \pa_y \chi + \beta_-\, \frac{\chi}{y}\Big) \Big) \rd x\rd y = : \NCt [w, \chi]\, ,
\ea
$$
and hence
$$
\FC[v, \psi;G] = \frac{\NCt [w, \chi]}{\DCt^+ [w]\, \DCt^- [\chi]} =: \FCt [w, \chi; G]\, .
$$
In view of (\ref{6.2.5}) we now replace $\chi$ by $w$ and obtain thus a lower bound on $C_c$:
\be \label{6.4.1}
\ba{l}
\disp C_c (A^+_R) = \inf_{w \in H^1_0 (A^+_R)}\  \sup_{\chi \in H^1_0 (A^+_R)}\  \FCt [w, \chi; A^+_R] \\[3ex]
\disp \qquad \qquad  \quad \geq \inf_{w \in H^1_0 (A^+_R)}\,  \FCt [w, w; A^+_R]\ \geq \inf_{w \in H^1_0 (A^+_\infty)}\, \FCt [w, w; A^+_\infty] \\[3ex]
\disp \qquad \qquad  \qquad \quad   = \inf_{w \in H^1_0 (H)}\,  \FCt [w, w; H]\ = \inf_{w \in H^1_0 (Q)}\, \FCt [w, w; Q]\, ,
\ea
\ee
where in the last line we made use of the scaling property (\ref{6.6}). As in subsection 6.2 the two-dimensional variational problem can be reduced to one dimension by expanding $w$ into the complete orthonormal system $\{s_n (x) : n \in \nat\}$ with $s_n (x)$ given by (\ref{6.1.12}): 
\be \label{6.4.2}
w(x,y) = \sum_{n=1}^\infty s_n (x)\, w_n (y)\, .
\ee
Inserting (\ref{6.4.2}) into $\DCt^+ [w]$ yields
$$
\ba{l}
\disp \DCt^+ [w] = \bigg(\sum_{n=1}^\infty\, \int_0^b \Big( d^2 k_n^2\, w_n^2 + {w'_n}\!^2 + \de_+\Big(\frac{w_n}{y}\Big)^2\,\Big)  \rd y\bigg)^{1/2} \\[3ex]
\disp \qquad \qquad = \bigg(\sum_{n=1}^\infty\, k_n \int_0^{k_n b} \Big( d^2  \wt_n^2 + {\wt'_n}\!^2 + \de_+\Big(\frac{\wt_n}{z}\Big)^2\,\Big)  \rd z\bigg)^{1/2} \\[3ex]
\disp \qquad \qquad \qquad =: \bigg(\sum_{n=1}^\infty k_n \Big(\DCt^+_n [\wt_n]\Big)^2 \bigg)^{1/2} ,
\ea
$$
where we have introduced the variables $z:= k_n y$ and $\wt_n (z) := w_n (y)$. An analogous calculation holds for $\DCt^- [w]$, whereas $\NCt [w, w]$ can be further simplified by integration by parts:
$$
\ba{l}
\disp \NCt [w,w] = \sum_{n=1}^\infty \,\int_0^b \Big( k_n^2\, w_n^2 + \Big(w'_n + \beta_+\, \frac{w_n}{y} \Big) \Big(w'_n + \beta_-\, \frac{w_n}{y} \Big) \Big) \rd y \\[3ex]
\disp \qquad  \qquad \quad = \sum_{n=1}^\infty \, \int_0^b \Big( k_n^2\, w_n^2 + {w'_n}\!^2 + \de_0 \Big(\frac{w_n}{y} \Big)^2 \, \Big) \rd y \\[3ex]
\disp \qquad  \qquad \qquad \quad = \sum_{n=1}^\infty \, k_n \int_0^{k_n b} \Big( \wt_n^2 + {\wt'_n}\!^2 + \de_0 \Big(\frac{\wt_n}{z} \Big)^2\, \Big) \rd z \\[3ex]
\disp \qquad \qquad \qquad \qquad \quad =: \sum_{n=1}^\infty \, k_n\ \NCt_n [\wt_n]
\ea
$$
with the abbreviation $\de_0 := \beta_+ \beta_- + (\beta_+ + \beta_-)/2$. Using the elementary inequality
$$
\disp \frac{\disp\sum_{n=1}^\infty a_n}{\disp \Big( \sum_{n=1}^\infty b_n^2\Big)^{1/2} \Big( \sum_{n=1}^\infty c_n^2\Big)^{1/2}} \, \geq \, \inf \bigg\{\frac{a_n}{b_n\, c_n} : n \in \nat \bigg\}\, ,\quad a_n \geq 0\, ,\, b_n >0\, ,\, c_n >0\, ,
$$
(\ref{6.4.1}) can thus be further estimated as follows
$$
\ba{l}
\disp C_c (A_R^+)\, \geq\, \inf_{w \in H_0^1 (Q)} \ \frac{\NCt[w,w]}{\DCt^+ [w]\, \DCt^- [w]} \geq\ \inf_{n \in \nat}\  \inf_{\wt_n \in H_0^1((0, k_n b))} \ \frac{\NCt_n[\wt_n]}{\DCt^+_n [\wt_n]\, \DCt^-_n [\wt_n]} \\[3ex]
\disp \qquad \qquad \quad \geq \, \inf_{\wt \in H_0^1((0, \infty))} \ \frac{\NCt_\infty[\wt]}{\DCt^+_\infty [\wt]\, \DCt^-_\infty [\wt]}\, ,
\ea
$$
where the index ``$\infty$'' indicates the interval $(0,\infty)$ of integration.
Introducing the ratios
$$
s:= \frac{\disp \int_0^\infty \wt^2 \,\rd z}{\disp \int_0^\infty {\wt'}\,\!^2 \,\rd z}\ , \qquad 
t:= \frac{1}{4}\ \frac{\disp \int_0^\infty \Big(\frac{\wt}{z}\Big)^2 \rd z}{\disp \int_0^\infty {\wt'}\,\!^2 \,\rd z}
$$
we obtain the final lower bound on $C_c(A^+_R)$:
\be \label{6.4.3}
C_c (A_R^+) \, \geq \inf_{\mbox{\scriptsize{$\ba{c}
0\leq s < \infty \\ 0\leq t \leq 1 
\ea$ }}} \ f(s,t)
\ee
with 
$$
f(s,t) := \frac{1 + s + 4 \, \de_0\, t}{(1 + d^2 s+ 4\, \de_+\, t)^{1/2} (1 + e^2 s + 4\, \de_-\, t)^{1/2} }
$$
and 
$$
\de_\pm = (\al + 1 \pm \ga)^2 -1\ , \quad \de_0 = (\al + 1)^2 - \ga^2 -1\, .
$$
Note that the range of $t$ is restricted by inequality (\ref{3.5}). The minimization in (\ref{6.4.3}) is elementary with minima obtained either at $s=0$ or $t=1$ (depending on the parameters $\al$, $\ga$, $d$, and $e$). Denoting this lower bound by $L\!B = L\!B (\al, \ga, d,e)$, it obviously holds for $\Ct_c (A_R^+)$ as well and, moreover, for any $R> 1$. 

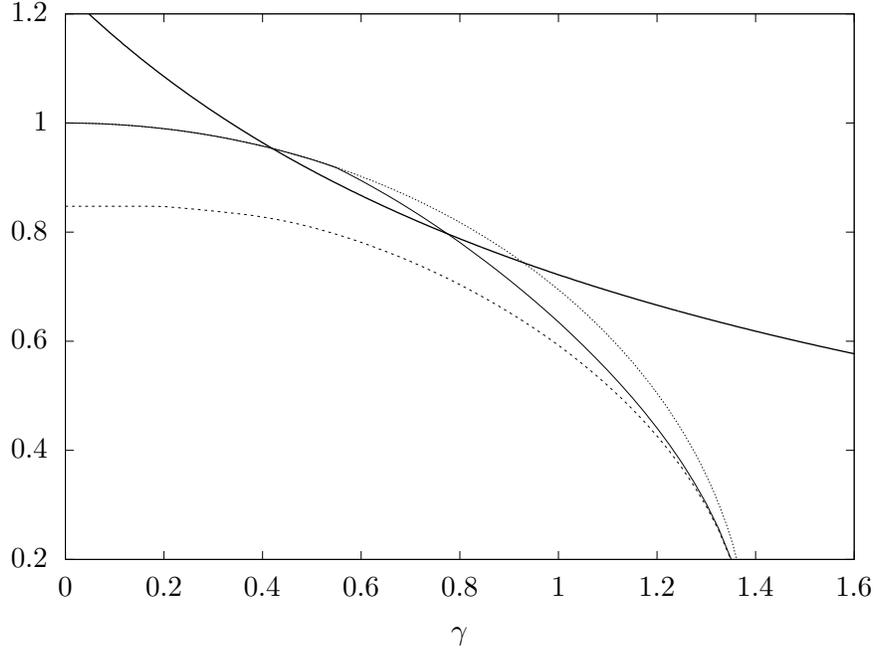
\begin{figure}
\begin{center}
\input{Fig63}
\caption{The lower bound $L\!B (\al, \ga, d,e)$ with $e=1.2$ and $d=1$ (dashed line) or $d= 0.1$ (solid thin line) together with the comparison function $C\!F(\al, \ga)$ (solid thick line) and the exact $x$-independent curve $\mu_0 = \big( 1- (\ga/(1 + \al))^2 \big)^{1/2}$ (dotted line). $\al = 0.39$ has been set throughout.}
\end{center}
\end{figure}
According to (\ref{5.18}), (\ref{5.19}) this bound has to be compared with the ``comparison function''
$$
C\!F(\al, \ga) := \frac{2 \sqrt{2}\, \| a_\rho - \al\|_\infty}{1 + \al + \ga}\, .
$$
Figure 6.3 displays the curves $\ga \mapsto L\!B (0.39, \ga, d, 1.2)$ for $d=1$ and 0.1 together with $\ga \mapsto C\!F(0.39, \ga)$, where $\al$ and $e$ have been fixed according to (\ref{5.3}), (\ref{5.4}), and (\ref{5.17}).
For $d=1$ the bounding curve lies strictly below the comparison function, for $d=0.1$, however, a ``$\ga$-window''opens up, where the lower bound $L\!B$ on $C_c(A_R^+)$ and $\Ct_c (A_R^+)$ clearly exceeds its comparison values. Further lowering of $d$ does not improve $L\!B$ in a significant way. Also shown in Fig.\ 6.3 is the exact $x$-independent curve $\ga \mapsto \mu_0$ in the limit $B=0$, which coincides with $L\!B (0.39, \ga, 0.1, 1.2)$ up to $\ga \approx 0.56$; this is in accordance with the fact that up to this value the minimum in (\ref{6.4.3}) is taken at $s=0$ and
$$
\inf_{0\leq t \leq 1}\, f(0,t)\, = \Big(1 - \Big(\frac{\ga}{1+ \al}\Big)^2\, \Big)^{1/2}\, .
$$
Above  $\ga \approx 0.56$ the minimum is located at $t=1$.

To be definite let us fix $\ga = 0.6$ and $d = 0.1$. We then have the proposition:
\begin{prop}
 Let
 \be \label{6.4.4}
 (\al, \ga, d,e) = (0.39, 0.6, 0.1, 1.2)\, ,
 \ee
 then inequalities (\ref{5.18}) and (\ref{5.19}) are satisfied for any $R> 1$ and for the quantity $\De$ from (\ref{5.20}) holds the lower bound
 \be \label{6.4.5}
 \De \geq 0.027 > 0\, .
 \ee
\end{prop}
Better lower bounds on $C_c$, $\Ct_c$ can be obtained according to proposition 6.2 by a numerical evaluation of the Euler-Lagrange equations (\ref{6.1.5}) (see appendix E) corroborating on one side the lower bound $L\!B$ and demonstrating on the other side that $L\!B$ generally is not sharp. 
\sect{Solution of $P_\Om$ in annuli}
This section contains the proof of theorem 2.1 based on Schauder's fixed point principle and the general solution of the linearized problem in section 5.

The following proposition provides a suitable version of the Leray-Schauder theorem (see Gilbarg and Trudinger 1998, p.\ 280, theorem 11.3).
\begin{prop}[fixed point principle]
 Let $T : \BC \ra \BC$ be a continuous, compact (or completely continuous) mapping of the Banach space $\BC$ into itself and let the set 
 \be \label{7.1}
 \big\{ u \in \BC : u = \la\, T u \mbox{ for some } \la \in [0,1]\big\}
 \ee
 be bounded. Then, $T$ has a fixed point $u_0 \in \BC$, i.e.
 \be \label{7.2}
 T u_0 = u_0\, .
 \ee
\end{prop}
In order to find a (nonlinear) mapping $T$ suitable for our purpose let us fix the parameters $\al$, $\ga$, $d$, and $e$ according to proposition 6.4 such that conditions (\ref{5.11}), (\ref{6.4a}), (\ref{5.18}), and (\ref{5.19}) are satisfied for any $R > 1$. Let $\BC := \LC_{p, \de} (A_R^+)$ with 
$$
\LC_{p, \de} (G) := clos \big( C_0^\infty (G) \, ,\, \|\cdot\,  y^{-\de} \|_{p, 0} \big)\, ,\qquad G \subset H
$$
and $p$ and $\de$ yet to be fixed. With $w \in \LC_{p, \de} (A_R^+)$ and 
\be \label{7.2a}
\aaa^\al := \big(a_\ze [w \rho^{-\al} - \Om]\, ,\, a_\rho [ w \rho^{-\al} - \Om] - \al\big)\, ,
\ee
where $a_\ze [\cdot]$ and $a_\rho [\cdot]$ are given by (\ref{2.24}), we solve eq.\ (\ref{5.2}) and obtain by proposition 5.2 a solution $v \in \HC^{(d)}_{\al + \ga} (A_R^+)$ with bound 
\be \label{7.3}
\|\na_d\, v\|_{\al + \ga} \leq \frac{\sqrt{2}}{\De\, d}\, \|\Om \|_{2 + \ga - \al} \leq \bigg(\frac{4 \, \pi}{(\ga- \al)(1 + \al - \ga)}\bigg)^{\!\frac{1}{2}} \frac{K}{\De\, d}\, .
\ee
The last inequality in (\ref{7.3}) follows by condition (\ref{2.30}) as in lemma A.1. Supposing that
\be \label{7.4}
\de < \frac{\al + \ga}{2} + \frac{1}{p}\, ,
\ee
we find by (\ref{3.20}) that 
\be \label{7.5}
\HC^{(d)}_{\al + \ga} (A_R^+)  \subset \LC_{p, \de } (A_R^+)\quad \mbox{for any }\ p \geq 2\, .
\ee
The mapping 
\be \label{7.6}
T : \LC_{p, \de} (A_R^+) \ra \LC_{p, \de} (A_R^+)\ , \quad w \mapsto v
\ee
is thus well defined and, moreover, continuous and compact.

To see the continuity of $T$ let $v$ and $\vt$ be two solutions of eq.\ (\ref{5.2}) with (\ref{7.2a}) and
$$
\widetilde{\aaa}^\al := \big(a_\ze [\wt \rho^{-\al} - \Om]\, ,\, a_\rho [ \wt \rho^{-\al} - \Om] - \al\big)\, ,
$$
respectively, $w$ and $\wt \in \LC_{p, \de}(A_R^+)$, $\psi \in C_0^\infty (A_R^+)$, and $\Om = \widetilde{\Om} \in \HC_{2 + \ga -\al} (A_R^+)$.
The difference of (\ref{5.2}) for $v$ and $\vt$ reads
$$
\ba{l}
\disp \int_{A_R^+} \na (v - \vt) \cdot \na \vii\, \rd \mu_\al  + \int_{A_R^+} \Big(\frac{v}{\rho}\, \aaa^\al [ w \rho^{-\al} - \Om] - \frac{\vt}{\rho}\, \aaa^\al [ \wt \rho^{-\al} - \Om]\Big)\cdot \na \vii\, \rd \mu_\al \\[2ex]
\disp \qquad \qquad =
\int_{A_R^+} \frac{\Om}{\rho} \Big( \aaa [w \rho^{-\al} - \Om] - \aaa[\wt \rho^{-\al} - \Om] \Big) \cdot \na \vii\, \rd \ze \rd \rho\, ,
\ea
$$
or, after some rearrangement,
$$
\ba{l}
\disp \int_{A_R^+} \Big(\na (v - \vt) + \frac{v - \vt}{\rho}\, \aaa^\al [ w \rho^{-\al} - \Om] \Big)\cdot \na \vii\, \rd \mu_\al \\[2ex]
\disp \qquad \qquad =
- \int_{A_R^+} \Big( \frac{\vt}{\rho} - \frac{\Om}{\rho}\, \rho^\al \Big) \big( \aaa [w \rho^{-\al} - \Om] - \aaa[\wt \rho^{-\al} - \Om] \big) \cdot \na \vii\, \rd \mu_\al\, .
\ea
$$
Using (\ref{5.12}) on the left-hand side and estimating the right-hand side by Cauchy-Schwarz's and H\"older's inequality one obtains
\be \label{7.7}
\ba{l}
\disp B[v - \vt, \psi ] \leq \int_{A_R^+} L\, |w - \wt|\, |\vt/\rho^{1 + \al} - \Om /\rho| \, |\na \psi|\, \rd \mu_\al \\[3ex]
\disp  \qquad  \leq 
L\, \big\| |w - \wt|\, |\vt/\rho^{1 + \al} - \Om /\rho| \big\|_{\al +\ga}\, \|\na \psi\|_{\al - \ga} \\[2ex]
\disp \qquad \qquad \leq \|(w - \wt) \rho^{-(\al + \ga)/2} \|_{p, 0}\, \big(\|\vt/\rho^{1+ \al} \|_{q, 0} + \|\Om/\rho\|_{q,0} \big) \, \|\na_e \psi\|_{\al - \ga}\, ,
\ea
\ee
where $\frac{1}{p} + \frac{1}{q} = \frac{1}{2}$ and $L\leq 1$ is some Lipschitz constant for $\aaa[\cdot]$. The critical step is to find a bound on the parantheses in the last line. A bound on $\Om$ is easy by (\ref{2.30}):
\be \label{7.8}
\Big\|\frac{\Om}{\rho}\Big\|_{q, 0} \leq | A_R^+|^{\frac{1}{q}}\,\Big\| \frac{\Om}{\rho}\Big\|_\infty \leq |A_R^+|^{\frac{p-2}{2 p}} K\, .
\ee
A bound on $\vt$, however, requires the refined inequality (\ref{3.18}) for ``small $q$'':
\be \label{7.9}
\Big\|\frac{\vt}{\rho^{1 + \al}}\Big\|_{q, 0} \leq C\, \|\na_d \, \vt\|_{\al + \ga}
\ee
provided that 
$$
1 + \al \leq \frac{\al + \ga}{2} + \frac{6- q}{2 q} \ \ \Longleftrightarrow \ \ \frac{1}{q} \geq \frac{1}{2} - \frac{\ga -\al}{6}\ \ \Longleftrightarrow \ \ p \geq \frac{6}{\ga - \al}\, ,
$$
which means for the parameter set (\ref{6.4.4}) that 
\be \label{7.10}
p \geq 29\, .
\ee
By (\ref{7.3}), inequality (\ref{7.7}) may thus be rewritten as 
\be \label{7.11}
\frac{B[v- \vt, \psi]}{\|\na_e \psi\|_{\al - \ga} } \leq \breve{C}\, K\,\|(w - \wt) \rho^{-(\al +\ga)/2} \|_{p,0}\, ,
\ee
where $\breve{C}$ denotes a constant depending on $p$, $\al$, $\ga$, $d$, $e$, and $R$. By approximation (\ref{7.11}) holds for any $\psi \in \HC^{(e)}_{\al - \ga} (A_R^+)$. In view of (\ref{7.11}) and (\ref{7.4}) we choose $\de := (\al + \ga)/2$ in (\ref{7.5}) and find by combination of (\ref{3.20}), (\ref{5.17a}) with (\ref{5.20}), and (\ref{7.11}):
$$
\ba{l}
\|(v -\vt) \rho^{-(\al + \ga)/2}\|_{p,0} \leq \Ch\, \|\na_d\, (v -\vt)\|_{\al + \ga}  \\[2ex]
\disp \qquad \qquad \qquad \leq \frac{\Ch}{\De}\, \sup_{\psi \neq 0}\, \frac{B[v- \vt, \psi]}{\|\na_e \psi\|_{\al - \ga} } \leq \frac{\Ch\, \breve{C} K}{\De}\,\|(w - \wt) \rho^{-(\al +\ga)/2} \|_{p,0}\, ,
\ea
$$
which implies continuity  of $T$ in $\LC_{p, (\al+ \ga)/2}\, (A_R^+)$.

To see the compactness of $T$ it is sufficient to prove compactness of the embedding (\ref{7.5}). This, however, is equivalent to the well-known compact emdedding $H^1_0 (G) \Subset L^p (G)$ for bounded $G \subset \real^2$ and $1 \leq p < \infty$. In fact, by (\ref{3.11}) we can conclude that a bounded sequence $(v_n) \subset \HC_\beta (A_R^+)$ implies boundedness of the sequence $(v_n\, \rho^{-\beta/2} ) \subset H_0^1 (A_R^+)$ and convergence of a sequence $(v_m\, \rho^{- \beta/2} ) \subset L^p (A_R^+)$ is equivalent to convergence of $(v_m) \subset \LC_{p, \beta/2}\, (A_R^+)$. 

Finally, boundedness of the set (\ref{7.1}) is immediate by (\ref{3.20}) and (\ref{7.3}):
$$
\ba{l}
\|w \rho^{-(\al + \ga)/2}\|_{p,0} \leq \Ch\, \|\na_d\, w\|_{\al + \ga} = \Ch\, \|\na_d (\la\, T w)\|_{\al + \ga} \\[2ex]
\disp \qquad \qquad \qquad \quad \leq \Ch\, \|\na_d\, v\|_{\al + \ga} \leq \bigg(\frac{4 \, \pi}{(\ga- \al)(1 + \al - \ga)}\bigg)^{\!\frac{1}{2}}\, \frac{\Ch\, K}{\De\, d}
\ea
$$
for any $w \in \LC_{p, (\al + \ga)/2}\, (A_R^+)$ satisfying $w = \la\, T w$ for some $\la \in [0,1]$.

The fixed point $v_0 = w_0$ according to proposition 7.1 is a solution of eq.\ (\ref{5.2}) with 
$$
\aaa^\al := \big(a_\ze [v_0 \rho^{-\al} - \Om]\, ,\, a_\rho [ v_0 \rho^{-\al} - \Om] - \al\big)\, ,
$$
hence $u_0 := v_0 \rho^{-\al}$ is a solution of the (nonlinear) equation (\ref{2.29}) in $A_R^+$, and by antisymmetric continuation in $A_R$.

The bound (\ref{2.31}) follows by (\ref{5.21a}), (\ref{6.4.4}), (\ref{6.4.5}), and (\ref{7.3}). Note that this bound depends (besides the in (\ref{6.4.4}) fixed parameters) on $K$ but not on $R$. This together with proposition 4.1 concludes the proof of theorem 2.1.
\sect{Solution of $P_\Om$ in the exterior plane}

The construction of solutions in the last section heavily depended on estimates that involved the diameter $R$ of the underlying domain, the bound (\ref{2.31}), however, does not. This allows us to construct the solution of problem $P_\Om (A_\infty)$ by suitably manipulating a sequence of solutions on $A_R$ with $R \ra \infty$, and to prove this way theorem 2.2.

More precisely let $\Om_n:= \Om\big|_{A_n}$ for $n \in \nat\setminus\{1\}$ and $\Om \in L^\infty (A_\infty)$ satisfying the bound (\ref{2.30}) for some $K > 0$. By theorem 2.1 we then have a sequence of solutions $(u_n) \subset H^1_{0 , as} (A_n)$ of the problems $P_{\Om_n} (A_n)$, $n \in \nat\setminus\{1\}$, with bounds
\be \label{8.1}
\|\na u_n \|_\beta  \leq \sqrt{C}\, , \qquad n \in \nat\setminus\{1\}\, .
\ee
Restricting $u_n$ onto $A_{2}$, the sequence $(u_n|_{A_{2}} )_{n \in \nat\setminus\{1\}}$ is bounded by (\ref{8.1}). There is thus a subsequence $(n^{(1)})  \subset \nat$ such that $\big(u_{n_i^{(1)}}|_{A_{2}} \big)_{i \in \nat}$ coverges weakly to some $u^{(1)} \in \HC_\beta  (A_{2})$ with bound $\|\na u^{(1)} \|_\beta \leq \sqrt{C}$. Applying this argument to the sequence $\big(u_{n_i^{(1)}}|_{A_{3}} \big)_{i \in \nat}$ one obtains a further subsequence $(n^{(2)} ) \subset (n^{(1)})$ such that $\big(u_{n_i^{(2)}}|_{A_{3}} \big)_{i \in \nat}$ converges weakly to some $u^{(2)} \in \HC_\beta (A_{3})$ with bound $\|\na u^{(2)} \|_\beta \leq \sqrt{C}$. Repeating this process one obtains a sequence $\big( (n^{(k)})\big)_{k \in \nat}$ of sequences with each element being a subsequence of its predecessor. Let $n_k^{(k)}$ denote the $k$th element of the subsequence $(n^{(k)})$ and $(n_k^{(k)})_{k \in \nat} =: n^*$ the diagonal sequence. By construction we have $(n_i^* )_{i \geq k} \subset (n^{(k)})$ and hence 
\be \label{8.2}
u_{n_i^*}\big|_{A_{1 + k}}  
\xrightharpoonup[\text{$i \ra \infty$}]{\text{$\| \na \cdot \|_\beta$}} 
u^{(k)} \in \HC_\beta (A_{1 + k}) \quad \mbox{ for any } \ k \in \nat
\ee
with
\be \label{8.3}
\|\na u^{(k)} \|_\beta \leq \sqrt{C}\, ,\qquad k \in \nat\, .
\ee
As $u^{(k)} |_{A_{1 + l}} = u^{(l)}$ for any $l \leq k$, the definition 
$$
u : A_\infty \ra \real \, ,\qquad u|_{A_{1 + k}} := u^{(k)} \, ,\qquad k \in \nat
$$
is unique and yields by (\ref{8.3}) a function  $\in \HC_\beta (A_\infty) \subset H^1_{loc, as} (A_\infty)$ with bound (\ref{2.33}).

In order to see that $u$ satisfies eq.\ (\ref{2.29}) in $A_\infty$ recall that by antisymmetry it is enough to satisfy (\ref{2.29}) in $A_\infty^+$. So, let $\vii \in C_0^\infty  (A_\infty^+)$ and $l \in \nat\setminus\{1\}$ such that in fact $\vii \in C_0^\infty  (A_l^+)$. For $u_{n^*_i}$ then holds
\be \label{8.4}
\int_{A^+_l} \na u_{n^*_i} \cdot \na \vii\, \rd \ze \rd \rho + \int_{A^+_l} \frac{u_{n^*_i}}{\rho}\, \aaa[u_{n^*_i} - \Om] \cdot \na \vii\, \rd \ze \rd \rho =
\int_{A^+_l} \frac{\Om}{\rho}\, \aaa[u_{n^*_i} - \Om] \cdot \na \vii\, \rd \ze \rd \rho\, .
\ee
Letting $i \ra \infty$  the passage to $u$ in the first term on the left-hand side of (\ref{8.4}) is justified by weak convergence. For the other two terms recall that by compactness weak convergence in $H^1_0 (A_l^+)$ implies strong convergence in $L^2 (A_l^+)$. Therefore, by the estimate
$$
\ba{l}
\disp \bigg| \int_{A_l^+} \Big( \frac{u_{n_i^*}}{\rho}\, \aaa [u_{n_i^*} - \Om ] - \frac{u}{\rho}\, \aaa[ u - \Om] \Big) \cdot \na \vii\, \rd \zeta \rd \rho \bigg| \\[3ex]
\disp \qquad \leq \int_{A_l^+} \Big(|u_{n_i^*}|\, \big|\aaa [u_{n_i^*} - \Om ] - \aaa[ u - \Om] \big| + | \aaa[ u - \Om]|\, |u_{n_i^*} - u| \Big) \Big| \frac{\na \vii}{\rho} \Big|\,  \rd \zeta \rd \rho \\[3ex]
\disp \qquad \qquad \leq M\, L\, \| u_{n_i^*}\|_0\, \|u_{n_i^*} - u\|_0 + M\sqrt{2}\, |A_l^+|^{1/2} \|u_{n_i^*} - u\|_0\, ,
\ea
$$
where $L$ denotes a Lipschitz constant for $\aaa[\cdot ]$ and $M$ a bound on $|\na \vii /\rho|$ in $supp\, \vii \Subset A_l^+$, the passage to $u$ in the second term on the left-hand side is justified, too, and a similar argument applies to the right-hand side.
This together with proposition 4.2 concludes the proof of theorem 2.2.
\sect{Solution of the signed direction problem without zeroes on the symmetry axis}
In this section we prove theorems 2.3 and 2.4 (with the stated restriction), which contain the solution of the axisymmetric signed direction problem (\ref{1.1}) in large spherical shells and in exterior space, respectively, by applying the results gathered so far and formulated essentially in theorems 1.1 and 1.2 and the appendices A and B. This application is straight forward; only to establish continuity of the solution up to the boundary requires some additional considerations. 

To apply theorem 2.1 let $R > R_0$ and let $\DD$ and $\DDh$ be continuous symmetric direction fields with rotation numbers $\ro$ and $\roh$ along the boundaries $S_1$ and $S_R$ of the annulus $A_R$, respectively. Let $\ro \geq \roh \geq 2$ and $\ro - \roh$ be an even number (including zero). Let, furthermore, $\{z_1, \ldots, z_{\ro - \roh}\}$ be a symmetric set of points in $A_R$, which do not lie on the symmetry axis.\footnote{In the case $\ro - \roh = 0$ the set is empty.}
By (\ref{A.6}) the associated zero-positions angle $\Psi$ then satisfies condition (\ref{2.30}). By (\ref{2.18}) and (\ref{2.19}), direction fields $\DD$ and $\DDh$, which satisfy the axis condition (\ref{2.34}) give rise to boundary functions $\phi$ and $\phih$, respectively, that satisfy (\ref{B.1}). By (\ref{B.3}) the harmonic interpolation $\Phi$ then satisfies condition (\ref{2.30}) as well. By theorem 2.1 we thus obtain a weak solution $u \in H_{0, as}^1 (A_R)$ of eq.\ (\ref{2.22}), which implies a weak solution $q = u + \Phi \in H_{as}^1 (A_R)$ of eq.\ (\ref{2.17}). Given $q$, the function $p$ is determined by (\ref{2.20}) up to a constant $p_0$. In fact, writing (\ref{2.20}) in Cartesian coordinates in the form
\be \label{9.1}
\ba{c}
\disp \pa_\ze\, p = \pa_\rho\, q + \frac{\sin (q - \Psi)}{q - \Psi}\, \frac{q -\Psi}{\rho}\, ,\\[2ex]
\disp \pa_\rho\, p = - \pa_\ze\, q + \frac{\cos (q - \Psi) - 1}{q - \Psi}\, \frac{q -\Psi}{\rho}
\ea
\ee
we find by (\ref{3.5}) the right-hand side in $L^2 (A_R)$. The integrability condition for (\ref{9.1}) is just (\ref{2.17}), which takes the weak form
$$
\int_{A_R} \bigg\{\Big[\pa_\ze\, q - \frac{\cos (q - \Psi) - 1}{q - \Psi}\, \frac{q -\Psi}{\rho}\Big] \pa_\ze \vii + 
\Big[\pa_\rho\, q + \frac{\sin (q - \Psi)}{q - \Psi}\, \frac{q -\Psi}{\rho}\Big] \pa_\rho \vii \bigg\} \rd \ze \rd \rho = 0
$$
for any $\vii \in C_0^\infty  (A_R)$. We obtain thus a symmetric function $p \in H^1 (A_R)$ that solves (\ref{2.20}) (see, e.g., Galdi 94, lemma 1.1, p.\ 101, and corollary 4.1, p.\ 53). Inserting $h$ as given by (\ref{2.14}) and $g= (p + i q)/2$ into (\ref{2.13}) then yields a symmetric weak solution $f \in H^1 (A_R)$ of eq.\ (\ref{2.8}) with boundary condition (\ref{2.12}), which by (\ref{2.7a}) is equivalent to an axisymmetric magnetic field $(B_\rho, B_\ze)$ satisfying (\ref{2.1}) a.e.\ in $A_R$. To obtain higher (interior) regularity let us switch from cylindrical to Cartesian coordinates in $\real^3$. The field
$$
\BB_C := (B_1, B_2, B_3) := (B_\rho \cos \theta, B_\rho \sin \theta, B_\ze ) \in \big(H^1(S\! S_R)\big)^3
$$
then satisfies the system 
$$
\na \ti \BB_C = 0\ , \qquad \na \cdot \BB_C = 0 \qquad \mbox{a.e. in } S\! S_R\, ,
$$
where $S\! S_R$ denotes the spherical shell $B_R \setminus \overline{B_1} \subset \real^3$ corresponding to $A_R \subset \real^2$. By the vector analysis identity
$$
\int_D \na \ti \BB_C \cdot \na \ti \ppsi\, \rd v + \int_D \na \cdot \BB_C \ \na \cdot \ppsi\, \rd v = \sum_{i = 1}^3 \int_D  \na B_i \cdot \na \psi_i \, \rd v\, , 
$$
where $\ppsi = (\psi_1, \psi_2, \psi_3) \in (C_0^\infty (D))^3$, one finds each component to be a weakly harmonic and hence harmonic function in $S\! S_R$. $B_i$ and hence $B_\ze$ and $B_\rho$ are thus smooth functions in their respective domains.
\begin{figure}
\begin{center}
\includegraphics[width=0.5\textwidth]{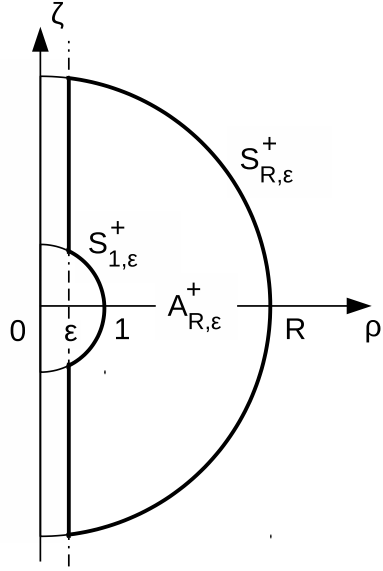}
\caption{$A_{R, \e}^+$ and the curved boundary components $S_{1, \e}^+$ and $S_{R, \e}^+$.}
\end{center}
\end{figure}

Improved boundary regularity, in particular $\BB \in C(\overline{A_R})$, needs some more subtle arguments, which in part are specific for two dimensions. As to the variable $q$ we rely on a boundary-regularity theorem asserting H\"older continuity of a weak solution up to the boundary provided that the boundary values themselves satisfy a H\"older condition (Ladyzhenskaya and Ural'tseva 1968, theorem 14.1, p.\ 201f.). For this purpose let us write eq.\ (\ref{2.17}) in the form 
$$
\Delta q = \pa_\ze \Big( \frac{1}{\rho} \big(\cos (q -\Psi) -1\big) \Big) - \pa_\rho \Big(\frac{1}{\rho}  \sin (q - \Psi) \Big) =: \pa_\ze F_\ze + \pa_\rho F_\rho
$$
with bounded $\FF$ on $A_{R, \e}^+$:
$$
|F_\ze | < \frac{2}{\e}\ ,\qquad |F_\rho| < \frac{1}{\e}\, .
$$
Here, $A_{R, \e}^+ := A_r^+ \cap \{\rho > \e\}$ denotes that part of $A_R^+$, which has a distance $\e > 0$ to the symmetry axis, and $S_{1, \e}^+$ and $S_{R, \e}^+$ denote the corresponding curved boundary components (see Fig.\ 9.1).

The direction fields $\DD$ and $\DDh$ are now assumed to be of H\"older class $C^{0,\beta} (S_1)$ and $C^{0,\beta} (S_R)$, respectively, for some $0< \beta \leq 1$, which implies the boundary functions $\phi$ and $\phih$ to be of the same class. 
By $q - \Phi  \in H_0^1 (A_R)$ we thus have
$$
trace\, q |_{S_{1,\e}^+} = \phi|_{S_{1,\e}^+}\ , \qquad trace\, q |_{S_{R,\e}^+} = \phih|_{S_{R,\e}^+}\, .
$$
To this situation the boundary-regularity theorem may be applied asserting $q \in C^{0, \al} (\overline{A_{R, \e}^+})$ for some $0 < \al \leq 1$, which depends on $\e$, $R$, $\|q \|_0$, and $\beta$, $\|\phi\|_{C^{0, \beta}}$, $\|\phih\|_{C^{0, \beta}}$. This implies, in particular, 
continuous assumption of the boundary values on $S^+_{1,\e} \cup S^+_{R,\e}$ for any $\e >0$, and hence on $S_1^+ \cup S_R^+$. 

This argument does not apply to $p$ since a-priori no boundary values are prescribed for this quantity. However, according to the governing equations for $p$, eqs.\ (\ref{9.1}), the variables $p$ and $q$ are complex conjugate variables up to a (in $A^+_{R, \e}$) bounded ``perturbation''. Without this perturbation Privaloff's theorem guarantees uniform H\"older continuity of both variables in a disc if only one variable takes continuously boundary values that are in some H\"older class at the boundary (see Bers et al.\ 1964, p.\ 279 or Kaiser 2010, lemma 2.4). In order to take advantage of Privaloff's theorem the following representation theorem (slightly adapted to our needs) is useful (see Bers and Nirenberg 1954 or Bers et al.\ 1964, p.\ 259).
\begin{thm}[Bers and Nirenberg]
 Let $G \subset B_R \subset \cpl$ and let $w \in H^1 (G)$ be a solution of the equation
 \be \label{9.2}
 \pa_\zb\, w = \ga \qquad \qquad \mbox{ a.e. in }\, G
 \ee
 with $\ga \in L^\infty (G)$. Then $w$ has the representation 
 \be \label{9.3}
 w(z, \zb) = f(z) + s_0 (z, \zb) \qquad \quad\mbox{ in } \, G\, ,
 \ee
 where $f$ is complex analytic in $G$ and where $s_0$ is uniformly H\"older continuous in $B_R$ and real on $S_R$.
\end{thm}
With 
$$
\Re z := \ze\, , \qquad \Im z := \rho\, , \qquad \Re w := p \, ,\qquad \Im w := q\, ,\qquad \Re \ga := F_\rho \, ,\qquad \Im \ga := F_\ze\, ,
$$
eqs.\ (\ref{9.1}) are clearly of type (\ref{9.2}) in $G:= A^+_{R, \e/2}$; (\ref{9.3}) then yields on $S_{R,\e}^+$:
$$
q|_{S^+_{R,\e}} = \Im w|_{S^+_{R,\e}} = \Im f|_{S^+_{R,\e}} .
$$
$\Im f$ is thus of class $C^{0, \al}$ on $S^+_{R, \e}$ and even analytic on (say) $\pa (\At^+_{R,\e}) \setminus S^+_{R, \e}$, where
$$
\At^+_{R,\e} := \big\{(\ze, \rho) \in H: (1 + \e)^2 < \ze^2+ \rho^2 < R^2 ,\, \rho > \e\big\}\, .
$$
Applying Privaloff's theorem on $f$ in the Lipschitz-domain\footnote{By a conformal transformation with Lipschitz-continuous extension on the boundary, $\At^+_{R, \e}$ is equivalent to a disc.} $\At^+_{R, \e}$ yields $f \in C^{0, \al} (\overline{\At^+_{R, \e}})$, which implies, in particular,
$$
p = \Re w = \Re f + \Re s_0
$$
to be H\"older continuous up to the boundary component $S^+_{R, \e}$. To see continuity up to the other boundary component $S^+_{1, \e}$ let us apply the inversion $z \mapsto \frac{1}{z}$ to (\ref{9.2}), which yields an equation of the same type in the inverted half-annulus
$$
\At^-_{1/R\, ,1} := \big\{(\ze, \rho) \in H: (1/R^2 < \ze^2+ \rho^2 < 1 ,\, \rho < 0 \big\}
$$
with outer boundary $S_1^-$. The above argument can thus be repeated establishing, finally, continuity of $p$ up to the boundary $S_1^+ \cup S_R^+$.

Concerning uniqueness let us remind that the number of zeroes is fixed by (\ref{2.6}) for any solution $f$ of problem $P_{\DD, \DDh}^c (A_R)$. 
Two solutions $f_1$ and $f_2$ with the same boundary direction fields $\DD^c$ and $\DDh^c$ and the same set of zeroes are thus represented by $u_1 = q_1 - \Phi$ and $u_2 = q_2 - \Phi$ satisfying (\ref{2.23}) with the same boundary function $\Phi$ and zero-positions angle $\Psi$, and hence the same $\Om = \Psi - \Phi$. According to proposition 4.1 there is $u_1 = u_2$ and the only nonuniqueness left is due to $p_0$, which in (\ref{9.1}) remains free. Therefore, $f_1$ and $f_2$ differ by the positive factor $e^{p_0/2}$. This proves theorem 2.3 in the case of no zeroes lying on the symmetry axis. 

With some obvious modifications all these arguments work in the unbounded case as well: $\roh + 1$ has to be replaced by the exact decay order $\dt$, $\DDh$ and $\phih$ can be discarded, and $H^1_{as} (A_R)$ is replaced by $H^1_{loc, as} (A_\infty)$, existence and uniqueness of $q$ then follow from theorem 2.2 and proposition 4.2, $p$ is again determined by eqs.\ (\ref{9.1}), and all the regularity arguments, local in character, apply to the unbounded case as well. This proves theorem 2.4 (again with the mentioned restriction on the zeroes).
\sect{Zeroes on the symmetry axis}
Zeroes on the symmetry axis are so far excluded. Overcoming this restriction by continuous deformation of known solutions is the subject of this section. 

Let us start with a lemma on a general invariance property of the rotation number that is tailored to our needs (see, e.g., Heinz 1959, theorem 3 or Kaiser 2010, lemma 3.2). 
\begin{lem}[deformation invariance]
 Let $\BB_\la$ be a continuous vector field on $\overline{A_R} \ti I$, where $\la$ varies in the closed interval $I$. Let $\BB_\la \neq 0$ on $S_\Rt \ti I$, $1 \leq \Rt \leq R$, and let $\rot_\la$ be the rotation number of $\BB_\la$ along $S_\Rt$. Then, $\rot_\la = const =: \rot$ for any $\la \in I$. 
\end{lem}
In the case that $\BB_\la$ is a harmonic vector field in $A_R$ or $A_\infty$, by (\ref{2.6}) and (\ref{2.7}) we can draw from lemma 10.1 some
consequences, which in a somewhat lax formulation can be summarized as follows:

(i) In the case that $B|_{(S_1 \cup S_R) \ti I} \neq 0$, zeroes of $\BB_\la$ can shift with $\la$, can split or coalesce, but they can never spring up or vanish in $A_R$.

(ii) If a zero with index $-n$, $n \in \nat$, crosses $S_\Rt$, $1 < \Rt < R$, at (say) $\la = \la_0$, according to (\ref{2.6}), $\rot_\la$ jumps at $\la_0$ by $+n$ in the case of ``emigration'' and by $-n$ in the case of ``immigration'' into $A_\Rt$.

(iii) In $A_\infty$ zeroes can vanish at infinity or originate from infinity at (say) $\la = \la_0$; according to (\ref{2.7}) the exact decay order $\dt$ increases or decreases at $\la_0$ according to the index of the zero.

(iv) If $\rot$ is fixed (as is the case in the direction problem for $\ro$ along $S_1$ and $\roh$ along $S_R$), zeroes cannot cross $S_\Rt$ (cannot originate or vanish at $S_1$ or $S_R$). However, zeroes can migrate from $A_R$ onto $S_1$ (or $S_R$), split there and migrate as ``border zeroes'' along $S_1$ (or $S_R$), and they can coalesce again and reenter into $A_R$ (see appendix E for a (two-dimensional) example of this behaviour).

In the following we denote by 
$$
\BB = \BB^\dt  = \BB^\dt [z_1, \ldots, z_{\ro - \dt +1}]
$$
the unique solutions of problem $P_\DD (A_\infty)$
according to theorem 2.4 with rotation number $\ro$, exact decay order $\dt$, and zeroes $z_1, \ldots, z_{\ro - \dt +1}$, which do not lie on the symmetry axis. In particular, let 
$$
\BB^\dt [\ldots, z, \zb]\, ,\qquad \BBt^\dt [\ldots, \zt, \ztb]\, ,\qquad  \BB^{\dt +2} [\ldots]
$$
be three solutions of $P_\DD (A_\infty)$ with $\ro - \dt -1$ common zeroes indicated by ``$\ldots$'' and, additionally, the simple zeroes $z$, $\zb$ in the first case and $\zt$, $\ztb$ in the second case. Let $\BB_{\lala} = \BB_{(\la, \lat, \lah)}$ be the linear combination
$$
\BB_\lala := \la\, \BB^\dt + \lat \, \BBt^\dt + \lah\, \BB^{\dt +2}\, ;
$$let, furthermore, $z_S = \ze_S$ and $\zt_S = \zet_S$ be two positions on the symmetry axis in $A_\infty$, and let $c_{\dt - 2}$ and $\ct_{\dt -2}$ be the leading coefficients in the representation (\ref{2.5}) of $\BB^\dt$ and $\BBt^\dt$, respectively.

For $z \neq \zt$ the fields $\BB^\dt$, $\BBt^\dt$, and $\BB^{\dt + 2}$ are clearly linearly independent and the linear systems (in the variables $\la$, $\lat$, $\lah$)
\be \label{10.1}
B_{\ze, \lala} (\ze_S, 0) = 0\, , \qquad \la\, c_{\dt -2} + \lat\, \ct_{\dt -2} = 0
\ee
and
\be \label{10.2}
B_{\ze, \lala} (\ze_S, 0) = 0\, , \qquad B_{\ze, \lala} (\zet_S, 0) = 0
\ee
have each up to (for the present arbitrary) factors unique solutions $\lala_1$ and $\lala_2$ such that $\BB_{\lala_1}$ has the $\ro - \dt -1$ common zeroes outside the symmetry axis and the simple zero $z_S$ on the symmetry axis, and $\BB_{\lala_2}$ has, besides the common zeroes, the simple zeroes $z_S$ and $\zt_S$ on the symmetry axis.
Note that by (\ref{10.1})$_b$ $\BB_{\lala_1}$ has the decay order $\dt + 1$, and by (\ref{10.1})$_a$ and (\ref{2.7}) this decay order is exact.

There is one question left open: Do $\BB_{\lala_{1,2}}$ satisfy the boundary condition $\BB_{\lala_{1,2}}|_{S_1} = a_{1,2}\, \DD$ for some {\em positive} amplitudes $a_{1,2}$? A-priori the linear combinations $\BB_{\lala_{1,2}}$ can have additional zeroes on $S_1$. To refute this presumption  let $\mu \mapsto \lala (\mu)$, $\mu \in [0,1]$ be a continuous path in $\la$-space connecting $\lala = (1,0,0)$ with $\lala_1$; $\mu \mapsto \BB_{\lala (\mu)}$ is thus a continuous deformation of $\BB_{\lala (0)} = \BB^\dt  [\ldots, z, \zb]$ with final state $\BB_{\lala (1)} = \BB^{\dt +1}_{\lala_1} [\ldots, z_S]$. During the deformation the common zeroes remain fixed, whereas $z$, $\zb$ will move. They can possibly enter $S_1$ and leave again; finally, however, one zero migrates to $z_S$, whereas the other vanishes at infinity. Additional zeroes cannot emerge, neither from $S_1$ (see (iv)) nor from infinity (see (iii)\footnote{Observe that by construction the decay order of $\BB_{\lala(\mu)}$ cannot fall below $\dt$.}). Therefore, $\BB_{\lala_1}|_{S_1} \neq 0$, i.e., $\BB_{\lala_1}$ is a signed solution of the direction problem. A similar argument applies to $\BB_{\lala_2}$, which is thus a signed solution of the direction problem, too.

Further zeroes can be shifted to the symmetry axis just by repeating the procedure with solutions, which possess already common zeroes on the symmetry axis. In the case that $\ro - \dt +1 = \nu$ is even, this way all zeroes can be shifted to arbitrary positions at the symmetry axis. In the case that $\ro  - \dt + 1 = \nu$ is odd, one simply solves the direction problem with $\dt - 1$ instead of $\dt$; we are then again in the even case and solution of the system (\ref{10.1}) provides us with a field that has an odd number of zeroes on the symmetry axis, an even number $\geq 0$ of zeroes outside, and the original exact decay order $\dt$. 

In the case of a bounded domain $A_R$, solutions of type (\ref{10.2}) still exist, i.e., an even number of zeroes can likewise be shifted to the symmetry axis; (\ref{10.1}), however, does no longer work: zeroes cannot be expelled from $A_R$ with its fixed direction fields at the boundaries $S_1$ and $S_R$ (see (iv)).

In conclusion these remarks remove the limitation in theorems 2.3 and 2.4 we had to impose in section 9.
\sect{The unsigned direction problem}
This concluding section discusses the unsigned direction problem, its relation to the signed problem, and proves, in particular, theorem 2.5. 

Let us start with recalling some definitions given in section 2: $S_\DD^\dt$ means the set of all solutions of the signed direction problem $P_\DD^s (A_\infty)$ with direction field $\DD$ and exact decay order $\dt$, and $L_\DD^\de$ means the linear space of all solutions of the unsigned problem $P_\DD^u (A_\infty)$ with decay order $\de$. Taking arbitrary linear combinations does clearly not respect the positivity of the amplitude function nor the exact decay order, but it respects the boundary condition up to a (locally varying) sign and the decay order. When taking only positive linear combinations, which defines a cone in $L_\DD^\de$, positivity of the amplitude and decay order are preserved, exactness of the decay order, however, can get lost. The set $C_\DD^\de$ of all solutions of $P_\DD^s (A_\infty)$ with decay order $\de$ is of this type.

The proof of theorem 2.5 is based on the following relations between these sets:
\be \label{11.1}
C_\DD^\de = \bigcup_{\de \leq \dt \leq \ro +1} S_\DD^\dt\, ,
\ee
\be \label{11.2}
L_\DD^\de = \left\{
\ba{cl}
\langle C_\DD^\de \rangle & \ \mbox{ if } \, \ro \geq \de -1 \\[1ex]
\{0\} & \ \mbox{ if }\, 2 \leq \ro < \de -1
\ea \right. ,
\ee
\be \label{11.3}
S_\DD^\dt \subset \langle S_\DD^\de \rangle \qquad \mbox{ for } \, 3 \leq \de \leq \dt \leq \ro +1\, ,
\ee
\vspace{-2ex}
\be \label{11.4}
dim \, \langle S_\DD^\de \rangle = \ro - \de + 2 \qquad \mbox{ if }\, 3 \leq \de \leq \ro +1 \, .
\ee
Here $\ro$ denotes again the rotation number of $\DD$ and $\langle S \rangle$ the real linear span of the set $S$.\\[-2ex]

(\ref{11.1}) is obvious by definition and (\ref{2.7}).\\[-2ex]

As to (\ref{11.2}): Let $\ro \geq \de -1$. To prove the nontrivial inclusion $L_\DD^\de \subset \langle C_\DD^\de \rangle$, let $\BB\in L_\DD^\de$ with $\BB|_{S_1} = a\, \DD$ and $a \in C(S_1)$. By theorem 2.4  we have $C_\DD^\de \neq \emptyset$. Let $\BBt \in C_\DD^\de$  with $\BBt|_{S_1} = {\tilde a } \, \DD$ and continuous ${\tilde a} > 0$. Choosing $\la > 0$ such that $a + \la {\tilde a} > 0$ on $S_1$, we have 
$$
\BB_\la := \BB + \la\, \BBt \in C_\DD^\de\, ,
$$
which implies for $\BB$ the representation
$$
\BB = \BB_\la - \la\, \BBt\, ,
$$
i.e., $\BB \in \langle C_\DD^\de \rangle$.

In the case $2 \leq \ro < \de - 1$ let $\BB \in L_\DD^\de$ and let $\BBt^{\ro +1}$ be the {\em unique} solution $\in S_\DD^{\ro + 1}$. Choose as above $\la > 0$ such that $a + \la {\tilde a} > 0$ on $S_1$. Then, by $\ro + 1 < \de $ the linear combination 
$$
\BB + \la \, \BBt^{\ro +1} =: \BB_\la^{\ro + 1}
$$
has exact decay order $\ro + 1$, i.e., we find $\BB_\la^{\ro +1} \in S_\DD^{\ro + 1}$ as well. By uniqueness we have $\BB_\la^{\ro +1} = \mu\, \BBt^{\ro + 1}$ for some $\mu > 0$, which implies $\BB = (\mu - \la) \BBt^{\ro + 1}$. Thus, either $\BB$ has the exact decay order $\ro + 1$, which contradicts $\BB \in L_\DD^\de$ with $\de > \ro + 1$, or $\BB \equiv 0$, hence $L_\DD^\de  = \{0 \}$.\\[-2ex]

As to (\ref{11.3}): Choose arbitrary $\BB^i \in S_\DD^i$, $i= \de, \ldots , \ro +1$; then, the set $S:= \{\BB^i : \de \leq i \leq \ro +1 \}$ is linearly independent as the $\BB^i$ differ pairwise by their asymptotic terms in the representation (\ref{2.5}). Defining 
$$
\BBt_i^\de := \BB^i + \BB^\de\, ,
$$
the set ${\widetilde S}:= \{\BBt^i : \de \leq i \leq \ro +1 \}$ is likewise linearly independent and, moreover, we have ${\widetilde S} \subset S_\DD^\de$. To prove (\ref{11.3}) we show that any $\BB^\dt  \in S_\DD^\dt$ for any $\de \leq \dt \leq \ro + 1$ has a representation in ${\widetilde S}$:
\be \label{11.5}
\BB^\dt = \sum_{i = \de}^{\ro +1} \la_i\, \BBt_i^\de \, , \qquad \la_i \in \real\, .
\ee
As $\BB^\dt$ has exact decay order $\dt$ and precisely $\ro -\dt +1$ zeroes, say $z_1, \ldots, z_{\ro - \dt +1}$, the $\la_i$ have to satisfy the following linear system of equations:
$$
\ba{cc}
\disp \sum_{i = \de}^{\ro + 1} \la_i\, \BBt_i^\de (z_n) = 0\, ,& \quad n = 1, \ldots ,\ro - \dt +1\, ,\\[3ex]
\disp \sum_{i = \de}^{\ro + 1} \la_i\, \ct_{i, n} = 0\, ,& \quad \de \leq n \leq \dt -1\, ,
\ea
$$
where $\ct_{i, n}$ is the coeficient $c_n$ in the representation (\ref{2.5}) of $\BBt_i^\de$. These $(\ro - \dt + 1) + (\dt - \de) = \ro - \de + 1$ equations for $\ro - \de + 2$ unknowns have always a nontrivial solution, say $\{ \mu_i : \de \leq i \leq \ro +1\}$. Let 
$$
{\widehat \BB}^\dt := \sum_{i = \de}^{\ro +1} \mu_i\, \BBt_i^\de \, .
$$
If $\BBh^\dt |_{S_1} \neq 0$ either $\BBh^\dt$ or $-\BBh^\dt$ belongs to $S_\DD^\dt$ and has, moreover, the same zeroes as $\BB^\dt$. By uniqueness we thus have $\BB^\de = \mu\, \BBh^\dt$ with $\mu \neq 0$, i.e., (\ref{11.5}) holds. If $\BBh^\dt$ has zeroes on $S_1$ consider
$$
\breve{\BB}^\dt := \BBh^\dt + \la\, \BB^\dt
$$
with $\la > 0$ such that $\breve{\BB}^\dt |_{S_1} \neq 0$, i.e., $\breve{\BB}^\dt \in S_\DD^\dt$. Again by uniqueness we find 
$$
\BB^\dt = \mu\, \breve{\BB}^\dt \ \Longleftrightarrow \ (1 - \mu \la) \BB^\dt = \mu\, \BBh^\dt
$$
with $\mu> 0$ and hence $1 - \mu \la \neq 0$, i.e., (\ref{11.5}) holds again.\\[-2ex]

As to (\ref{11.4}): This is an immediate consequence of the preceeding arguments: the set ${\widetilde S}$ is linearly independent by construction and generates $S_\DD^\de$ according to (\ref{11.5}). ${\widetilde S}$ with its $\ro -\de + 2$ elements is thus a basis of $S_\DD^\de$.\\[-2ex]

The proof of theorem 2.5 now follows by (\ref{11.1})--(\ref{11.4}) and the following concluding remark.\\[-2ex]

\noi {\em Remark}: In the ``trivial case'' $\de > \ro +1$, the condition $\ro \geq 2$ as stated in (\ref{11.2})$_2$ is not necessary; the argument given there, however, must be altered. Equation (\ref{2.7}) holds in fact for any $\ro \in \ints$, which implies immediately $S_\DD^\de = \emptyset$. To prove $L_\DD^\de  = \emptyset$ assume $0 \not \equiv \BB \in L_\DD^\de$, which then has necessarily zeroes on $S_1$. Let $\e>0$ such that $\BB \neq  0$ on $\overline{B_{1 + \e}} \setminus \overline{B_1}$ and compare $\ro$ with the rotation number $\ro_\e$ of $\BB$ along $S_{1+ \e}$. For $\e \ra 0$, $\ro_\e$ remains constant; however, it will differ from $\ro$ by $-\frac{1}{2} \ti\!$ number of sign reversals of $\BB$ along $S_1$ (for an illustration see Fig.\ E.2). Thus, $\ro_\e < \ro$ and applying (\ref{2.7}) on $\BB$ in $\real^2\setminus B_{1+ \e}$ yields a contradiction. The assumption $\BB \not \equiv 0$ has thus to be dismissed.
\setcounter{equation}{0}
\setcounter{section}{0}
\renewcommand{\theequation}{\thesection.\arabic{equation}}
\renewcommand{\thesection}{\Alph{section}}

\sect{The zero-positions angle $\Psi$}
Let us write the analytic function (\ref{2.14}) in the form
\be \label{A.1}
h(z) = z^{-\ro} \prod_{n =1}^N (z- z_n) (z - \zb_n) = z^{-\roh} \prod_{n =1}^N \Big(1- \frac{z_n}{z}\Big) \Big(1 - \frac{\zb_n}{z}\Big)\, ,\quad \ro - \roh = 2\, N\, ,\ N \in \nat_0
\ee
with points $\{z_1, \ldots , z_N\} \subset A_\infty^+$, which satisfy, in particular, $\Im z_n = \rho_n > 0$. Decomposing $\hb/h$ into real and imaginary parts one obtains
\be \label{A.2}
\frac{\hb}{h} = v(\ze, \rho) + i\, w(\ze , \rho)\, , \qquad (\ze, \rho) \in A_\infty
\ee
with 
\be \label{A.3}
v= \frac{\hb^2 + h^2}{2 \,|h|^2}\, ,\qquad w = \frac{\hb^2 - h^2}{2 i\, |h|^2}\, , \qquad v^2 + w^2 = 1\, .
\ee
We now define 
\be \label{A.4}
\Psi (\ze, \rho) := \arctan \frac{w(\ze, \rho)}{v (\ze, \rho)}\, ,
\ee
where $\arctan$ means the principal branch of the inverse $\tan$ function with $\arctan 0 = 0$. This definition resolves the ambiguity in the choice of $\Psi$ and fixes the discontinuities of $\Psi$ at the zeroes of $v$. Let us illustrate this by the special case $h_0 = z^{-\ro}$. With $z = r e^{i \vi}$ one obtains
$$
\frac{\hb_0}{h_0} = e^{2 i \ro \vi} = \cos 2 \ro\, \vi + i \sin 2 \ro\, \vi
$$
and hence\footnote{Recall that tilde denotes dependence on polar coordinates: $\Pst (r, \vi) := \Psi 
(r \cos \vi, r \sin \vi)$.}
$$
\Pst_0  = \arctan (\tan 2 \ro\, \vi)\, ;
$$
$\Pst_0$ jumps by $\pi$ along all rays with angles $\vi_\nu = \frac{\pi}{4 \ro} \, (1 + 2 \nu)$, where $\nu$ is an integer and $\vi_\nu \in (-\pi, \pi)$ (see Fig.\ A.1).
\begin{figure}
\begin{center}
\includegraphics[width=0.4\textwidth]{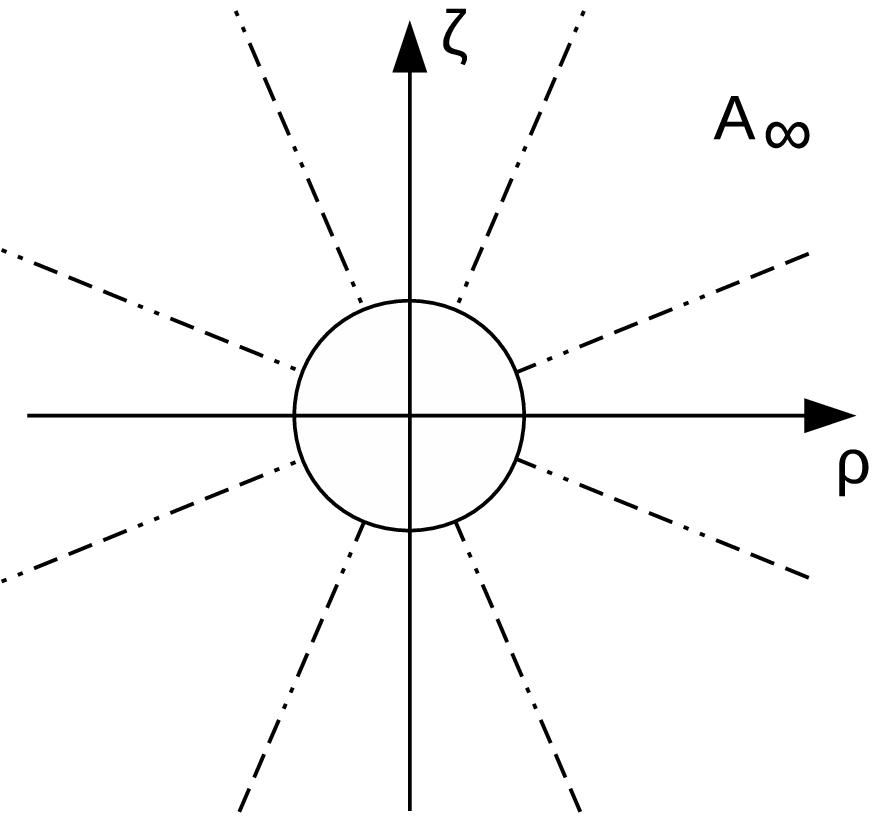}
\caption{Lines of discontinuity of $\Psi$ for $\ro = \roh = 2$.}
\end{center}
\end{figure}
Moreover, by the reflection symmetry $\Pst_0 (\cdot , \pi - \vi) = - \Pst_0(\cdot , \vi)$ and the estimate $\vi \leq \frac{\pi}{2} \sin \vi$ on $[0, \pi/2]$, one finds the estimate
\be \label{A.5}
|\Pst_0 (\cdot, \vi)| \leq \pi \ro\, |\sin \vi| \quad \mbox{in }\; A_R\, .
\ee
According to the following lemma inequality (\ref{A.5}) (and some consequences) hold for general $\Psi$.
\begin{lem}
 Let $\Psi$ as given by (\ref{A.1}) -- (\ref{A.4}). Then the following inequalities hold:
 \be \label{A.6}
 |\Pst ( \cdot , \vi) |\leq K\, |\sin \vi|\qquad \mbox{in } \ A_\infty \, ,
 \ee 
 \be \label{A.7}
 |\Psi ( \cdot , \rho)/\rho |\leq K \qquad \mbox{in } \ A_\infty \, ,
 \ee
 \be \label{A.8}
 \int_{A_\infty} \Psi^2\, \rho^{-(2 + \beta)}\, \rd \ze \rd \rho \leq \frac{2 \pi K^2}{\beta (1 - \beta)}\, , \qquad 0 < \beta < 1 
  \ee
  for some $K> 0$.
\end{lem}
\textsc{Proof:} Inequality (\ref{A.7}) is an immediate consequence of (\ref{A.6}):
$$
|\Psi ( \cdot , \rho)/\rho | = |\Pst (r, \vi)/ r \sin \vi | \leq K/ r \leq K \qquad \mbox{in } \ A_\infty \, ,
$$ 
and (\ref{A.8}) follows by integration:
$$
\ba{l}
\disp \int_{A_\infty} \Psi^2\, \rho^{-(2 + \beta)}\, \rd \ze \rd \rho = 
2 \int_0^\pi \int_1^\infty \frac{\Pst^2 (r, \vi)}{r^{2+ \beta}\, |\sin \vi|^{2+ \beta}}\, r\, \rd r \rd \vi \\[3ex]
\disp \qquad \leq 4 K^2 \int_1^\infty r^{-(1 + \beta)}\, \rd r\, \int_0^{\pi/2} |\sin \vi|^{- \beta} \rd \vi \leq \frac{2 \pi K^2}{\beta (1 - \beta)}\, , 
\ea
$$
where we made use of $\sin \vi \geq 2 \vi/\pi$ on $[0, \pi/2]$.

It remains to prove (\ref{A.6}). For this purpose we decompose $A_\infty^+$ into wedge-type components $A_1 \cup A_2 \cup A_3$ with
$$
\ba{l}
\disp A_1 := \big\{ (r,\vi): 1 < r < \infty\, ,\ \vi_0 \leq \vi \leq \pi - \vi_0 \big\} , \\[2ex]
\disp A_2 := \big\{ (r,\vi): 1 < r < \Rt\, ,\ 0 < \sin \vi < \sin \vi_0 \big\} , \\[2ex]
\disp A_3 := \big\{ (r,\vi): r \geq \Rt\, ,\ 0 < \sin \vi < \sin \vi_0 \big\}
\ea
$$
(see Fig.\ A.2).
\begin{figure}
\begin{center}
\includegraphics[width=0.5\textwidth]{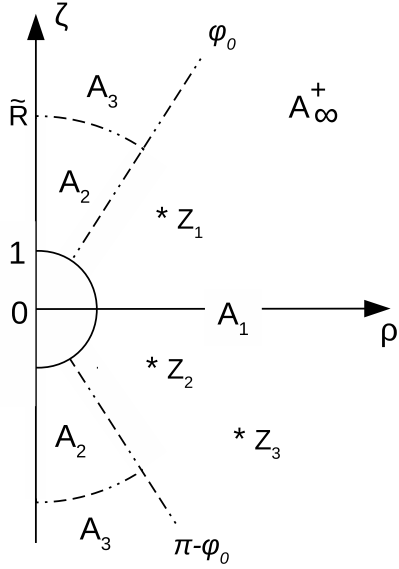}
\caption{The decomposition $A_\infty^+ = A_1 \cup A_2 \cup A_3$.}
\end{center}
\end{figure}
Here $\Rt$ (and $\e$) are chosen such that 
\be \label{A.9}
\Rt > \frac{1 }{\e} \max_{1 \leq n \leq N} |z_n|\, ,\qquad \frac{8 \e}{1 - 3 \e} < \frac{1 }{2 C(N)} < 1 \, ,
\ee
where $C(N)$ is some combinatorial constant that appears in inequality (\ref{A.11}). Note that by (\ref{A.3}) $v(\ze, \rho)$ and $w(\ze, \rho)$ are real analytic functions in a neighbourhood of the symmetry axis $\{\rho = 0\}$ with $v(\cdot, 0) = 1$ and $w(\cdot , 0) =0$. $\vi_0$ can thus be chosen such that $v \geq 1/2$ and $|w| \leq c\, \rho$ in $A_2$ for some $c > 0$. Moreover, it is assumed that 
\be \label{A.9a}
\vi_0 < \frac{\pi}{8\, \roh}\, .
\ee
Under these conditions inequality (\ref{A.6}) then holds trivially in $A_1 \cup A_2$:
$$
\ba{l}
\disp |\Pst (\cdot, \vi) | \leq \frac{\pi}{2 \sin \vi_0} \, \sin \vi \qquad \mbox{in } \ A_1\, , \\[3ex]
\disp |\Pst (\cdot, \vi) | \leq \frac{w}{v} \leq 2c \, \Rt \sin \vi \qquad \mbox{in } \ A_2\, .
\ea
$$
For $A_3$ we need a more careful consideration. Let 
$$
g_n := \frac{\big( \overline{1 - (z_n/z)\big)\big(1 - (\zb_n/ z)}\big) }{\big(1 - (z_n/ z)\big)\big(1 - (\zb_n/ z)\big)}\, ,\qquad n = 1, \ldots , N\, .
$$
By $z = r e^{i \vi}$, $z_n = r_n\,
e^{i \vi_n}$, $g_n$ may be expressed as 
$$
g_n = \frac{1 - 2\, (r_n/r)\, e^{i \vi} \cos \vi_n + (r_n/ r)^2\, e^{2 i \vi}}{1 - 2 \, (r_n/ r)\, e^{- i \vi} \cos \vi_n + (r_n/ r)^2\, e^{- 2 i \vi}}
= 1 + i\,  \frac{ - 4\, (r_n/r) \sin \vi \cos \vi_n + 2 (r_n/ r)^2\, \sin 2 \vi}{1 - 2 \, (r_n/ r)\, e^{- i \vi} \cos \vi_n + (r_n/ r)^2\, e^{- 2 i \vi}}\, ,
$$
which allows by (\ref{A.9}) the estimate

\be \label{A.10}
|g_n - 1| \leq \frac{ 8 \e}{1 - 3 \e} \sin \vi = \eta \sin \vi\, .
\ee
By (\ref{A.10}) and the abbreviation $8 \e/(1 - 3 \e) =: \eta$ it follows for $g_N := \prod_{ n =1}^N g_n$:
\be \label{A.11}
\ba{l}
\disp \big|g_N^2 - 1\big| = \bigg| \prod_{n=1}^N g_n^2 - 1\bigg| = \bigg|\prod_{n=1}^N \big( (g_n - 1)(g_n + 1) + 1\big) -1 \bigg| \\[3ex]
\disp \qquad \leq (3 \eta \sin \vi + 1)^N - 1 \leq C(N)\, \eta \sin \vi < \frac{1}{2} \sin \vi \leq \frac{1}{2}
\ea
\ee
and, furthermore,
\be \label{A.12}
\bigg|\frac{g_N^2 - 1}{g_N^2 + 1}\bigg| \leq \frac{(1/2) \sin \vi}{2 - 1/2} = \frac{1}{3} \sin \vi \leq \frac{1}{3}\, .
\ee
Expressing now $w/v$ by $g_N$, using at that (\ref{A.1}) and (\ref{A.3}), one obtains
\be \label{A.13}
\ba{l}
\disp \frac{w}{v} = - i\, \frac{(\hb /h)^2 - 1}{(\hb/ h)^2 + 1} = - i\, \frac{e^{4 i \roh\, \vi} g_N^2 - 1}{e^{4 i \roh\, \vi} g_N^2 + 1} \\[3ex]
\disp \qquad = \frac{(g_N^2 + 1) \sin 2 \roh\, \vi - i \, (g_N^2 - 1) \cos 2 \roh\, \vi }{ (g_N^2 + 1) \cos 2 \roh\, \vi + i \, (g_N^2 - 1) \sin 2 \roh\, \vi} = \frac{\tan 2 \roh\, \vi - i\, (g_N^2 - 1)/(g_N^2 + 1)}{1 + i \tan 2 \roh\, \vi\,  (g_N^2 -1)/(g_N^2 + 1)}\, .
\ea
\ee
Combining (\ref{A.4}), (\ref{A.12}), and (\ref{A.13}) yields, finally, the desired estimate for $\Psi$ in $A_3$:
$$
|\Psi| \leq \Big|\frac{w}{v}\Big| \leq \frac{\tan 2 \roh \, \vi + (1/3) \sin \vi}{1 - (1/3) \tan 2 \roh\, \vi} \leq \frac{3}{2} \Big(2 \sqrt{2} \, \roh + \frac{1}{3}\Big) \sin \vi\, ,
$$
where in the last estimate we made use of (\ref{A.9a}).
\qed
\sect{Harmonic interpolation of the boundary data}
The solution of the Dirichlet problem in bounded domains with regular boundary for continuous boundary data is a standard problem (see, e.g., Gilbarg and Trudinger 1998, p.\ 26). The only not quite standard issue is antisymmetry and some estimates related to it. We summarize useful properties in the following proposition. 
\begin{prop}
 Let $\phi$ and $\phih : [-\pi , \pi] \ra \real$ be continuous, antisymmetric functions with $\phi (-\pi) = \phi (\pi)$ and axis condition 
 \be \label{B.1}
 \phi (\vi) = O (\vi)\quad \mbox{ for }\  \vi \ra 0\, , \qquad \phi (\pi - \vi)  = O (\vi) \quad  \mbox{for }\  \vi \searrow 0\, ,
 \ee
 and analogously for $\phih$. Then there is a unique antisymmetric solution $\Phi \in C^2 (A_R) \cap C(\overline{A_R})$ of the Dirichlet problem\footnote{Recall that tilde denotes dependence on polar coordinates: $\Pht (r, \vi) := \Phi (r \cos \vi, r \sin \vi)$.} 
\be \label{B.2}
\left.\ba{c}
\Delta \Phi = 0 \qquad \mbox{ in }\, A_R\, ,\\[1ex]
\Pht (1, \cdot)  = \phi\, ,\qquad \Pht (R, \cdot) = \phih\, ,
\ea
\right\}
\ee
which, moreover, satisfies the estimates 
 \be \label{B.3}
 |\Pht ( \cdot , \vi) |\leq K\, |\sin \vi|\qquad \mbox{in } \ A_R \, ,
 \ee 
 \be \label{B.4}
 |\Phi ( \cdot , \rho)/\rho |\leq K \qquad \mbox{in } \ A_R \, ,
 \ee
 \be \label{B.5}
 \int_{A_R} \Phi^2\, \rho^{-(2 + \beta)}\, \rd \ze \rd \rho \leq \frac{2 \pi K^2}{\beta (1 - \beta)}\, , \qquad 0 < \beta < 1 
  \ee
  for some $K> 0$, which does not depend on $R$.
 \end{prop}
\textsc{Proof:} Let $\Phi^+$ be the unique solution of the Dirichlet problem 
$$
\ba{c}
\Delta \Phi^+ = 0 \qquad \mbox{ in }\, A_R^+\, ,\\[1ex]
\Pht^+ (1, \vi )  = \phi (\vi)\, ,\quad \Pht^+ (R, \vi) = \phih( \vi) \quad \mbox{for } \, 0 \leq \vi \leq \pi\, , \\[1ex]
\Pht^+ (r, 0 )  = \Pht^+ (r, \pi) = 0 \quad \mbox{for } \, 1 \leq r \leq R
\ea
$$
in the half-annulus $A_R^+ := A_R \cap \{\rho > 0\}$ with continuous boundary data on the regular boundary $\pa A_R^+$. By antisymmetric continuation onto $A_R$ one obtains a function that is harmonic in $A_R^+ \cup A_R^-$ and, by means of the mean-value characterization of harmonic functions, in fact in $A_R$. This is the asserted antisymmetric solution of (\ref{B.2}).

The bound (\ref{B.3}) holds by (\ref{B.1}) on the boundaries $S_1$ and $S_R$ with some constant $K> 0$. To extend this bound onto $A_R$ consider the comparison function 
$$
\Phb^+ (r, \vi) := \frac{K}{1 + R}\, \Big( r + \frac{R}{r} \Big) \sin \vi\, ,
$$
which satisfies the Dirichlet problem
$$
\ba{c}
\Delta \Phb^+ = 0 \qquad \mbox{ in }\, A_R^+\, ,\\[1ex]
\Phb^+ (1, \vi )  =  \Phb^+ (R, \vi) = K \sin\vi \quad \mbox{for } \, 0 \leq \vi \leq \pi\, , \\[1ex]
\Phb^+ (r, 0 )  = \Phb^+ (r, \pi) = 0 \quad \mbox{for } \, 1 \leq r \leq R\, .
\ea
$$
Applying the maximum principle to $\Pht^+ \mp \Phb^+$ then yields 
$$
|\Pht^+ (r, \vi)| \leq  \frac{K}{1 + R}\, \Big( r + \frac{R}{r} \Big) \sin \vi \leq K \sin \vi \qquad  \mbox{in }\, A_R^+ \, ,
$$
which is (\ref{B.3}). 
(\ref{B.4}) and (\ref{B.5}) follow from (\ref{B.3}) as in appendix A.
\qed

Analogous results hold in the unbounded case, which by Kelvin's transformation (see, e.g., Axler et al., p.\ 54f)
can be reduced to a bounded problem.
\begin{prop}
 Let $\phi : [-\pi , \pi] \ra \real$ be a continuous, antisymmetric function with $\phi (-\pi) = \phi (\pi)$ and axis condition (\ref{B.1}).
 Then there is a unique antisymmetric solution $\Phi \in C^2 (A_\infty) \cap C(\overline{A_\infty})$ of the exterior Dirichlet problem 
$$
\ba{c}
\Delta \Phi = 0 \qquad \mbox{ in }\, A_\infty\, ,\\[1ex]
\Pht (1, \cdot)  = \phi\, ,\qquad |\Pht (r, \cdot)| \ra 0 \quad \mbox{for }\, r \ra \infty\, ,
\ea
$$
which, moreover, satisfies the estimates (\ref{B.3}) -- (\ref{B.5}) (with $A_R$ replaced by $A_\infty$) for some $K>0$.
\end{prop}
\textsc{Proof:} The most simple way to obtain a solution of the asserted kind is to solve again a bounded problem, viz., 
$$
\ba{c}
\Delta \Phi'^+ = 0 \qquad \mbox{ in }\, B_1^+\, ,\\[1ex]
\Pht'^+ (1, \vi )  = - \phi (\vi) \quad \mbox{for } \, 0 \leq \vi \leq \pi\, , \\[1ex]
\Pht'^+ (r, 0 )  = \Pht'^+ (r, \pi) = 0 \quad \mbox{for } \, 0 \leq r \leq 1\, ,
\ea
$$
to continue the solution $\Phi'^+$ antisymmetrically onto $B_1$, and to reflect $\Phi'$ by Kelvin's transformation into $A_\infty$. The estimates (\ref{B.3}) -- (\ref{B.5}) are proved as in the bounded case, where $\overline{\Phi}'^+ := (K/r) \sin \vi$ is a suitable comparison function in $A_\infty^+$.
\qed
\sect{Proof of the generalized Lax-Milgram criterion}
The proof of proposition 5.1 is an exercise in standard functional analysis and is given here for completeness. 

Fixing some $u \in \BC$ in the bilinear form $B(\cdot, \cdot)$ yields a linear form on $\BtC$. The mapping 
$$
A: \BC \ra \BtC' \, , \qquad u \mapsto B(u, \cdot)
$$
is thus well defined with $Au$ acting on $\BtC$ according to 
\be \label{C.1}
\langle A u, \ut\rangle  = B(u, \ut) \quad \mbox{for any }\, \ut \in \BtC\, .
\ee
This mapping is bounded by (\ref{C.1}), (\ref{5.6}), and (\ref{5.5}): 
$$
\|A u\|_{\BtC'} = \sup_{0 \neq \ut \in \BtC} \frac{B(u, \ut)}{\|\ut\|_{\BtC}} \leq K\, \|u\|_\BC \, ,
$$
and injective by (\ref{5.7}):
\be \label{C.2}
\|A u\|_{\BtC'} = \sup_{0 \neq \ut \in \BtC} \frac{B(u, \ut) }{\|\ut\|_{\BtC}} \geq c\, \|u\|_\BC\, .
\ee
The range $A(\BC)$ of $A$ is thus a closed set. 

In order to prove $A(\BC) = \BtC'$ let us assume $A(\BC) \subsetneqq \BtC'$. According to a well-known consequence of Hahn-Banach's principle there is a nontrivial functional $F_0 \in (\BtC')'$ vanishing  on $A(\BC)$ (see, e.g., Kato 1976, p.\ 135). By reflexivity of $\BtC$, $F_0$ is associated to some $\ut_0 \in \BtC$ such that   
\be \label{C.3}
\langle F_0 , f\rangle  = \langle f, \ut_0\rangle \qquad \mbox{ for any } \, f \in \BtC'\, .
\ee
By $\| F_0 \|_{\BtC''} \neq 0$ and (\ref{C.3}) we find for $\ut_0$:
\be \label{C.4}
0 < \|F_0\|_{\BtC''} = \sup_{0 \neq f \in \BtC'} \frac{\langle F_0 , f\rangle}{\|f\|_{\BtC'}} = \sup_{0 \neq f \in \BtC'} \frac{\langle f , \ut_0\rangle}{\|f\|_{\BtC'}} \leq \|\ut_0\|_{\BtC}\, .
\ee
On the other hand, by $F_0\big|_{A(\BC)} = 0$ and (\ref{C.1}) we have for any $u \in \BC$:
$$
0 = \langle F_0, A u \rangle = \langle A u, \ut_0\rangle = B(u, \ut_0)\, ,
$$
and hence by (\ref{5.8}):
$$
\|\ut_0\|_{\BtC} = 0\, ,
$$
which contradicts (\ref{C.4}).

$A: \BC \ra \BtC'$ is thus an isomorphism, which guarantees the unique solvability of the equation $$
A u = f\, ,
$$
which by (\ref{C.1}) is equivalent to eq.\ (\ref{5.9}). Finally, the estimate (\ref{5.10}) follows from (\ref{C.2}).
\qed
\sect{The $x$-independent min-max problem}
We give here a direct (i.e.\ without use of Euler-Lagrange techniques) solution of the min-max problem (\ref{6.1}), (\ref{6.2}) corroborating the result (\ref{6.2.4}). 

In the $x$-independent case problem (\ref{6.1}), (\ref{6.2}) in the interval $(b_0, b_1)$ reads
\be \label{D.0}
\inf_{v \neq 0}\ \sup_{\psi \neq 0}\, \FC_0 [v, \psi] =: \mu_0
\ee
with 
\be \label{D.1}
\FC_0 [v, \psi] := \frac{\disp \int_{b_0}^{b_1} v'\, \psi'\, y^{-\al}\, \rd y}{\disp \Big(\int_{b_0}^{b_1} v'^{\, 2}\, y^{-(\al + \ga)}\, \rd y \Big)^{1/2} \Big(\int_{b_0}^{b_1} \psi'^{\, 2}\, y^{-(\al - \ga)}\, \rd y \Big)^{1/2} }
\ee
and 
\be \label{D.2}
y^{-(\al + \ga)/2}\, v \in H_0^1 ((b_0,b_1)) \, , \qquad y^{-(\al - \ga)/2}\, \psi  \in H_0^1 ((b_0,b_1)) \, .
\ee
(\ref{D.1}) may be simplified by the substitution $(y/b_1)^{1 + \al} =: z$; $\FC_0$ then takes the form
\be \label{D.3}
\FCt_0 [v, \psi] := \frac{\disp \int_{B}^{1} v'\, \psi'\, \rd z}{\disp \Big(\int_{B}^{1} v'^{\, 2}\, z^{-\gt}\, \rd z \Big)^{1/2} \Big(\int_{B}^{1} \psi'^{\, 2}\, z^{\gt}\, \rd z \Big)^{1/2} }
\ee
with $\gt := \ga/(1 + \al)$ and $B:= (b_0 /b_1)^{1 + \al}$, and prime denoting derivation with respect to $z$.

To determine the maximizing $\psi$ for given $v$ let us estimate the numerator in (\ref{D.3}) in the form
\be \label{D.4}
\int_B^1 v' \psi' \, \rd z = \int_B^1 (v' + C)\, \psi'\, \rd z \leq \bigg( \int_B^1 (v' + C)^2\, z^{-\gt}\, \rd z \bigg)^{\frac{1}{2}} \bigg( \int_B^1 \psi'^{\, 2} \, z^\gt \, \rd z \bigg)^{\frac{1}{2}} ,
\ee
where we made use of the boundary conditions on $\psi$. The maximum is reached if (\ref{D.4}) is an equality, i.e.
$$
z^{\gt/2}\, \psi' = \la \, z^{-\gt/2} (v' + C)
$$
with some $\la \in \real$. Given $v$ with $z^{- \gt/2} \, v \in H_0^1 ((B, 1))$, which corresponds to (\ref{D.2})$_a$, the choice
\be \label{D.5}
\psi[v](z) := \la \int_z^1 (v' + C)\, \zt^{- \gt} \, \rd \zt \ ,\qquad C:= -\frac{1 - \gt}{1 - B^{1 - \gt}} \int_B^1 v' z^{- \gt}\, \rd z
\ee
is admissible, i.e.\ $z^{ \gt/2} \, \psi[v] \in H_0^1 ((B, 1))$.\footnote{In the case that $B=0$, $\gt < 1$ is assumed.} Inserting (\ref{D.5}) into (\ref{D.4}) yields 
$$
\ba{l}
\disp \sup_{\psi \neq 0}\, \frac{\disp \int_B^1 v' \psi' \, \rd z}{\disp \Big(\int_B^1 \psi'^{\, 2} \, z^\gt\, \rd z\Big)^{1/2}} = 
\frac{\disp \int_B^1 v' \psi[v]' \, \rd z}{\disp \Big(\int_B^1 \psi[v]'^{\, 2} \, z^\gt\, \rd z\Big)^{1/2}} 
 = \bigg(\int_B^1 (v' + C)^2\, z^{-\gt} \, \rd z \bigg)^{\frac{1}{2}} \\[3ex]
 \disp \qquad \qquad = \bigg(\int_B^1 v'^{\, 2} z^{-\gt} \, \rd z 
 -\frac{1 - \gt}{1 - B^{1 - \gt}} \Big(\int_B^1 v' z^{- \gt}\, \rd z\Big)^2
\bigg)^{\frac{1}{2}}
\ea
$$
and hence
\be \label{D.6}
\ba{l}
\disp \inf_{v \neq 0}\, \sup_{\psi \neq 0}\, \FCt_0 [v, \psi] = \inf_{v \neq 0} \Bigg(1 -\frac{1 - \gt}{1 - B^{1 - \gt}}\, \frac{\disp \Big(\int_B^1 v' z^{- \gt}\, \rd z\Big)^2}{\disp \int_B^1 v'^{\, 2} z^{-\gt}\, \rd z} \Bigg)^{\frac{1}{2}} \\[3ex]
\disp \qquad \qquad = \Bigg(1 -\frac{1 - \gt}{1 - B^{1 - \gt}}\ \sup_{v \neq 0}\, \frac{\disp \Big(\int_B^1 v' z^{- \gt}\, \rd z\Big)^2}{\disp \int_B^1 v'^{\, 2} z^{-\gt}\, \rd z} \Bigg)^{\frac{1}{2}} .
\ea
\ee
Similarly, using the boundary conditions on $v$ we can estimate
\be \label{D.7}
\ba{l}
\disp \int_B^1 v' z^{-\gt} \, \rd z = \int_B^1 v' z^{-\gt/2} \big( z^{-\gt/2} - C\, z^{\gt/2} \big) \rd z \\[2ex]
\disp \qquad \qquad \leq \bigg( \int_B^1 v'^{\, 2} z^{-\gt}\, \rd z \bigg)^{\frac{1}{2}} \bigg(\frac{1 - B^{1- \gt}}{1 - \gt} - 2 C (1 - B) + C^2\, \frac{1 - B^{1 + \gt}}{1 + \gt}  \bigg)^{\frac{1}{2}} .
\ea
\ee
Choosing $C := (1 - B)(1 + \gt)/(1 - B^{1 - \gt})$, (\ref{D.7}) implies
\be \label{D.8}
 \sup_{v \neq 0}\, \frac{\disp \Big(\int_B^1 v' z^{- \gt}\, \rd z\Big)^2}{\disp \int_B^1 v'^{\, 2} z^{-\gt}\, \rd z} \leq \frac{1 - B^{1 - \gt}}{1 - \gt} - (1 -B)^2\, \frac{1 + \gt}{1 - B^{1 + \gt}}\, ,
 \ee
 and the admissible function
 $$
 v_0 := 1 - z - \frac{1 - B}{ 1 - B^{1 + \gt}}\, (1 - z)^{1 + \gt}
 $$
 demonstrates that inequality (\ref{D.8}) is in fact an equality. Combining (\ref{D.0}), (\ref{D.6}), and (\ref{D.8}) we thus obtain 
 $$
 \mu_0^2 = 1 - \frac{1 - \gt}{1 - B^{1 - \gt}} \Big( \frac{1 - B^{1 - \gt}}{1 - \gt} - (1 -B)^2\, \frac{1 + \gt}{1 - B^{1 + \gt}} \Big) = (1 - B)^2\, \frac{1 - \gt}{1 - B^{1 - \gt}}\, \frac{1 + \gt}{1 - B^{1 + \gt}}\, ,
 $$
 which is (\ref{6.2.2}).
\sect{Numerical solution of the Euler-Lagrange equations associated to the min-max problem}
In subsection 6.2 the product ansatz (\ref{6.1.11}), (\ref{6.1.12}) in the rectangular domain $Q$ allowed the reduction of the two-dimensional min-max problem to a sequence of one-dimensional problems labeled by the index $n \in \nat$.
The corresponding one-dimensional Euler-Lagrange equations are obtained by inserting (\ref{6.1.11}) into (\ref{6.1.5}). Using, moreover, the symmmetric variables $w$ and $\chi$ as given by (\ref{6.2.1}) one obtains the following second-order system of ordinary differential equations on the interval $(0,b)$:
\be \label{E.1}
\left.\ba{l}
\disp \chi_n'' - \frac{\ga}{y}\, \chi_n' - \Big(\frac{(\al + \ga + 2)(\al - \ga)}{4\, y^2} + k_n^2\Big) \chi_n \\[2ex]
\disp \qquad \qquad \qquad - \mu_n \Big[w_n'' - \Big(\frac{(\al + \ga +2)( \al + \ga)}{4\, y^2} + d^2 k_n^2\Big) w_n \Big] = 0\, ,\\[3ex]
\disp w_n'' + \frac{\ga}{y}\, w_n' - \Big(\frac{(\al - \ga + 2)(\al + \ga)}{4\, y^2} + k_n^2\Big) w_n \\[2ex]
\disp \qquad \qquad \qquad - \mu_n \Big[\chi_n'' - \Big(\frac{(\al - \ga +2)( \al - \ga)}{4\, y^2} + e^2 k_n^2\Big) \chi_n \Big] = 0
\ea \right\}
\ee
together with the boundary conditions
\be \label{E.2}
w_n(0) = \chi_n (0) = 0\, ,\qquad w_n(b) = \chi_n (b) = 0
\ee
and eigenvalue $\mu_n$. 

For a numerical treatment of system (\ref{E.1}) two more transformations are appropriate. Introducing the variable $z:= k_n y$ yields the system
\be \label{E.3}
\left.\ba{l}
\disp \chi'' - \frac{\ga}{z}\, \chi' - \Big(\frac{(\al + \ga + 2)(\al - \ga)}{4\, z^2} + 1 \Big) \chi \\[2ex]
\disp \qquad \qquad \qquad - \mu \Big[w'' - \Big(\frac{(\al + \ga +2)( \al + \ga)}{4\, z^2} + d^2 \Big) w \Big] = 0\, ,\\[3ex]
\disp w'' + \frac{\ga}{z}\, w' - \Big(\frac{(\al - \ga + 2)(\al + \ga)}{4\, z^2} + 1 \Big) w \\[2ex]
\disp \qquad \qquad \qquad - \mu \Big[\chi'' - \Big(\frac{(\al - \ga +2)( \al - \ga)}{4\, z^2} + e^2 \Big) \chi \Big] = 0 \, ,
\ea \right\}
\ee
where we omitted the index $n$ at $\chi$, $w$, and $\mu$ and 
where prime denotes differentiation with respect to the variable $z \in [0, k_n b]$. Note that the parameters $a$, $b$, and $n$ are now combined into the single parameter 
$$
\bti := k_n b = n\, \frac{\pi}{2}\, \frac{b}{a}\, ,
$$
which limits the interval $[0, \bti]$ of integration in (\ref{E.3}). Finally, the change of variables $f_\pm := w \pm \chi$ separates the highest ($\triangleq $ second) derivatives; in these variables (\ref{E.3}) takes the form
\be \label{E.4}
\left.\ba{l}
\disp (1- \mu) f_+'' = (1- \mu)\, \frac{\al (\al +2)}{4\, z^2} \, f_+ - (1 + \mu)\, \frac{\ga^2}{4 \, z^2} \, f_+ + \Big(1 - \frac{\mu}{2}\, (d^2+ e^2)\Big) f_+ \\[2ex]
\disp \qquad \qquad \qquad - \frac{\ga}{z}\, f_-' + \big(1 - \mu(\al + 1)\big) \frac{\ga}{2\, z^2} \, f_- - \frac{\mu}{2}\, (d^2 -e^2)  f_- \, ,\\[3ex]
\disp (1+ \mu) f_-'' = (1+ \mu)\, \frac{\al (\al +2)}{4\, z^2} \, f_- - (1 - \mu)\, \frac{\ga^2}{4 \, z^2} \, f_- + \Big(1 + \frac{\mu}{2}\, (d^2+ e^2)\Big) f_- \\[2ex]
\disp \qquad \qquad \qquad - \frac{\ga}{z}\, f_+' + \big(1 + \mu(\al + 1)\big) \frac{\ga}{2\, z^2} \, f_+ + \frac{\mu}{2}\, (d^2 -e^2)  f_+ \, .
\ea \right\}
\ee
Equations (\ref{E.4}) together with the boundary conditions 
\be \label{E.5}
f_\pm (b_0) = f_\pm (\bti) = 0
\ee
are solved in the interval $[b_0, \bti]$ by a shooting method. To this end trajectories are calculated starting from both ends of the interval towards the center for a given value of $\mu$. The integration is done with the standard Runge-Kutta-Fehlberg method, consisting of a fourth-order Runge-Kutta solver with fifth-order error prediction to control the stepsize. The matching of the trajectories at the center of the interval amounts to the vanishing of a $4\ti 4$ determinant computed from the trajectories at the center. 
The search for zeroes of this determinant depending on $\mu$ is done by an interval method with respect to $\mu$: we scan systematically the interval $[\varepsilon, 1]$ for subintervals exhibiting a sign change at the endpoints and refine in the positive case the subinterval subsequently. In order not to miss higher-order zeroes we looked, additionally, for minima of the modulus of the determinant. These zeroes correspond to the eigenvalues of $\mu$ and, in particular, the lowest one corresponds to $\mu_{min}$.
\begin{figure}
\begin{center}
\input{FigE1}
\caption{$\mu_0$ versus $\ga$ for $\al = 0.39$ and $b_0/b_1 = 0$, $10^{-6}$ (cross), $10^{-4}$ (square), $10^{-2}$ (asterisk). Solid lines denote analytical curves, symbols denote numerical values.}
\end{center}
\end{figure}
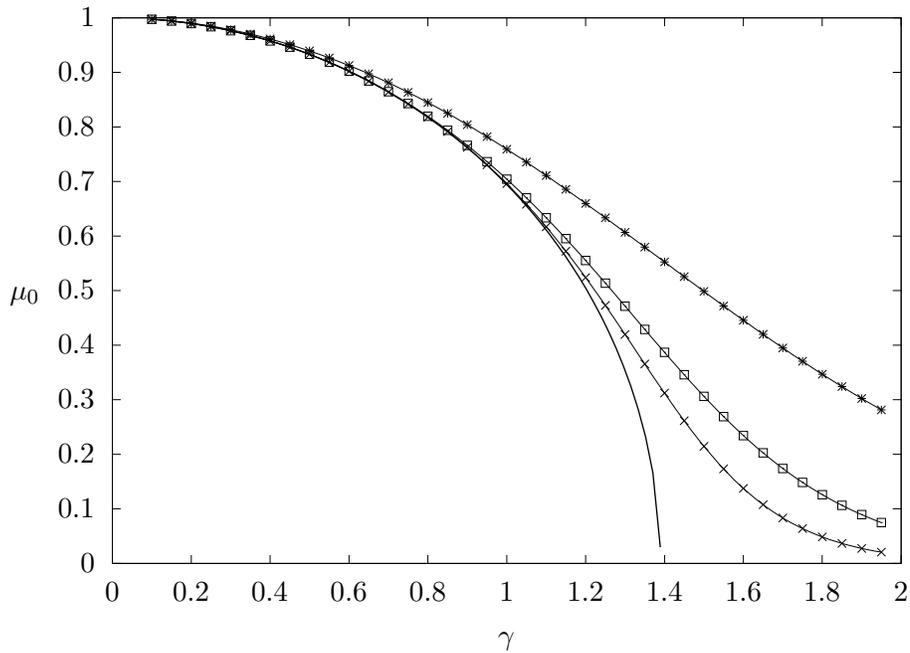 

To validate the numerical code let us first consider the  $x$-independent case as described by eqs.\ (\ref{6.2.2}). These can be obtained from (\ref{E.4}) just by dropping the coefficients $1 \mp \frac{\mu}{2} (d^2 + e^2)$ and $\frac{\mu}{2} (d^2- e^2)$. Figure E.1 demonstrates good agreement of this way numerically determined points in a $\ga$--$\mu_0$--diagram with curves $\ga \mapsto \mu_0$ as described by the analytical result (\ref{6.2.4}) for $\al = 0.39$ and three representative ratios $b_0 /b_1$. 
For all three ratios the deviation of the numerical data from the analytical curve is beyond the resolution of the plot. 

The full equations (\ref{E.4}) with boundary conditions (\ref{E.5}) have been run with $\al = 0.39$, $e =1.2$, and a variety of values for $d$, $b_0$, and $\bti$. With respect to $\bti$ we find $\mu_{min}$ to be a monotonically decreasing function.
Figure E.2 displays results in a $\ga$--$\mu_{min}$--diagram for $d= 0.1$, $b_0= 10^{-3}$, and $\bti = 10^2$ together with the comparison function $C\!F(0.39, \ga)$ and the lower bound $L\!B (0.39, \ga, 0.1, 1.2)$.
For $\ga \gtrsim 0.6$, $\mu_{min}$ is above $L\!B$, which means an improvement of the lower bound $L\!B$. 
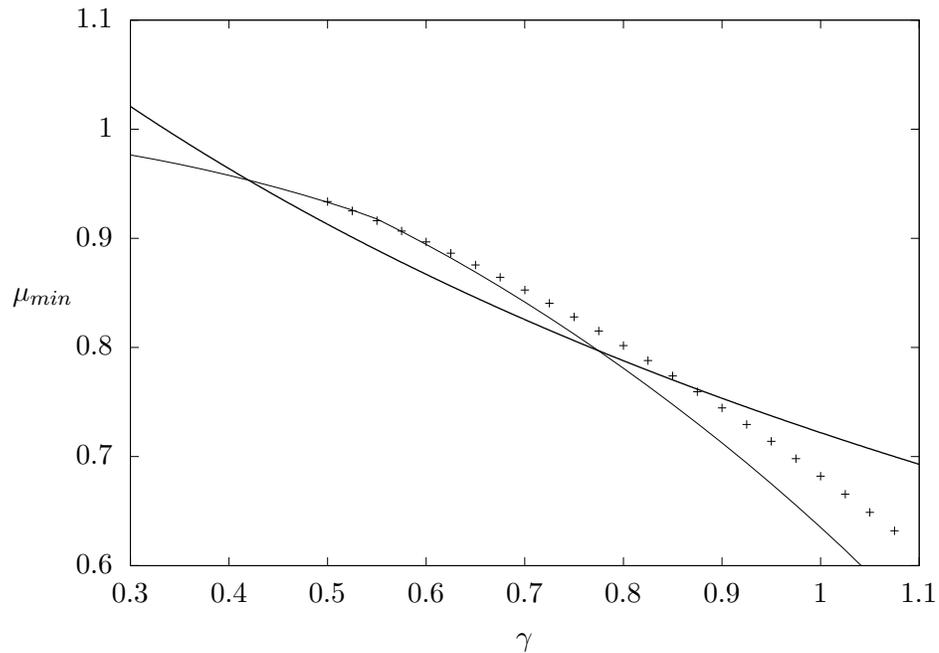
\begin{figure}
\begin{center}
\input{FigE2}
\caption{$\mu_{min}$ versus $\ga$ (crosses) together with the comparison function $C\!F(\al, \ga)$ (thick line) and the lower bound $L\!B (\al, \ga, d,e)$ (thin line) for $\al = 0.39$, $d = 0.1$, $e = 1.2$, $b_0 = 10^{-3}$, and $\bti = 10^2$.}
\end{center}
\end{figure}
Further decreasing $b_0$ or increasing $\bti$ does not lower $\mu_{min}$ in a significant way, which indicates that $(10^{-3}, 10^2)$ is a sufficient approximation of the interval $(0, \infty)$. For $0.5 \lesssim \mu_{min} \lesssim 0.6$, $\mu_{min}$ and $L\!B$ coincide within numerical accuracy. No results are shown for $\ga$ below $0.5$ because of numerical instability of our code in this range (a problem that we did not further pursue). 

Let us finally recall that in the case that $\mu_{min}$ is unique, which holds in the $x$-independent case but has not been proven in the general case, $\mu_{min}$ represents the (in the limits of numerical accuracy) exact value $C_c = \Ct_c$ of the min-max problem. 
\sect{An analytic solution with single zero in two dimensions}
This appendix presents the explicit solution of the direction problem for a simple direction field in two dimensions.\footnote{The two-dimensional field mimics an axisymmetric field in three dimensions, where polynomial solutions of type (\ref{F.2}) are not at our disposal.} The solution contains a single zero, whose position is governed by a parameter $\la$. The example illustrates the migration of this position with $\la$ and, in particular, the splitting of the zero when hitting the boundary.
\begin{figure}
\begin{center}
\includegraphics[width=0.4\textwidth]{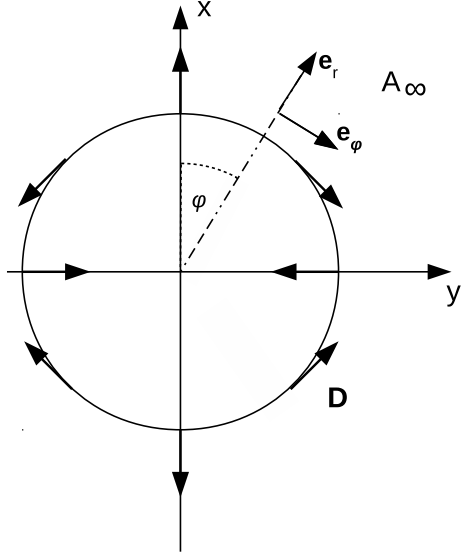}
\caption{The two-dimensional direction field (\ref{F.1}) with rotation number $\ro = 3$.}
\end{center}
\end{figure}
Let $(r, \phi)$ be polar coordinates in $\real^2$ with basis vectors $\eee_r$ and $\eee_\vi$, and $A_\infty$ the exterior plane with boundary $S_1$. Furthermore, let 
\be \label{F.1}
\DD := \cos 2 \vi \, \eee_r + \sin 2 \vi \, \eee_\vi
\ee
be a (symmetric) direction field on $S_1$ (see Fig.\ F.1) and 
\be \label{F.2}
\ba{l}
\disp \BBt_\la (r,\vi) := \na \Upsilon_\la (r, \vi) := \na \bigg( 3 \la\, \frac{\cos 2 \vi}{2 r^2} - (1-\la^2) \Big(\frac{\cos \vi}{r} + \frac{\cos 3 \vi}{3 r^3} \Big)\bigg) \\[3ex]
\disp \qquad \qquad =\Big[(1-\la^2) \Big(\frac{\cos \vi}{r^2} + \frac{\cos 3 \vi}{ r^4} \Big) - 3 \la\, \frac{\cos 2 \vi}{ r^3}\Big] \eee_r \\[2ex]
\disp \qquad \qquad \qquad + \Big[(1-\la^2) \Big(\frac{\sin \vi}{r^2} + \frac{\sin 3 \vi}{ r^4} \Big) - 3 \la\, \frac{\sin 2 \vi}{ r^3}\Big] \eee_\vi 
\ea
\ee
a vector field on $A_\infty$.

An elementary calculation shows that 
$$
\De \Upsilon_\la = 0 \qquad \mbox{ in }\ A_\infty
$$
and 
$$
\ba{l}
\disp \BBt_\la (1, \vi)  = \big[(1-\la^2) (\cos \vi + \cos 3 \vi) - 3 \la \cos 2 \vi \big] \eee_r \\[2ex]
\disp \qquad \qquad \qquad + \big[(1-\la^2) (\sin \vi + \sin 3 \vi) - 3 \la \sin 2 \vi \big] \eee_\vi \\[3ex]
\disp \qquad  =\big[2(1-\la^2) \cos \vi  \cos 2 \vi - 3 \la \cos 2 \vi \big] \eee_r 
+ \big[2(1-\la^2) \cos \vi  \sin 2 \vi - 3 \la \sin 2 \vi \big] \eee_\vi \\[2ex]
\disp \qquad \qquad \qquad = \big(2(1 - \la^2) \cos \vi - 3 \la \big) \DD =: a_\la (\vi)\, \DD\, . 
\ea
$$
$\BBt_\la$ is thus a harmonic vector field in $A_\infty$ that solves the (unsigned) direction problem.

\begin{figure}
\begin{center}
\includegraphics[width=0.4\textwidth]{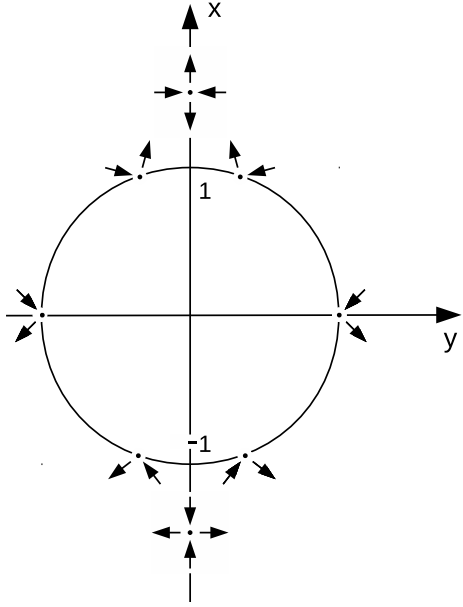}
\caption{Zeroes of $\BBt_\la$ for 5 different values of $\la$. Arrows at the zeroes indicate in- and outgoing field lines.}
\end{center}
\end{figure}
Let us consider the parameter range $\la \in (-1, 1)$ , where $\BBt_\la$ has the exact decay order $\dt = 2$. As long as $a_\la \neq 0$ the 2D-analogue of eq.\ (\ref{2.7}) predicts $\ro -\dt = 3-2 =1$ zero of $\BBt_\la$ in $A_\infty$, where $\ro$ is the rotation number of $\DD$ (see Kaiser 2010). From (\ref{F.2}) one obtains explicitly for $z_0 = x_0 + i y_0 = r_0 (\cos \vi_0 + i \sin \vi_0)$:
$$
(r_0, \vi_0) = \left\{ 
\ba{lc}
( \mu + \sqrt{\mu^2 - 1}\, , 0) & \mu >1 \\[1ex]
(1\, ,\, \arccos \mu ) &-1 \leq \mu \leq 1 \\[1ex]
(-\mu + \sqrt{\mu^2 - 1}\, , \pi) & \mu < -1
\ea \right. ,
$$
where $\mu := 3 \la /(2 - 2 \la^2)$. With $\la$ running through the interval $[-1, 1]$ from $1$ to $-1$ one finds the following behaviour of $z_0$ (see Fig.\ F.2): $z_0$ moves from $+ \infty$ ($\la = 1$) along the $x$-axis, hits the boundary at $x =1$ ($\la = 1/2$), splits into two ``boundary-zeroes'', which move symmetrically along $S_1$ and coalesce
again at $x =-1$ ($\la = -1/2$); for $\la < -1/2$, $z_0$ enters again $A_\infty$ and moves along the negative $x$-axis up to $- \infty$ ($\la = -1$). 

Note that eq.\ (\ref{2.7}) remains valid in the presence of boundary-zeroes, if these are included on the left-hand side, however weighted by a factor $1/2$.
%
%
%
%
%

\section*{Acknowledgement}
The authors would like to thank A.\ Tilgner for substantial help in preparing the numerical code for the computations in appendix E.  
%
%
%
%
%
%
%
%
%
\section*{References}


Axler, S., Bourdon, P., Ramey, W., {\it Harmonic function theory}. Springer, New York 1992.\\[1ex]
Bers, L., John, F., Schechter, M., {\it Partial Differential Equations}. 
Interscience Publishers, New York 1964.\\[1ex]
Bers, L., Nirenberg, L., {\it On a representation theorem for linear elliptic systems with discontinuous coefficients and its applications}. Atti del convegno internazionale sulle equazioni alle derivate parzialli, pp. 111--140, Trieste 1954.\\[1ex]
Courant, R., Hilbert, D., {\it Methods of mathematical physics} vol.\ I. Interscience, New York 1953.\\[1ex]
Folland, G. B., {\it Introduction To Partial Differential Equations}. Princeton University Press,
Princeton, N.Y. 1995.\\[1ex]
Galdi, G. P., {\it An Introduction to the Mathematical Theory of the Navier-Stokes Equations Vol. I. Linearized Steady Problems}. Springer
Tracts in Natural Philosophy Vol. 38, Springer-Verlag, Berlin, Heidelberg, New York 1994.\\[1ex]
Gilbarg, D., Trudinger, N.\ S., {\it Elliptic Partial Differential Equations 
of Second Order}. Springer, Berlin, Heidelberg, New York 2nd edition 1998.\\[1ex]
Gubbins, D., Herrero-Bervera, E. (Eds.), {\it Encyclopedia of geomagnetism and paleomagnetism}, p.\ 466f and p.\ 679f.
Springer, Berlin, Heidelberg, New York 2007. \\[1ex]
Hardy, G. H., Littlewood, J. E. and P$\acute{{\rm o}}$lya, G., {\em Inequalities}. Cambridge University Press 1973. \\[1ex]
Heinz, E., {\it An elementary analytic theory of the degree of mapping in $n$-dimensional space}. Journal of Mathematics and
Mechanics {\bf 8}, 231--247 (1959).\\[1ex]
Hulot, G., Khokhlov, A., Le Mou\"el, J.L., {\it
Uniqueness of  mainly dipolar magnetic fields recovered from directional
data}. Geophys.\ J.\ Int. {\bf 129}, 347--354 (1997). \\[1ex]
Kaiser, R., {\it The geomagnetic direction problem: the two-dimensional and the three-dimensional axisymmetric cases}.
SIAM J.\ Math.\ Anal. {\bf 42}, 701-728 (2010). \\[1ex]
Kaiser, R., {\it Uniqueness and non-uniqueness in the non-axisymmetric direction problem}.
Q.\ Jl.\ Mech.\ Appl.\ Math. {\bf 65}, 347-360 (2012). \\[1ex]
Kaiser, R., Neudert, M., {\it A non-standard boundary value problem related to 
geomagnetism}.  Quarterly of Applied Mathematics {\bf 62}, 423--457 (2004).\\[1ex]
Kato, T., {\it Perturbation Theory for Linear Operators}. Springer, Berlin, Heidelberg, New York 2nd edition 1976.\\[1ex]
Ladyzhenskaya, O. A., Ural'tseva, N. N., {\it Linear and Quasilinear Elliptic Equations}. Academic Press, New York, London 1968.\\[1ex]
Lieberman, G. M., {\it Oblique Derivative Problems for Elliptic Equations}. World Scientific, Singapore 2013.\\[1ex]
Martensen, E., {\it Potentialtheorie}. Teubner-Verlag, Stuttgart 1968. \\[1ex]
Merrill, R. T., McElhinny, M. W., {\it The Earth's Magnetic Field
(Its History, Origin and Planetary Perspective)}. Academic Press, London 1983. \\[1ex]
Paneah, B. P., {\it The Oblique Derivative Problem. The Poincar\'e Problem. Mathematical Topics, vol. 17}. Wiley-VCH-Verlag, Berlin 2000.\\[1ex]
Proctor, M. R. E., Gubbins, D., {\it Analysis of geomagnetic directional data}.
Geophys.\ J.\ Int. {\bf 100}, 69--77 (1990). \\[1ex]

\end{document}

%% file: Fig22cos.tex
\begingroup
  \makeatletter
  \providecommand\color[2][]{%
    \GenericError{(gnuplot) \space\space\space\@spaces}{%
      Package color not loaded in conjunction with
      terminal option `colourtext'%
    }{See the gnuplot documentation for explanation.%
    }{Either use 'blacktext' in gnuplot or load the package
      color.sty in LaTeX.}%
    \renewcommand\color[2][]{}%
  }%
  \providecommand\includegraphics[2][]{%
    \GenericError{(gnuplot) \space\space\space\@spaces}{%
      Package graphicx or graphics not loaded%
    }{See the gnuplot documentation for explanation.%
    }{The gnuplot epslatex terminal needs graphicx.sty or graphics.sty.}%
    \renewcommand\includegraphics[2][]{}%
  }%
  \providecommand\rotatebox[2]{#2}%
  \@ifundefined{ifGPcolor}{%
    \newif\ifGPcolor
    \GPcolorfalse
  }{}%
  \@ifundefined{ifGPblacktext}{%
    \newif\ifGPblacktext
    \GPblacktexttrue
  }{}%
  \let\gplgaddtomacro\g@addto@macro
  \gdef\gplbacktext{}%
  \gdef\gplfronttext{}%
  \makeatother
  \ifGPblacktext
    \def\colorrgb#1{}%
    \def\colorgray#1{}%
  \else
    \ifGPcolor
      \def\colorrgb#1{\color[rgb]{#1}}%
      \def\colorgray#1{\color[gray]{#1}}%
      \expandafter\def\csname LTw\endcsname{\color{white}}%
      \expandafter\def\csname LTb\endcsname{\color{black}}%
      \expandafter\def\csname LTa\endcsname{\color{black}}%
      \expandafter\def\csname LT0\endcsname{\color[rgb]{1,0,0}}%
      \expandafter\def\csname LT1\endcsname{\color[rgb]{0,1,0}}%
      \expandafter\def\csname LT2\endcsname{\color[rgb]{0,0,1}}%
      \expandafter\def\csname LT3\endcsname{\color[rgb]{1,0,1}}%
      \expandafter\def\csname LT4\endcsname{\color[rgb]{0,1,1}}%
      \expandafter\def\csname LT5\endcsname{\color[rgb]{1,1,0}}%
      \expandafter\def\csname LT6\endcsname{\color[rgb]{0,0,0}}%
      \expandafter\def\csname LT7\endcsname{\color[rgb]{1,0.3,0}}%
      \expandafter\def\csname LT8\endcsname{\color[rgb]{0.5,0.5,0.5}}%
    \else
      \def\colorrgb#1{\color{black}}%
      \def\colorgray#1{\color[gray]{#1}}%
      \expandafter\def\csname LTw\endcsname{\color{white}}%
      \expandafter\def\csname LTb\endcsname{\color{black}}%
      \expandafter\def\csname LTa\endcsname{\color{black}}%
      \expandafter\def\csname LT0\endcsname{\color{black}}%
      \expandafter\def\csname LT1\endcsname{\color{black}}%
      \expandafter\def\csname LT2\endcsname{\color{black}}%
      \expandafter\def\csname LT3\endcsname{\color{black}}%
      \expandafter\def\csname LT4\endcsname{\color{black}}%
      \expandafter\def\csname LT5\endcsname{\color{black}}%
      \expandafter\def\csname LT6\endcsname{\color{black}}%
      \expandafter\def\csname LT7\endcsname{\color{black}}%
      \expandafter\def\csname LT8\endcsname{\color{black}}%
    \fi
  \fi
    \setlength{\unitlength}{0.0500bp}%
    \ifx\gptboxheight\undefined%
      \newlength{\gptboxheight}%
      \newlength{\gptboxwidth}%
      \newsavebox{\gptboxtext}%
    \fi%
    \setlength{\fboxrule}{0.5pt}%
    \setlength{\fboxsep}{1pt}%
\begin{picture}(7200.00,5040.00)%
    \gplgaddtomacro\gplbacktext{%
      \csname LTb\endcsname
      \put(1386,704){\makebox(0,0)[r]{\strut{}\LARGE -1}}%
      \put(1386,1733){\makebox(0,0)[r]{\strut{}\LARGE -0.5}}%
      \put(1386,2762){\makebox(0,0)[r]{\strut{}\LARGE 0}}%
      \put(1386,3790){\makebox(0,0)[r]{\strut{}\LARGE 0.5}}%
      \put(1386,4819){\makebox(0,0)[r]{\strut{}\LARGE 1}}%
      \put(2047,484){\makebox(0,0){\strut{}\LARGE -20}}%
      \put(3104,484){\makebox(0,0){\strut{}\LARGE -10}}%
      \put(4161,484){\makebox(0,0){\strut{}\LARGE 0}}%
      \put(5218,484){\makebox(0,0){\strut{}\LARGE 10}}%
      \put(6275,484){\makebox(0,0){\strut{}\LARGE 20}}%
    }%
    \gplgaddtomacro\gplfronttext{%
      \csname LTb\endcsname
      \put(700,2761){\makebox(0,0){\strut{}\LARGE $a_\ze$}}%
      \put(4160,154){\makebox(0,0){\strut{}\LARGE $x$}}%
    }%
    \gplbacktext
    \put(0,0){\includegraphics{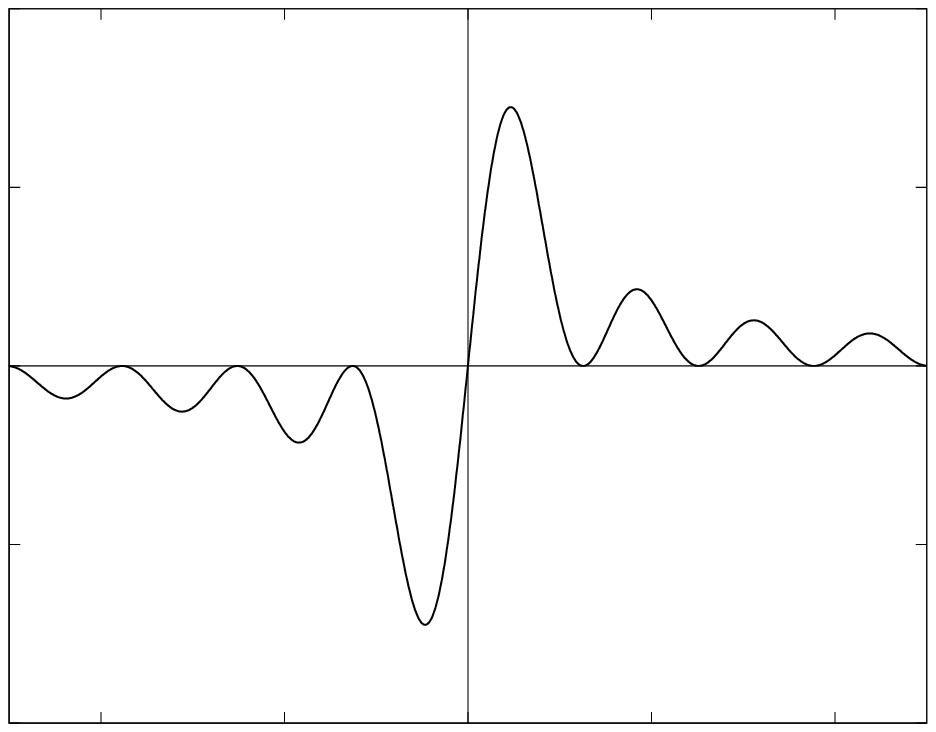}}%
    \gplfronttext
  \end{picture}%
\endgroup

%% file: Fig22sin.tex
\begingroup
  \makeatletter
  \providecommand\color[2][]{%
    \GenericError{(gnuplot) \space\space\space\@spaces}{%
      Package color not loaded in conjunction with
      terminal option `colourtext'%
    }{See the gnuplot documentation for explanation.%
    }{Either use 'blacktext' in gnuplot or load the package
      color.sty in LaTeX.}%
    \renewcommand\color[2][]{}%
  }%
  \providecommand\includegraphics[2][]{%
    \GenericError{(gnuplot) \space\space\space\@spaces}{%
      Package graphicx or graphics not loaded%
    }{See the gnuplot documentation for explanation.%
    }{The gnuplot epslatex terminal needs graphicx.sty or graphics.sty.}%
    \renewcommand\includegraphics[2][]{}%
  }%
  \providecommand\rotatebox[2]{#2}%
  \@ifundefined{ifGPcolor}{%
    \newif\ifGPcolor
    \GPcolorfalse
  }{}%
  \@ifundefined{ifGPblacktext}{%
    \newif\ifGPblacktext
    \GPblacktexttrue
  }{}%
  \let\gplgaddtomacro\g@addto@macro
  \gdef\gplbacktext{}%
  \gdef\gplfronttext{}%
  \makeatother
  \ifGPblacktext
    \def\colorrgb#1{}%
    \def\colorgray#1{}%
  \else
    \ifGPcolor
      \def\colorrgb#1{\color[rgb]{#1}}%
      \def\colorgray#1{\color[gray]{#1}}%
      \expandafter\def\csname LTw\endcsname{\color{white}}%
      \expandafter\def\csname LTb\endcsname{\color{black}}%
      \expandafter\def\csname LTa\endcsname{\color{black}}%
      \expandafter\def\csname LT0\endcsname{\color[rgb]{1,0,0}}%
      \expandafter\def\csname LT1\endcsname{\color[rgb]{0,1,0}}%
      \expandafter\def\csname LT2\endcsname{\color[rgb]{0,0,1}}%
      \expandafter\def\csname LT3\endcsname{\color[rgb]{1,0,1}}%
      \expandafter\def\csname LT4\endcsname{\color[rgb]{0,1,1}}%
      \expandafter\def\csname LT5\endcsname{\color[rgb]{1,1,0}}%
      \expandafter\def\csname LT6\endcsname{\color[rgb]{0,0,0}}%
      \expandafter\def\csname LT7\endcsname{\color[rgb]{1,0.3,0}}%
      \expandafter\def\csname LT8\endcsname{\color[rgb]{0.5,0.5,0.5}}%
    \else
      \def\colorrgb#1{\color{black}}%
      \def\colorgray#1{\color[gray]{#1}}%
      \expandafter\def\csname LTw\endcsname{\color{white}}%
      \expandafter\def\csname LTb\endcsname{\color{black}}%
      \expandafter\def\csname LTa\endcsname{\color{black}}%
      \expandafter\def\csname LT0\endcsname{\color{black}}%
      \expandafter\def\csname LT1\endcsname{\color{black}}%
      \expandafter\def\csname LT2\endcsname{\color{black}}%
      \expandafter\def\csname LT3\endcsname{\color{black}}%
      \expandafter\def\csname LT4\endcsname{\color{black}}%
      \expandafter\def\csname LT5\endcsname{\color{black}}%
      \expandafter\def\csname LT6\endcsname{\color{black}}%
      \expandafter\def\csname LT7\endcsname{\color{black}}%
      \expandafter\def\csname LT8\endcsname{\color{black}}%
    \fi
  \fi
    \setlength{\unitlength}{0.0500bp}%
    \ifx\gptboxheight\undefined%
      \newlength{\gptboxheight}%
      \newlength{\gptboxwidth}%
      \newsavebox{\gptboxtext}%
    \fi%
    \setlength{\fboxrule}{0.5pt}%
    \setlength{\fboxsep}{1pt}%
\begin{picture}(7200.00,5040.00)%
    \gplgaddtomacro\gplbacktext{%
      \csname LTb\endcsname
      \put(1518,704){\makebox(0,0)[r]{\strut{}\LARGE -1}}%
      \put(1518,1733){\makebox(0,0)[r]{\strut{}\LARGE -0.5}}%
      \put(1518,2762){\makebox(0,0)[r]{\strut{}\LARGE 0}}%
      \put(1518,3790){\makebox(0,0)[r]{\strut{}\LARGE 0.5}}%
      \put(1518,4819){\makebox(0,0)[r]{\strut{}\LARGE 1}}%
      \put(2165,484){\makebox(0,0){\strut{}\LARGE -20}}%
      \put(3196,484){\makebox(0,0){\strut{}\LARGE -10}}%
      \put(4227,484){\makebox(0,0){\strut{}\LARGE 0}}%
      \put(5257,484){\makebox(0,0){\strut{}\LARGE 10}}%
      \put(6288,484){\makebox(0,0){\strut{}\LARGE 20}}%
    }%
    \gplgaddtomacro\gplfronttext{%
      \csname LTb\endcsname
      \put(800,2761){\makebox(0,0){\strut{}\LARGE $a_\ro$}}%
      \put(4226,154){\makebox(0,0){\strut{}\LARGE $x$}}%
    }%
    \gplbacktext
    \put(0,0){\includegraphics{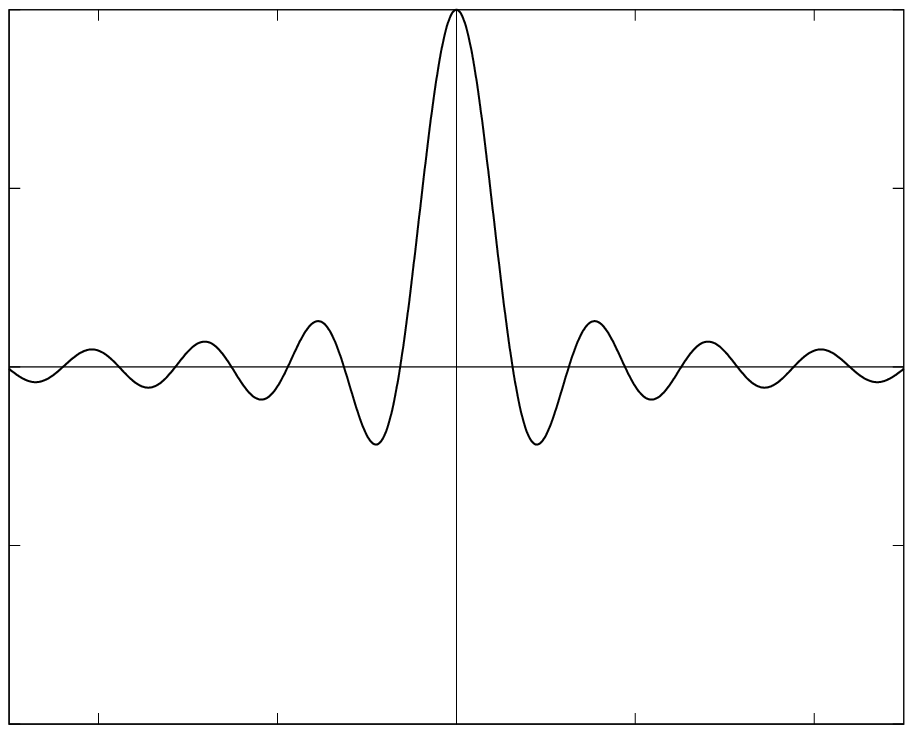}}%
    \gplfronttext
  \end{picture}%
\endgroup

%% file: Fig41.tex
\begingroup
  \makeatletter
  \providecommand\color[2][]{%
    \GenericError{(gnuplot) \space\space\space\@spaces}{%
      Package color not loaded in conjunction with
      terminal option `colourtext'%
    }{See the gnuplot documentation for explanation.%
    }{Either use 'blacktext' in gnuplot or load the package
      color.sty in LaTeX.}%
    \renewcommand\color[2][]{}%
  }%
  \providecommand\includegraphics[2][]{%
    \GenericError{(gnuplot) \space\space\space\@spaces}{%
      Package graphicx or graphics not loaded%
    }{See the gnuplot documentation for explanation.%
    }{The gnuplot epslatex terminal needs graphicx.sty or graphics.sty.}%
    \renewcommand\includegraphics[2][]{}%
  }%
  \providecommand\rotatebox[2]{#2}%
  \@ifundefined{ifGPcolor}{%
    \newif\ifGPcolor
    \GPcolorfalse
  }{}%
  \@ifundefined{ifGPblacktext}{%
    \newif\ifGPblacktext
    \GPblacktexttrue
  }{}%
  \let\gplgaddtomacro\g@addto@macro
  \gdef\gplbacktext{}%
  \gdef\gplfronttext{}%
  \makeatother
  \ifGPblacktext
    \def\colorrgb#1{}%
    \def\colorgray#1{}%
  \else
    \ifGPcolor
      \def\colorrgb#1{\color[rgb]{#1}}%
      \def\colorgray#1{\color[gray]{#1}}%
      \expandafter\def\csname LTw\endcsname{\color{white}}%
      \expandafter\def\csname LTb\endcsname{\color{black}}%
      \expandafter\def\csname LTa\endcsname{\color{black}}%
      \expandafter\def\csname LT0\endcsname{\color[rgb]{1,0,0}}%
      \expandafter\def\csname LT1\endcsname{\color[rgb]{0,1,0}}%
      \expandafter\def\csname LT2\endcsname{\color[rgb]{0,0,1}}%
      \expandafter\def\csname LT3\endcsname{\color[rgb]{1,0,1}}%
      \expandafter\def\csname LT4\endcsname{\color[rgb]{0,1,1}}%
      \expandafter\def\csname LT5\endcsname{\color[rgb]{1,1,0}}%
      \expandafter\def\csname LT6\endcsname{\color[rgb]{0,0,0}}%
      \expandafter\def\csname LT7\endcsname{\color[rgb]{1,0.3,0}}%
      \expandafter\def\csname LT8\endcsname{\color[rgb]{0.5,0.5,0.5}}%
    \else
      \def\colorrgb#1{\color{black}}%
      \def\colorgray#1{\color[gray]{#1}}%
      \expandafter\def\csname LTw\endcsname{\color{white}}%
      \expandafter\def\csname LTb\endcsname{\color{black}}%
      \expandafter\def\csname LTa\endcsname{\color{black}}%
      \expandafter\def\csname LT0\endcsname{\color{black}}%
      \expandafter\def\csname LT1\endcsname{\color{black}}%
      \expandafter\def\csname LT2\endcsname{\color{black}}%
      \expandafter\def\csname LT3\endcsname{\color{black}}%
      \expandafter\def\csname LT4\endcsname{\color{black}}%
      \expandafter\def\csname LT5\endcsname{\color{black}}%
      \expandafter\def\csname LT6\endcsname{\color{black}}%
      \expandafter\def\csname LT7\endcsname{\color{black}}%
      \expandafter\def\csname LT8\endcsname{\color{black}}%
    \fi
  \fi
    \setlength{\unitlength}{0.0500bp}%
    \ifx\gptboxheight\undefined%
      \newlength{\gptboxheight}%
      \newlength{\gptboxwidth}%
      \newsavebox{\gptboxtext}%
    \fi%
    \setlength{\fboxrule}{0.5pt}%
    \setlength{\fboxsep}{1pt}%
\begin{picture}(5668.00,2834.00)%
    \gplgaddtomacro\gplbacktext{%
      \csname LTb\endcsname
      \put(726,704){\makebox(0,0)[r]{\strut{}$0$}}%
      \put(726,1086){\makebox(0,0)[r]{\strut{}$0.2$}}%
      \put(726,1468){\makebox(0,0)[r]{\strut{}$0.4$}}%
      \put(726,1849){\makebox(0,0)[r]{\strut{}$0.6$}}%
      \put(726,2231){\makebox(0,0)[r]{\strut{}$0.8$}}%
      \put(726,2613){\makebox(0,0)[r]{\strut{}$1$}}%
      \put(858,484){\makebox(0,0){\strut{}$0$}}%
      \put(1741,484){\makebox(0,0){\strut{}$1$}}%
      \put(2623,484){\makebox(0,0){\strut{}$2$}}%
      \put(3506,484){\makebox(0,0){\strut{}$3$}}%
      \put(4388,484){\makebox(0,0){\strut{}$4$}}%
      \put(5271,484){\makebox(0,0){\strut{}$5$}}%
    }%
    
    \gplgaddtomacro\gplfronttext{%
      \csname LTb\endcsname
      \put(198,1658){\makebox(0,0){\strut{}$f$}}%
      \put(3064,154){\makebox(0,0){\strut{}$x$}}%
      \put(4550,1908){\makebox(0,0){\strut{}$\sim\ln x$}}%
    }%
    \gplbacktext
    \put(0,0){\includegraphics{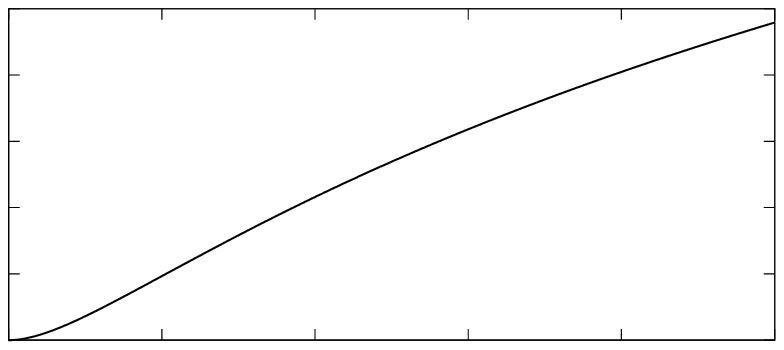}}%
    \gplfronttext
  \end{picture}%
\endgroup

%% file: Fig62.tex
\begingroup
  \makeatletter
  \providecommand\color[2][]{%
    \GenericError{(gnuplot) \space\space\space\@spaces}{%
      Package color not loaded in conjunction with
      terminal option `colourtext'%
    }{See the gnuplot documentation for explanation.%
    }{Either use 'blacktext' in gnuplot or load the package
      color.sty in LaTeX.}%
    \renewcommand\color[2][]{}%
  }%
  \providecommand\includegraphics[2][]{%
    \GenericError{(gnuplot) \space\space\space\@spaces}{%
      Package graphicx or graphics not loaded%
    }{See the gnuplot documentation for explanation.%
    }{The gnuplot epslatex terminal needs graphicx.sty or graphics.sty.}%
    \renewcommand\includegraphics[2][]{}%
  }%
  \providecommand\rotatebox[2]{#2}%
  \@ifundefined{ifGPcolor}{%
    \newif\ifGPcolor
    \GPcolorfalse
  }{}%
  \@ifundefined{ifGPblacktext}{%
    \newif\ifGPblacktext
    \GPblacktexttrue
  }{}%
  \let\gplgaddtomacro\g@addto@macro
  \gdef\gplbacktext{}%
  \gdef\gplfronttext{}%
  \makeatother
  \ifGPblacktext
    \def\colorrgb#1{}%
    \def\colorgray#1{}%
  \else
    \ifGPcolor
      \def\colorrgb#1{\color[rgb]{#1}}%
      \def\colorgray#1{\color[gray]{#1}}%
      \expandafter\def\csname LTw\endcsname{\color{white}}%
      \expandafter\def\csname LTb\endcsname{\color{black}}%
      \expandafter\def\csname LTa\endcsname{\color{black}}%
      \expandafter\def\csname LT0\endcsname{\color[rgb]{1,0,0}}%
      \expandafter\def\csname LT1\endcsname{\color[rgb]{0,1,0}}%
      \expandafter\def\csname LT2\endcsname{\color[rgb]{0,0,1}}%
      \expandafter\def\csname LT3\endcsname{\color[rgb]{1,0,1}}%
      \expandafter\def\csname LT4\endcsname{\color[rgb]{0,1,1}}%
      \expandafter\def\csname LT5\endcsname{\color[rgb]{1,1,0}}%
      \expandafter\def\csname LT6\endcsname{\color[rgb]{0,0,0}}%
      \expandafter\def\csname LT7\endcsname{\color[rgb]{1,0.3,0}}%
      \expandafter\def\csname LT8\endcsname{\color[rgb]{0.5,0.5,0.5}}%
    \else
      \def\colorrgb#1{\color{black}}%
      \def\colorgray#1{\color[gray]{#1}}%
      \expandafter\def\csname LTw\endcsname{\color{white}}%
      \expandafter\def\csname LTb\endcsname{\color{black}}%
      \expandafter\def\csname LTa\endcsname{\color{black}}%
      \expandafter\def\csname LT0\endcsname{\color{black}}%
      \expandafter\def\csname LT1\endcsname{\color{black}}%
      \expandafter\def\csname LT2\endcsname{\color{black}}%
      \expandafter\def\csname LT3\endcsname{\color{black}}%
      \expandafter\def\csname LT4\endcsname{\color{black}}%
      \expandafter\def\csname LT5\endcsname{\color{black}}%
      \expandafter\def\csname LT6\endcsname{\color{black}}%
      \expandafter\def\csname LT7\endcsname{\color{black}}%
      \expandafter\def\csname LT8\endcsname{\color{black}}%
    \fi
  \fi
    \setlength{\unitlength}{0.0500bp}%
    \ifx\gptboxheight\undefined%
      \newlength{\gptboxheight}%
      \newlength{\gptboxwidth}%
      \newsavebox{\gptboxtext}%
    \fi%
    \setlength{\fboxrule}{0.5pt}%
    \setlength{\fboxsep}{1pt}%
\begin{picture}(7200.00,5040.00)%
    \gplgaddtomacro\gplbacktext{%
      \csname LTb\endcsname
      \put(726,704){\makebox(0,0)[r]{\strut{}$0$}}%
      \put(726,1116){\makebox(0,0)[r]{\strut{}$0.1$}}%
      \put(726,1527){\makebox(0,0)[r]{\strut{}$0.2$}}%
      \put(726,1939){\makebox(0,0)[r]{\strut{}$0.3$}}%
      \put(726,2350){\makebox(0,0)[r]{\strut{}$0.4$}}%
      \put(726,2762){\makebox(0,0)[r]{\strut{}$0.5$}}%
      \put(726,3173){\makebox(0,0)[r]{\strut{}$0.6$}}%
      \put(726,3585){\makebox(0,0)[r]{\strut{}$0.7$}}%
      \put(726,3996){\makebox(0,0)[r]{\strut{}$0.8$}}%
      \put(726,4408){\makebox(0,0)[r]{\strut{}$0.9$}}%
      \put(726,4819){\makebox(0,0)[r]{\strut{}$1$}}%
      \put(858,484){\makebox(0,0){\strut{}$0$}}%
      \put(2344,484){\makebox(0,0){\strut{}$0.5$}}%
      \put(3831,484){\makebox(0,0){\strut{}$1$}}%
      \put(5317,484){\makebox(0,0){\strut{}$1.5$}}%
      \put(6803,484){\makebox(0,0){\strut{}$2$}}%
    }%
    \gplgaddtomacro\gplfronttext{%
      \csname LTb\endcsname
      \put(198,2761){\makebox(0,0){\strut{}$\mu_0$}}%
      \put(3830,154){\makebox(0,0){\strut{}$\gamma$}}%
    }%
    \gplbacktext
    \put(0,0){\includegraphics{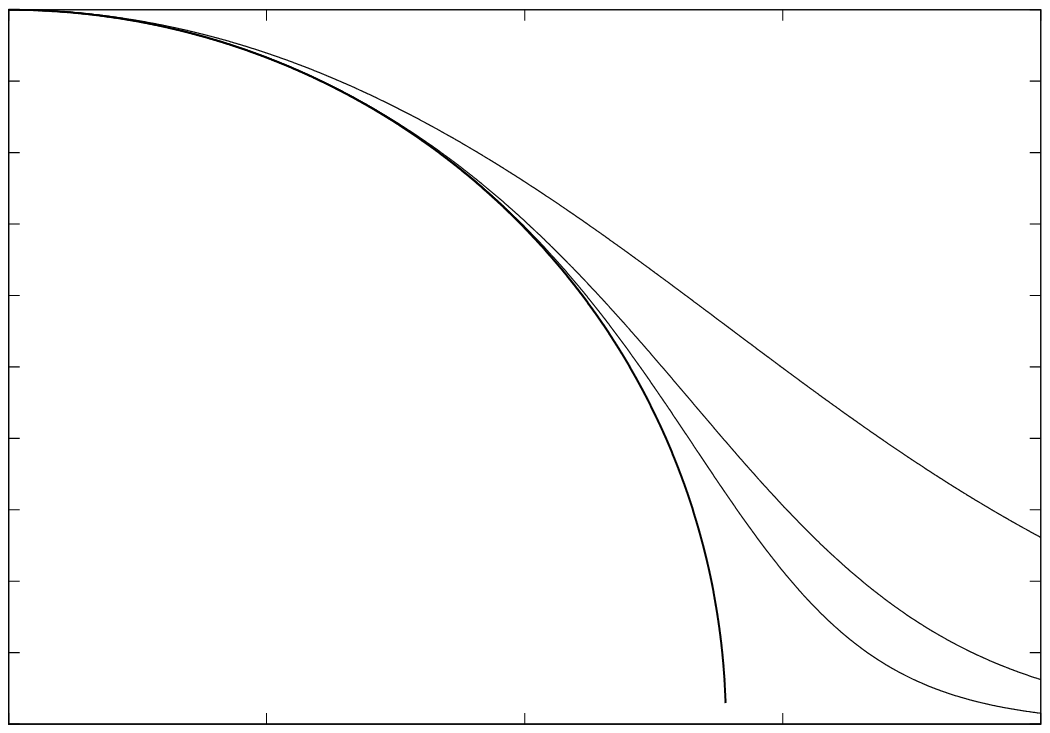}}%
    \gplfronttext
  \end{picture}%
\endgroup

%% file: Fig63.tex
\begingroup
  \makeatletter
  \providecommand\color[2][]{%
    \GenericError{(gnuplot) \space\space\space\@spaces}{%
      Package color not loaded in conjunction with
      terminal option `colourtext'%
    }{See the gnuplot documentation for explanation.%
    }{Either use 'blacktext' in gnuplot or load the package
      color.sty in LaTeX.}%
    \renewcommand\color[2][]{}%
  }%
  \providecommand\includegraphics[2][]{%
    \GenericError{(gnuplot) \space\space\space\@spaces}{%
      Package graphicx or graphics not loaded%
    }{See the gnuplot documentation for explanation.%
    }{The gnuplot epslatex terminal needs graphicx.sty or graphics.sty.}%
    \renewcommand\includegraphics[2][]{}%
  }%
  \providecommand\rotatebox[2]{#2}%
  \@ifundefined{ifGPcolor}{%
    \newif\ifGPcolor
    \GPcolorfalse
  }{}%
  \@ifundefined{ifGPblacktext}{%
    \newif\ifGPblacktext
    \GPblacktexttrue
  }{}%
  \let\gplgaddtomacro\g@addto@macro
  \gdef\gplbacktext{}%
  \gdef\gplfronttext{}%
  \makeatother
  \ifGPblacktext
    \def\colorrgb#1{}%
    \def\colorgray#1{}%
  \else
    \ifGPcolor
      \def\colorrgb#1{\color[rgb]{#1}}%
      \def\colorgray#1{\color[gray]{#1}}%
      \expandafter\def\csname LTw\endcsname{\color{white}}%
      \expandafter\def\csname LTb\endcsname{\color{black}}%
      \expandafter\def\csname LTa\endcsname{\color{black}}%
      \expandafter\def\csname LT0\endcsname{\color[rgb]{1,0,0}}%
      \expandafter\def\csname LT1\endcsname{\color[rgb]{0,1,0}}%
      \expandafter\def\csname LT2\endcsname{\color[rgb]{0,0,1}}%
      \expandafter\def\csname LT3\endcsname{\color[rgb]{1,0,1}}%
      \expandafter\def\csname LT4\endcsname{\color[rgb]{0,1,1}}%
      \expandafter\def\csname LT5\endcsname{\color[rgb]{1,1,0}}%
      \expandafter\def\csname LT6\endcsname{\color[rgb]{0,0,0}}%
      \expandafter\def\csname LT7\endcsname{\color[rgb]{1,0.3,0}}%
      \expandafter\def\csname LT8\endcsname{\color[rgb]{0.5,0.5,0.5}}%
    \else
      \def\colorrgb#1{\color{black}}%
      \def\colorgray#1{\color[gray]{#1}}%
      \expandafter\def\csname LTw\endcsname{\color{white}}%
      \expandafter\def\csname LTb\endcsname{\color{black}}%
      \expandafter\def\csname LTa\endcsname{\color{black}}%
      \expandafter\def\csname LT0\endcsname{\color{black}}%
      \expandafter\def\csname LT1\endcsname{\color{black}}%
      \expandafter\def\csname LT2\endcsname{\color{black}}%
      \expandafter\def\csname LT3\endcsname{\color{black}}%
      \expandafter\def\csname LT4\endcsname{\color{black}}%
      \expandafter\def\csname LT5\endcsname{\color{black}}%
      \expandafter\def\csname LT6\endcsname{\color{black}}%
      \expandafter\def\csname LT7\endcsname{\color{black}}%
      \expandafter\def\csname LT8\endcsname{\color{black}}%
    \fi
  \fi
    \setlength{\unitlength}{0.0500bp}%
    \ifx\gptboxheight\undefined%
      \newlength{\gptboxheight}%
      \newlength{\gptboxwidth}%
      \newsavebox{\gptboxtext}%
    \fi%
    \setlength{\fboxrule}{0.5pt}%
    \setlength{\fboxsep}{1pt}%
\begin{picture}(7200.00,5040.00)%
    \gplgaddtomacro\gplbacktext{%
      \csname LTb\endcsname
      \put(726,704){\makebox(0,0)[r]{\strut{}$0.2$}}%
      \put(726,1527){\makebox(0,0)[r]{\strut{}$0.4$}}%
      \put(726,2350){\makebox(0,0)[r]{\strut{}$0.6$}}%
      \put(726,3173){\makebox(0,0)[r]{\strut{}$0.8$}}%
      \put(726,3996){\makebox(0,0)[r]{\strut{}$1$}}%
      \put(726,4819){\makebox(0,0)[r]{\strut{}$1.2$}}%
      \put(858,484){\makebox(0,0){\strut{}$0$}}%
      \put(1601,484){\makebox(0,0){\strut{}$0.2$}}%
      \put(2344,484){\makebox(0,0){\strut{}$0.4$}}%
      \put(3087,484){\makebox(0,0){\strut{}$0.6$}}%
      \put(3831,484){\makebox(0,0){\strut{}$0.8$}}%
      \put(4574,484){\makebox(0,0){\strut{}$1$}}%
      \put(5317,484){\makebox(0,0){\strut{}$1.2$}}%
      \put(6060,484){\makebox(0,0){\strut{}$1.4$}}%
      \put(6803,484){\makebox(0,0){\strut{}$1.6$}}%
    }%
    \gplgaddtomacro\gplfronttext{%
      \csname LTb\endcsname
      \put(3830,154){\makebox(0,0){\strut{}$\gamma$}}%
    }%
    \gplbacktext
    \put(0,0){\includegraphics{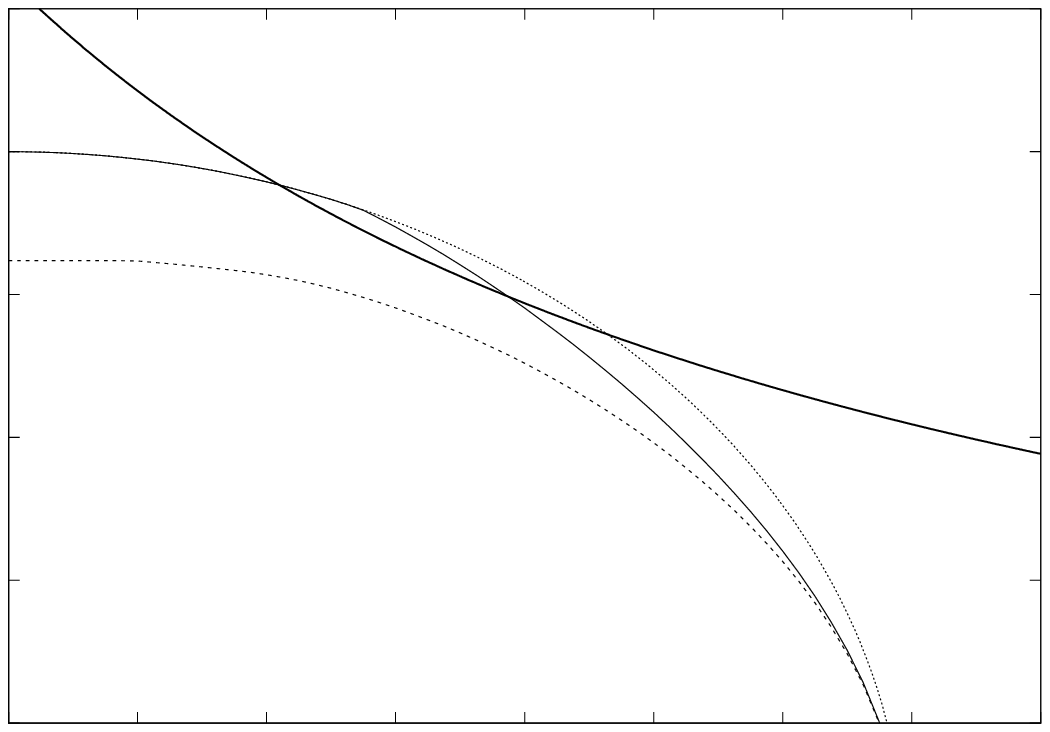}}%
    \gplfronttext
  \end{picture}%
\endgroup

%% file: FigE1.tex
\begingroup
  \makeatletter
  \providecommand\color[2][]{%
    \GenericError{(gnuplot) \space\space\space\@spaces}{%
      Package color not loaded in conjunction with
      terminal option `colourtext'%
    }{See the gnuplot documentation for explanation.%
    }{Either use 'blacktext' in gnuplot or load the package
      color.sty in LaTeX.}%
    \renewcommand\color[2][]{}%
  }%
  \providecommand\includegraphics[2][]{%
    \GenericError{(gnuplot) \space\space\space\@spaces}{%
      Package graphicx or graphics not loaded%
    }{See the gnuplot documentation for explanation.%
    }{The gnuplot epslatex terminal needs graphicx.sty or graphics.sty.}%
    \renewcommand\includegraphics[2][]{}%
  }%
  \providecommand\rotatebox[2]{#2}%
  \@ifundefined{ifGPcolor}{%
    \newif\ifGPcolor
    \GPcolorfalse
  }{}%
  \@ifundefined{ifGPblacktext}{%
    \newif\ifGPblacktext
    \GPblacktexttrue
  }{}%
  \let\gplgaddtomacro\g@addto@macro
  \gdef\gplbacktext{}%
  \gdef\gplfronttext{}%
  \makeatother
  \ifGPblacktext
    \def\colorrgb#1{}%
    \def\colorgray#1{}%
  \else
    \ifGPcolor
      \def\colorrgb#1{\color[rgb]{#1}}%
      \def\colorgray#1{\color[gray]{#1}}%
      \expandafter\def\csname LTw\endcsname{\color{white}}%
      \expandafter\def\csname LTb\endcsname{\color{black}}%
      \expandafter\def\csname LTa\endcsname{\color{black}}%
      \expandafter\def\csname LT0\endcsname{\color[rgb]{1,0,0}}%
      \expandafter\def\csname LT1\endcsname{\color[rgb]{0,1,0}}%
      \expandafter\def\csname LT2\endcsname{\color[rgb]{0,0,1}}%
      \expandafter\def\csname LT3\endcsname{\color[rgb]{1,0,1}}%
      \expandafter\def\csname LT4\endcsname{\color[rgb]{0,1,1}}%
      \expandafter\def\csname LT5\endcsname{\color[rgb]{1,1,0}}%
      \expandafter\def\csname LT6\endcsname{\color[rgb]{0,0,0}}%
      \expandafter\def\csname LT7\endcsname{\color[rgb]{1,0.3,0}}%
      \expandafter\def\csname LT8\endcsname{\color[rgb]{0.5,0.5,0.5}}%
    \else
      \def\colorrgb#1{\color{black}}%
      \def\colorgray#1{\color[gray]{#1}}%
      \expandafter\def\csname LTw\endcsname{\color{white}}%
      \expandafter\def\csname LTb\endcsname{\color{black}}%
      \expandafter\def\csname LTa\endcsname{\color{black}}%
      \expandafter\def\csname LT0\endcsname{\color{black}}%
      \expandafter\def\csname LT1\endcsname{\color{black}}%
      \expandafter\def\csname LT2\endcsname{\color{black}}%
      \expandafter\def\csname LT3\endcsname{\color{black}}%
      \expandafter\def\csname LT4\endcsname{\color{black}}%
      \expandafter\def\csname LT5\endcsname{\color{black}}%
      \expandafter\def\csname LT6\endcsname{\color{black}}%
      \expandafter\def\csname LT7\endcsname{\color{black}}%
      \expandafter\def\csname LT8\endcsname{\color{black}}%
    \fi
  \fi
    \setlength{\unitlength}{0.0500bp}%
    \ifx\gptboxheight\undefined%
      \newlength{\gptboxheight}%
      \newlength{\gptboxwidth}%
      \newsavebox{\gptboxtext}%
    \fi%
    \setlength{\fboxrule}{0.5pt}%
    \setlength{\fboxsep}{1pt}%
\begin{picture}(7200.00,5040.00)%
    \gplgaddtomacro\gplbacktext{%
      \csname LTb\endcsname
      \put(726,704){\makebox(0,0)[r]{\strut{}$0$}}%
      \put(726,1116){\makebox(0,0)[r]{\strut{}$0.1$}}%
      \put(726,1527){\makebox(0,0)[r]{\strut{}$0.2$}}%
      \put(726,1939){\makebox(0,0)[r]{\strut{}$0.3$}}%
      \put(726,2350){\makebox(0,0)[r]{\strut{}$0.4$}}%
      \put(726,2762){\makebox(0,0)[r]{\strut{}$0.5$}}%
      \put(726,3173){\makebox(0,0)[r]{\strut{}$0.6$}}%
      \put(726,3585){\makebox(0,0)[r]{\strut{}$0.7$}}%
      \put(726,3996){\makebox(0,0)[r]{\strut{}$0.8$}}%
      \put(726,4408){\makebox(0,0)[r]{\strut{}$0.9$}}%
      \put(726,4819){\makebox(0,0)[r]{\strut{}$1$}}%
      \put(858,484){\makebox(0,0){\strut{}$0$}}%
      \put(1453,484){\makebox(0,0){\strut{}$0.2$}}%
      \put(2047,484){\makebox(0,0){\strut{}$0.4$}}%
      \put(2642,484){\makebox(0,0){\strut{}$0.6$}}%
      \put(3236,484){\makebox(0,0){\strut{}$0.8$}}%
      \put(3831,484){\makebox(0,0){\strut{}$1$}}%
      \put(4425,484){\makebox(0,0){\strut{}$1.2$}}%
      \put(5020,484){\makebox(0,0){\strut{}$1.4$}}%
      \put(5614,484){\makebox(0,0){\strut{}$1.6$}}%
      \put(6208,484){\makebox(0,0){\strut{}$1.8$}}%
      \put(6803,484){\makebox(0,0){\strut{}$2$}}%
    }%
    \gplgaddtomacro\gplfronttext{%
      \csname LTb\endcsname
      \put(198,2761){\makebox(0,0){\strut{}$\mu_0$}}%
      \put(3830,154){\makebox(0,0){\strut{}$\gamma$}}%
    }%
    \gplbacktext
    \put(850,700){\includegraphics{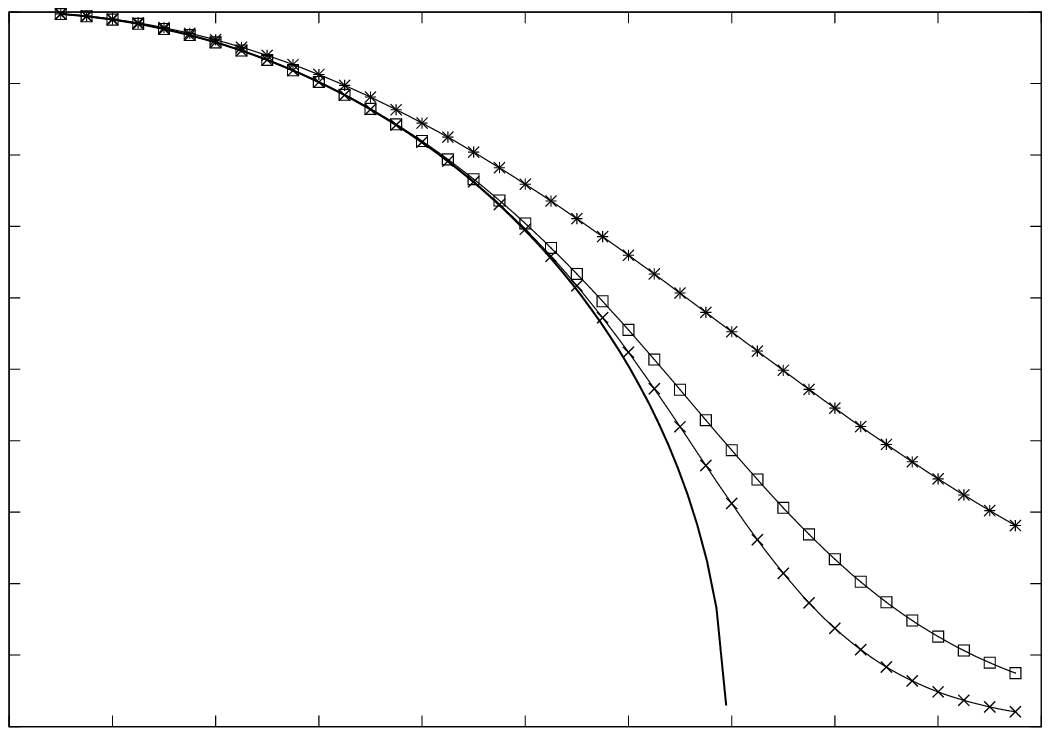}}%
    \gplfronttext
  \end{picture}%
\endgroup

%% file: FigE2.tex
\begingroup
  \makeatletter
  \providecommand\color[2][]{%
    \GenericError{(gnuplot) \space\space\space\@spaces}{%
      Package color not loaded in conjunction with
      terminal option `colourtext'%
    }{See the gnuplot documentation for explanation.%
    }{Either use 'blacktext' in gnuplot or load the package
      color.sty in LaTeX.}%
    \renewcommand\color[2][]{}%
  }%
  \providecommand\includegraphics[2][]{%
    \GenericError{(gnuplot) \space\space\space\@spaces}{%
      Package graphicx or graphics not loaded%
    }{See the gnuplot documentation for explanation.%
    }{The gnuplot epslatex terminal needs graphicx.sty or graphics.sty.}%
    \renewcommand\includegraphics[2][]{}%
  }%
  \providecommand\rotatebox[2]{#2}%
  \@ifundefined{ifGPcolor}{%
    \newif\ifGPcolor
    \GPcolorfalse
  }{}%
  \@ifundefined{ifGPblacktext}{%
    \newif\ifGPblacktext
    \GPblacktexttrue
  }{}%
  \let\gplgaddtomacro\g@addto@macro
  \gdef\gplbacktext{}%
  \gdef\gplfronttext{}%
  \makeatother
  \ifGPblacktext
    \def\colorrgb#1{}%
    \def\colorgray#1{}%
  \else
    \ifGPcolor
      \def\colorrgb#1{\color[rgb]{#1}}%
      \def\colorgray#1{\color[gray]{#1}}%
      \expandafter\def\csname LTw\endcsname{\color{white}}%
      \expandafter\def\csname LTb\endcsname{\color{black}}%
      \expandafter\def\csname LTa\endcsname{\color{black}}%
      \expandafter\def\csname LT0\endcsname{\color[rgb]{1,0,0}}%
      \expandafter\def\csname LT1\endcsname{\color[rgb]{0,1,0}}%
      \expandafter\def\csname LT2\endcsname{\color[rgb]{0,0,1}}%
      \expandafter\def\csname LT3\endcsname{\color[rgb]{1,0,1}}%
      \expandafter\def\csname LT4\endcsname{\color[rgb]{0,1,1}}%
      \expandafter\def\csname LT5\endcsname{\color[rgb]{1,1,0}}%
      \expandafter\def\csname LT6\endcsname{\color[rgb]{0,0,0}}%
      \expandafter\def\csname LT7\endcsname{\color[rgb]{1,0.3,0}}%
      \expandafter\def\csname LT8\endcsname{\color[rgb]{0.5,0.5,0.5}}%
    \else
      \def\colorrgb#1{\color{black}}%
      \def\colorgray#1{\color[gray]{#1}}%
      \expandafter\def\csname LTw\endcsname{\color{white}}%
      \expandafter\def\csname LTb\endcsname{\color{black}}%
      \expandafter\def\csname LTa\endcsname{\color{black}}%
      \expandafter\def\csname LT0\endcsname{\color{black}}%
      \expandafter\def\csname LT1\endcsname{\color{black}}%
      \expandafter\def\csname LT2\endcsname{\color{black}}%
      \expandafter\def\csname LT3\endcsname{\color{black}}%
      \expandafter\def\csname LT4\endcsname{\color{black}}%
      \expandafter\def\csname LT5\endcsname{\color{black}}%
      \expandafter\def\csname LT6\endcsname{\color{black}}%
      \expandafter\def\csname LT7\endcsname{\color{black}}%
      \expandafter\def\csname LT8\endcsname{\color{black}}%
    \fi
  \fi
    \setlength{\unitlength}{0.0500bp}%
    \ifx\gptboxheight\undefined%
      \newlength{\gptboxheight}%
      \newlength{\gptboxwidth}%
      \newsavebox{\gptboxtext}%
    \fi%
    \setlength{\fboxrule}{0.5pt}%
    \setlength{\fboxsep}{1pt}%
\begin{picture}(7200.00,5040.00)%
    \gplgaddtomacro\gplbacktext{%
      \csname LTb\endcsname
      \put(726,704){\makebox(0,0)[r]{\strut{}$0.6$}}%
      \put(726,1527){\makebox(0,0)[r]{\strut{}$0.7$}}%
      \put(726,2350){\makebox(0,0)[r]{\strut{}$0.8$}}%
      \put(726,3173){\makebox(0,0)[r]{\strut{}$0.9$}}%
      \put(726,3996){\makebox(0,0)[r]{\strut{}$1$}}%
      \put(726,4819){\makebox(0,0)[r]{\strut{}$1.1$}}%
      \put(858,484){\makebox(0,0){\strut{}$0.3$}}%
      \put(1601,484){\makebox(0,0){\strut{}$0.4$}}%
      \put(2344,484){\makebox(0,0){\strut{}$0.5$}}%
      \put(3087,484){\makebox(0,0){\strut{}$0.6$}}%
      \put(3830,484){\makebox(0,0){\strut{}$0.7$}}%
      \put(4574,484){\makebox(0,0){\strut{}$0.8$}}%
      \put(5317,484){\makebox(0,0){\strut{}$0.9$}}%
      \put(6060,484){\makebox(0,0){\strut{}$1$}}%
      \put(6803,484){\makebox(0,0){\strut{}$1.1$}}%
    }%
    \gplgaddtomacro\gplfronttext{%
      \csname LTb\endcsname
      \put(198,2761){\makebox(0,0){\strut{}$\mu_{min}$}}%
      \put(3830,154){\makebox(0,0){\strut{}$\gamma$}}%
    }%
    \gplbacktext
    \put(0,0){\includegraphics{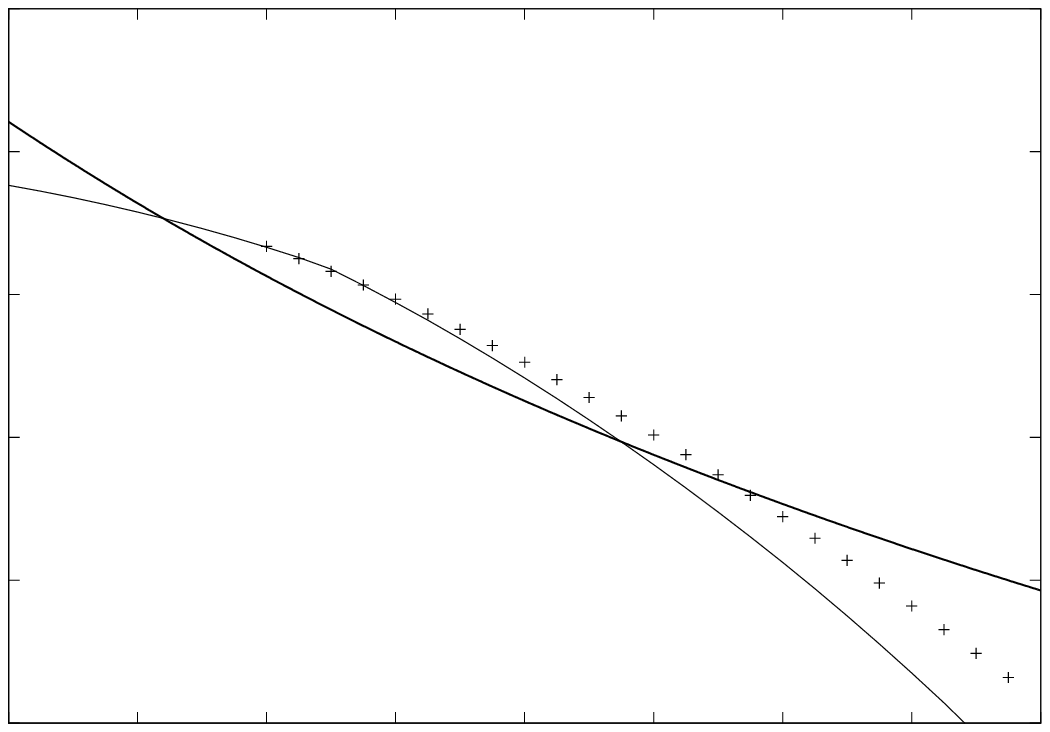}}%
    \gplfronttext
  \end{picture}%
\endgroup